\documentclass[%
 aip,
 reprint,%
]{revtex4-1}

\usepackage{graphicx}% Include figure files
\usepackage{dcolumn}% Align table columns on decimal point
\usepackage{bm}% bold math
\usepackage[utf8]{inputenc}
\usepackage[T1]{fontenc}
\usepackage{mathptmx}
\usepackage{etoolbox}

\usepackage{algorithm}
\usepackage{algorithmic}
\usepackage{bm}
\usepackage{enumitem}
\usepackage{xr}
\usepackage{lipsum}
\usepackage{hyperref}
\usepackage[mathcal]{euscript}
\usepackage{graphicx}
\usepackage{epstopdf}
\usepackage{tikz}
\usepackage[outline]{contour}
\usepackage{onimage}
\usepackage{amsfonts}
\usepackage{amsmath}
\usepackage{amssymb}
\usepackage{amsthm}
\usepackage{amsopn}
\usepackage{latexsym}
\usepackage{empheq}
\usepackage[caption=false]{subfig}
\ifpdf
  \DeclareGraphicsExtensions{.eps,.pdf,.png,.jpg}
\else
  \DeclareGraphicsExtensions{.eps}
\fi

\DeclareMathOperator{\E}{\ \!\mathbb{E}}

\DeclareMathOperator{\qf}{qf}
\DeclareMathOperator{\Tr}{Tr}

\DeclareMathOperator{\Range}{Range}

\DeclareMathOperator{\Null}{Null}

\DeclareMathOperator{\row}{row}
\DeclareMathOperator{\grad}{\nabla}
\DeclareMathOperator{\td}{\mathrm{d}}
\DeclareMathOperator{\ddt}{\frac{\td}{\td t}}

\DeclareMathOperator{\D}{d}

\newcommand{\R}{\mathbb{R}}

\newcommand{\mcal}[1]{{\mathcal{#1}}}
\newcommand{\mscr}[1]{{\mathcal{#1}}}

\makeatletter
\renewcommand*\env@matrix[1][*\c@MaxMatrixCols c]{%
  \hskip -\arraycolsep
  \let\@ifnextchar\new@ifnextchar
  \array{#1}}
\makeatother

\newcommand*\widefbox[1]{\fbox{\hspace{2em}#1\hspace{2em}}}

\newtheorem{theorem}{Theorem}
\newtheorem{lemma}{Lemma}
\newtheorem{proposition}{Proposition}
\newtheorem{corollary}{Corollary}

\theoremstyle{remark}
\newtheorem{example}{Example}
\newtheorem{remark}{Remark}

\makeatletter
\def\@email#1#2{%
 \endgroup
 \patchcmd{\titleblock@produce}
  {\frontmatter@RRAPformat}
  {\frontmatter@RRAPformat{\produce@RRAP{*#1\href{mailto:#2}{#2}}}\frontmatter@RRAPformat}
  {}{}
}%
\makeatother
\begin{document}

\preprint{AIP/123-QED}

\title[Autoencoders for Projection]{Learning Nonlinear Projections for Reduced-Order Modeling of Dynamical Systems using Constrained Autoencoders}
% Force line breaks with \\
\author{Samuel E. Otto}
\email{ottoncdr@uw.edu}
\affiliation{AI Institute in Dynamic Systems, University of Washington, Seattle, WA, USA}
\author{Gregory R. Macchio}%
% \email{gm0796@princeton.edu}%
\author{Clarence W. Rowley}
% \email{cwrowley@princeton.edu}
\affiliation{Mechanical and Aerospace Engineering, Princeton University, NJ, USA}

\date{\today}% It is always \today, today,

\begin{abstract}
Recently developed reduced-order modeling techniques aim to approximate nonlinear dynamical systems on low-dimensional manifolds learned from data.
This is an effective approach for modeling dynamics in a post-transient regime where the effects of initial conditions and other disturbances have decayed.
However, modeling transient dynamics near an underlying manifold, as needed for real-time control and forecasting applications, is complicated by the effects of fast dynamics and nonnormal sensitivity mechanisms.
To begin to address these issues, we introduce a parametric class of nonlinear projections described by constrained autoencoder neural networks in which both the manifold and the projection fibers are learned from data.
Our architecture uses invertible activation functions and biorthogonal weight matrices to ensure that the encoder is a left inverse of the decoder. 
We also introduce new dynamics-aware cost functions that promote learning of oblique projection fibers that account for fast dynamics and nonnormality.
To demonstrate these methods and the specific challenges they address, we provide a detailed case study of a three-state model of vortex shedding in the wake of a bluff body immersed in a fluid, which has a two-dimensional slow manifold that can be computed analytically.
In anticipation of future applications to high-dimensional systems, we also propose several techniques for constructing computationally efficient reduced-order models using our proposed nonlinear projection framework.
This includes a novel sparsity-promoting penalty for the encoder that avoids detrimental weight matrix shrinkage via computation on the Grassmann manifold.
\end{abstract}

\maketitle

\begin{quotation}
% The ``lead paragraph'' is encapsulated with the \LaTeX\ 
% \verb+quotation+ environment and is formatted as a single paragraph before the first section heading. 
% (The \verb+quotation+ environment reverts to its usual meaning after the first sectioning command.) 
% Note that numbered references are allowed in the lead paragraph.
% %
% The lead paragraph will only be found in an article being prepared for the journal \textit{Chaos}.
% According to the AIP author instructions (https://publishing.aip.org/resources/researchers/author-instructions/#prep), ``The first paragraph of the article should be the Lead Paragraph and contain the main points of the article, providing the “big picture” in a way that can be understood by non-specialist readers.''
Reduced-order modeling involves constructing a low-dimensional approximation of a high-dimensional dynamical system in order to enable tasks such as rapid forecasting, state estimation/tracking from streaming observations, and feedback control.
Autoencoders are a type of neural network that achieves dimensionality reduction by first compressing (encoding) and then reconstructing (decoding) high-dimensional state vectors.
We introduce a novel autoencoder architecture that can be used to project dynamical systems onto learned low-dimensional submanifolds of the state space.
Unlike prior work, we are able to learn appropriate projection fibers consisting of states that project to the same point on the manifold.
These fibers are crucial for obtaining accurate forecasts for states that do not lie on the manifold.
We introduce two new dynamics-aware cost functions to learn projections with appropriate fibers for reduced-order modeling of dynamical systems.
Finally, we compare our approach to standard architectures and cost functions on a simple slow-fast system.
\end{quotation}

\section{Introduction}
\label{sec:intro}
Dynamical systems arising from discretized continuum equations such as those governing fluid flows are often too high-dimensional to be used for real-time forecasting, state estimation, and control applications.
Simplified reduced-order models (ROMs) can be constructed by projecting the original dynamical system, referred to as the full-order model (FOM), into lower-dimensional spaces.
For reviews of existing methods, see \citet{Ghadami2022data, Rowley2017model, Benner2017model, Rozza2022advanced}.

The simplest approach is to employ linear projections onto subspaces.
By far the most widely used method of this type is Proper Orthogonal Decomposition (POD), which is also known as Principal Component Analysis (PCA).
As pointed out by \citet{Ohlberger2016reduced}, the effectiveness of linear projections for dimensionality reduction is limited by how closely relevant trajectories of the system can be approximated in linear subspaces.
This can be quantified using various measures such as Kolmogorov $n$-width and the decay rates of singular values obtained by POD.
Advection-dominated fluid flows exhibit spatially translating coherent structures that are notoriously difficult to model in low-dimensional subspaces.
This has motivated the development of techniques for projecting dynamics of fluid flows onto low-dimensional curved manifolds.
A recent approach by \citet{Lee2020model} with subsequent extensions in \citet{Romor2023nonlinear} projects dynamics orthogonally onto a manifold learned from data using a convolutional autoencoder neural network.
Similarly, \citet{Anderson2022evolution} project dynamics orthogonally onto a user-specified manifold with an interpretable parametrization.
Another approach utilized by \citet{Geelen2023operator} and \citet{Benner2022quadratic} is to project dynamics onto a manifold expressed as a graph over a POD subspace, in a direction orthogonal to that subspace.

A common feature of the above approaches for nonlinear projection-based model reduction is the use of orthogonal projection.
However, the ``direction of projection'' determined by the projection fibers is of critical importance for modeling transient dynamics.
To understand this, we first note that projection is unnecessary when modeling the dynamics of a system after transients have decayed onto an attracting submanifold.
Indeed, one can find an embedding of the underlying manifold from post-transient data and then learn the dynamics in the embedding space.
Essentially any embedding will do since in this case one only cares about the system's behavior on the manifold.
This is the principle behind successful data-driven methods for approximating dynamics near spectral submanifolds by \citet{Cenedese2022data} and other low-dimensional manifolds learned from data using auotencoders as in \citet{Fresca2021comprehensive, Conti2023reduced, Champion2019data}.
On the other hand, projection is needed to account for the ways in which perturbed trajectories settle back onto attracting manifolds.
This is critical for modeling the effect of actuation and control because input signals can result in such perturbations.

To understand the importance of the projection fibers from a geometric point of view, consider the basin of an attracting normally hyperbolic invariant manifold.
The basin is known to have an ``asymptotic rate foliation'' with leaves consisting of initial states that approach the same trajectory on the manifold\cite{Fenichel1974asymptotic,Fenichel1977asymptotic,Kuehn2015multiple}.
Specifically, this is the trajectory of the base point of intersection of the leaf with the invariant manifold.
% We can account for the transient effects of perturbations and ensure that trajectories of FOM settle onto corresponding trajectories of the ROM by projecting along the leaves of the asymptotic rate foliation.
Using a projection that collapses each leaf to its base point ensures that trajectories of the FOM settle onto the projected trajectories of the ROM at a rate determined by the fast time scale.
If the projection fibers were different, then there could be a persistent or growing error between the FOM and the ROM.
Near the manifold, the correct affine projections vary spatially according to a nonlinear partial differential equation derived by \citet{Roberts1989appropriate}.
While this equation is difficult to solve analytically, \citet{Roberts2000computer} uses computer algebra to find series expansions for spatially varying modes defining the projection near equilibria.
We illustrate the importance of the direction of projection using several examples including slow-fast systems with attracting slow manifolds \cite{Kuehn2015multiple}.

The direction of projection is also important for modeling nonnormal dynamical systems such as those arising from shear-dominated fluid flows \cite{Trefethen-93, Schmid2001stability}.
Linear dynamical systems governed by nonnormal operators can give rise to phenomena including transient growth and high sensitivity to state variables that remain small along trajectories \cite{Trefethen-93, Embree2005spectra}.
The term ``nonnormality'' has also been used to characterize nonlinear systems exhibiting these phenomena either due to nonnormal linearized dynamics, or other nonlinear effects such as asymmetric nonlinear coupling between states.
% The linearized dynamics of these systems are governed by nonnormal operators, giving rise to phenomena including transient growth and sensitivity to state variables that remain small along trajectories \cite{Trefethen-93, Embree2005spectra}.
Model reduction methods for linear systems \cite{Antoulas2005approximation, Gugercin2008H2, Rowley2005model} yield oblique projections that account for nonnormality.
To shed light on the utility of oblique projections, consider the oblique projections found using Balanced Truncation (BT) \cite{Moore1981principal}.
The BT projection coincides with state variable truncation in a coordinate system where the observability and controllability Gramians of a linear dynamical system are equal and diagonal.
In nonlinear systems exhibiting nonnormal dynamics, oblique linear projections have also proven to be useful for reduced-order modeling \cite{Benner2018H2, Benner2017balanced, Ahuja2010feedback, Barbagallo2009closed, Ilak2010model, Illingworth2011feedback, Otto2022optimizing, Otto2022model}.
For example, Covariance Balancing Reduction using Adjoint Snapshots (CoBRAS) \cite{Otto2022model} replaces the controllability Gramian in BT with a covariance matrix of states along nonlinear trajectories.
The observability Gramian is replaced by a gradient covariance matrix measuring the sensitivity of future outputs to state perturbations. 
The nonlinear balancing method introduced by Scherpen \cite{Scherpen1993balancing} yields a nonlinear oblique projection.
This projection is constructed by truncating nonlinear coordinates in which functions measuring nonlinear observability and controllability are balanced in the neighborhood of a fixed point.
This particular projection is computationally expensive to compute for high-dimensional systems, though significant progress on this issue has been made by \citet{Kramer2022nlBT1, Kramer2023nlBT2} using local series expansion methods.
In each of these cases, an oblique projection is needed in order to balance competing requirements to capture controllability and observabilty or state variance and sensitivity in nonnormal nonlinear systems.

In order to construct accurate reduced-order models of the systems described above, we introduce a large parametric class of nonlinear oblique projectons defined by constrained autoencoder neural networks.
% In this paper we introduce a large parametric class of nonlinear oblique projections defined by constrained autoencoder neural networks.
Autoencoders consist of an encoder neural network that reduces the dimension of an input vector followed by a decoder neural network that aims to reconstruct the original vector \cite{Goodfellow2016deep}.
While the decoder can be trained to parametrize a manifold, the encoder does not generally recover the correct coordinates.
This means that autoencoders do not generally define projections.
Because of this issue, \citet{Lee2020model} neglect the encoder after training and project the dynamics orthogonally onto the manifold defined by the decoder.
A main contribution of our work is to introduce constraints on the architecture of an autoencoder so that it defines a nonlinear oblique projection.
Specifically, we ensure that the process of decoding followed by encoding is always the identity.
To do this, we introduce a pair of smooth activation functions which are inverses.
One is used in the encoder and the other is used in the decoder.
We also enforce bi-orthogonality constraints between the weight matrices defining corresponding layers of the encoder and decoder.
Related architectures include neural networks with orthogonality constraints as employed by \citet{Lezcano2019cheap} and the invertible neural networks developed by \citet{Dinh2014nice, Dinh2016density, Kingma2018glow} with applications to inverse problems by \citet{Ardizzone2018analyzing}.
Using our approach, we are able to utilize the encoder and its tangent map to construct nonlinear projection-based reduced-order models with learned projection fibers.
Specifically, the projection fibers can now be oblique and vary over the learned manifold.

The standard loss function used to train autoencoders minimizes the distance between data vectors and their reconstructions after applying the encoder and decoder.
Minimizing this loss encourages the encoder to learn a direction of projection that is orthogonal to the learned manifold.
In order to learn oblique projections for constructing accurate ROMs, we introduce two new loss functions leveraging trajectory data from the FOM and its governing equations.
The first loss function we introduce combines the usual reconstruction error with the error between the time derivative of trajectories projected onto the learned manifold and the time derivative of the ROM at the projected points.
This promotes learning of a manifold that lies near the training data and a direction of projection that yields correct time derivatives for the reduced-order model.
The second loss function is closely related to the gradient-weighted objective minimized by CoBRAS \cite{Otto2022model}.
Specifically, we weight the differences between data vectors and their reconstructions using the autoencoder against gradients of random projections of the FOM's output along trajectories.
This allows the network to learn the directions along which the state data can be safely projected onto the learned manifold while having minimal effect on future outputs of the system.
We also introduce a sparsity-promoting penalty for the weight matrices of the encoder.
In a similar manner to the Discrete Empirical Interpolation Method (DEIM) \cite{Chaturantabut2010nonlinear}, sparsifying the encoder provides computational speedups for ROMs of systems with sparse coupling between state variables.
In order to avoid the shrinkage and other biases concomitant with the standard $\ell^1$ penalization, our penalty is invariant under the action of invertible matrices applied from the left to the weight matrix.

The remainder of the paper is organized as follows.  In Section~\ref{sec:nlprom} we discuss how projection-based reduced-order models are constructed, and provide an example illustrating the importance of the fiber.  In Section~\ref{sec:architecture} we describe the architecture of our autoencoder, including the invertible activation functions and bi-orthogonality constraints that ensure that our autoencoder is a projection, and in Section~\ref{sec:cost_funs} we discuss the loss functions mentioned above.  We provide a detailed study of a simple model problem with three states, and a two-dimensional slow manifold in Section~\ref{sec:noack_case_study}.  Finally, in Section~\ref{sec:ROM_assembly} we discuss the construction of computationally efficient models.

\section{Nonlinear projection-based reduced-order modeling}
\label{sec:nlprom}

\begin{figure}
    \centering
    \scalebox{0.6}{
    \begin{tikzpicture}[scale=0.95,>=stealth]
        
        % TikZ Coding Notes: 
        %   1. Defining things in TikZ is like stacking paper on a desk
        %   2. "\uncover<#->" controls when things in the figure are displayed.
        
        % Defining Fill Style For Planes
        \tikzset{facestyle/.style={fill=black!15,opacity=1}}
        \tikzset{facestyle2/.style={fill=black!30,opacity=1}}
        
        % Defining Line Widths
        \tikzset{guide/.style={thin}}
        
        % Defining Angle Of W-Plane
        \def\th{25}
        
        % Drawing Ambient Space
        \draw (-6,-3.5) rectangle (6,3.75);
        \draw (-6,3.75) node[below right] {$\mcal{U}\subset\mcal{N}$};
        
        % Drawing Low Half Of P^{-1}(x)
        \draw (0,0) .. controls (1.5, -3) and (2.5, -3.5) .. (4,-3.25);
        
        % Drawing Manifold
        \def\dyy{-0.25}
        \def\dyyy{1.75}
        \filldraw[white] (1.5, -1.5) circle (26.75pt);
        \filldraw[white] (0.5,-1) circle (31pt);
        \draw (-5.5,0+\dyy+\dyyy) .. controls (-1-0.5,2+\dyy+\dyyy) and (1+0.5,-0.5+\dyy+\dyyy) .. (5.5,2+\dyy+\dyyy);
        \draw (-5.5,-3+\dyy) .. controls (-1-0.5,-1+\dyy) and (1+0.5,-3.5+\dyy) .. (5.5,-1+\dyy);
        \draw (-5.5,0+\dyy+\dyyy) -- (-5.5,-3+\dyy);
        \draw (5.5,2+\dyy+\dyyy)-- (5.5,-1+\dyy);
        \draw (-5.5,0+\dyy+\dyyy) node[below right] {$\hat{\mathcal{M}} = \Range(P)$};
        
        % V-Plane
        \draw[facestyle] (-2,0,-1) -- (-2,0,2) -- (2,0,2) -- (2,0,-1) -- cycle;
        \draw (1.5,-0.70) node[right] {$T_x \hat{\mathcal{M}} = \Range(\D P(x))$};
        
        % Upper Portion Of Tilted Plane
        \draw[facestyle2] ({2*cos(\th)},{2*sin(\th)},-1) -- ({2*cos(\th)},{2*sin(\th)},2) -- (0,0,2) -- (0,0,-1) -- cycle;
        \draw ({2*cos(\th)},{2*sin(\th)},-1) node[below right] {$\Range(\D  P(x)^T)$};
        
        % Upper Half Of Normals To planes
        \draw[thick] (0,0,0) -- ({-2*sin(\th)},{2*cos(\th)},0);
        \draw({-1.370*sin(\th)},{1.370*cos(\th)},0) node[above right] {$\Null(\D P(x)) = T_x P^{-1}(x) \: $};
        
        % Drawing Upper Half Of P^{-1}(x)
        \draw (-4,3.25) .. controls (-2.5, 3.5) and (-1.5, 3) .. (0, 0);
        \draw (4,-3.25) node[above right] {$P^{-1}(x)$};
        
        % Defining The (x,y,z) Of X Vector
        \def\x{-1.0} \def\y{1} \def\z{1} 
        \draw[->,very thick] (0,0,0) -- (\x,\y,\z);
        \draw (\x,\y,\z) node[above] {$v$};
        
        % Drawing The Oblique Projection Onto W-Space
        \draw[->,very thick] (0,0,0) -- ({\x + \y*tan(\th)},0,\z);
        \draw ({\x + \y*tan(\th)},0,\z) node[left] {$\D P(x)v$};
        \draw[guide] (\x,\y,\z) -- ({\x + \y*tan(\th)},0,\z);
        
        % Drawing x
        \filldraw (0,0) circle (1pt) node[right] {$x$};
      
    \end{tikzpicture}
    }
    \caption{The anatomy of a smooth projection as characterized by Thm.~1.15 in \citet{Michor2008topics}.}
    \label{fig:smooth_projection}
\end{figure}
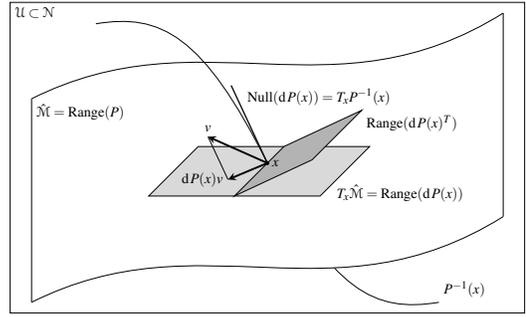

We consider a full-order model (FOM) described by a dynamical system
\begin{equation}
\begin{aligned}
    \ddt x &= f(x, u) \qquad x(0) = x_0 \\
    y &= g(x),
\end{aligned}
\label{eqn:FOM}
\end{equation}
with state variable $x(t)\in \R^n$, output observations $y(t)\in \R^m$, and inputs $u(t)$ taking values in an arbitrary space.
In many systems of interest the dynamics of the FOM can be accurately described on a low-dimensional submanifold $\mcal{M}\subset \R^n$ of the state space.
This can happen when the dynamics cause states to rapidly approach $\mcal{M}$, or when the system's output is insensitive to state variables normal to $\mcal{M}$ in some coordinate system.
Our goal is to identify a suitable manifold and construct a reduced-order model (ROM) of the form
\begin{equation}
\begin{aligned}
    \ddt \hat{x} &= \hat{f}(\hat{x}, u) \qquad \hat{x}(0) = \hat{x}_0 \\
    \hat{y} &= \hat{g}(\hat{x}),
\end{aligned}
\label{eqn:ROM}
\end{equation}
whose state $\hat{x}$ evolves on $\mcal{M}$ and whose output $\hat{y}$ approximates the output $y$ of the FOM over some set of inputs and initial conditions of interest.
One approach is to construct a smooth projection $P:\R^n \to \R^n$, that is an idempotent map $P \circ P = P$, and apply its tangent map to the FOM, yielding
\begin{equation}
\begin{aligned}
    \ddt \hat{x} &= \hat{f}_P(\hat{x},u) := \D P(\hat{x}) f(\hat{x}, u) \qquad \hat{x}(0) = P(x_0) \\
    \hat{y} &= g(\hat{x}).
\end{aligned}
\label{eqn:P_ROM}
\end{equation}
Geometrically, Theorem.~1.15 in \citet{Michor2008topics} (see Figure~\ref{fig:smooth_projection}) says that if $P$ is a smooth idempotent map on a connected manifold~$\mcal{N}$ (in our case $\mcal{N} = \R^n$) then the image set $\hat{\mcal{M}} = \Range(P)$ is automatically a smooth, closed, and connected submanifold of $\mcal{N}$.
\emph{We aim to find a projection whose image manifold accurately captures trajectories of interest from the FOM over a range of initial conditions and input signals.}
Moreover, the theorem shows that there is an open neighborhood $\mcal{U}$ of $\hat{\mcal{M}}$ in $\mcal{N}$ on which the tangent map $\D P(x)$ has constant rank equal to the dimension of $\hat{\mcal{M}}$.
In any such neighborhood $\mcal{U}$ of $\hat{\mcal{M}}$ where $\D P(x)$ has constant rank, the fiber $\left. P\right\vert_{\mcal{U}}^{-1}(x) = \{ p\in\mcal{U} \ : \ P(p) = x \}$ of each $x\in\hat{\mcal{M}}$ is a closed submanifold of $\mcal{U}$ with dimension complementary to $\hat{\mcal{M}}$ in $\mcal{N}$ and intersecting $\hat{\mcal{M}}$ transversally at $x$.
The tangent map $\D P (x)$ at $x\in\hat{\mcal{M}}$ is the linear projection on $T_x\mcal{N}$ whose range is $T_x\hat{\mcal{M}}$ and whose nullspace is tangent to the fiber, that is $\Null \D P(x) = T_x P^{-1}(x)$.
This characterization of smooth projections is depicted in Figure~\ref{fig:smooth_projection}.
Since the ROM in (\ref{eqn:P_ROM}) is obtained by modifying the FOM along the fibers of the projection,
\emph{we aim to design the projection so that varying initial states along the fibers has little affect on the system's output signal over a desired prediction horizon.}

One approach described by \citet{Lee2020model} is to parametrize a smooth submanifold $\hat{\mcal{M}} \subset \R^n$ and to define a projection by mapping $x$ to the nearest point on $\hat{\mcal{M}}$.
Such a projection is well-defined, smooth, and has constant rank in a neighborhood of $\hat{\mcal{M}}$ in $\R^n$ thanks to the tubular neighborhood theorem (Theorem~6.24 in \citet{Lee2013introduction}).
The fibers of this projection are orthogonal to $\hat{\mcal{M}}$ and the corresponding tangent map $\D P (x)$ is the orthogonal projection onto $T_x\hat{\mcal{M}}$.
While projecting the dynamics of the FOM orthogonally onto the tangent space of the learned manifold minimizes the projection error $\left\Vert f(x,u) - \D P(x) f(x,u) \right\Vert_x$ at each $x\in\hat{\mcal{M}}$, it can lead to large errors in the dynamics of the ROM described by (\ref{eqn:P_ROM}).
Therefore, we argue that the direction of projection as determined by the fibers of $P$ and their tangent spaces at intersections with $\Range(P)$ are important ingredients for constructing accurate nonlinear projection-based reduced-order models via (\ref{eqn:P_ROM}).

The following toy example illustrates why the projection fibers are important for modeling the dynamics of slow-fast systems (see \citet{Kuehn2015multiple}) using data-driven methods.
\begin{example}[Sources of projection error in a slow-fast system]
    \label{ex:toy_slow_fast}
    Consider the two dimensional system,
    \begin{equation}
        \begin{aligned}
            \dot{x}_1 &= \lambda x_1 (1 - x_1^2) \\
            \varepsilon \dot{x}_2 &= x_1^2 -x_2,
            \label{eq:pitchfork}
        \end{aligned}
    \end{equation}
    where $\lambda, \varepsilon >0$ and $\varepsilon^{-1} \gg \lambda$.
    There are two asymptotically stable fixed points at $(\pm 1, 1)$ and one unstable fixed point at $(0,0)$. 
    For small $\varepsilon$, (\ref{eq:pitchfork}) has an attracting slow invariant manifold containing the fixed points and lying near the critical manifold $x_2 = x_1^2$.
    Using Theorem~11.1.1 in \citet{Kuehn2015multiple}, we can express the slow manifold as a graph $x_2 = h_\epsilon(x_1)$ whose expansion in $\varepsilon$ is given by
    \begin{multline}
        % \begin{aligned}
        h_\epsilon(x_1) = x_1^2 + 2\lambda \left(x_1^4 - x_1^2 \right) \varepsilon \\
        + 4 \lambda^2 \left( 2 x_1^6 - 3 x_1^4 + x_1^2 \right) \varepsilon^2 
        + \mathcal{O}\left(\varepsilon^3\right).
        \label{eq:pitchfork_slow_manifold}
        % \end{aligned}
    \end{multline}
    Here, we use the parameter values $\lambda = 0.1$ and $\varepsilon = 0.1$.

    In Figure~\ref{fig:ic_proj_prob} we consider an initial condition (blue $+$) not lying on the slow manifold and two initial conditions (red $+$) resulting from different projections onto the slow manifold.
    The fast dynamics of $x_2$ cause the resulting trajectory to approach the slow manifold vertically.
    In the left panel, the trajectory of the orthogonally projected initial condition has a large phase error on the slow manifold, with the two trajectories only approaching each other at the slow rate $\mcal{O}(e^{-2\lambda t})$ as $t\to \infty$, as shown in Figure~\ref{fig:ic_proj_rates}.
    On the other hand, the trajectory of the vertically projected initial condition has zero phase error, with the two trajectories converging at the fast rate $\mcal{O}(e^{-t/\varepsilon})$. 
    % The aforementioned asymptotic analysis is confirmed in Figure \ref{fig:ic_proj_rates}.
    \begin{figure}
        \centering
        \includegraphics[scale=0.25]{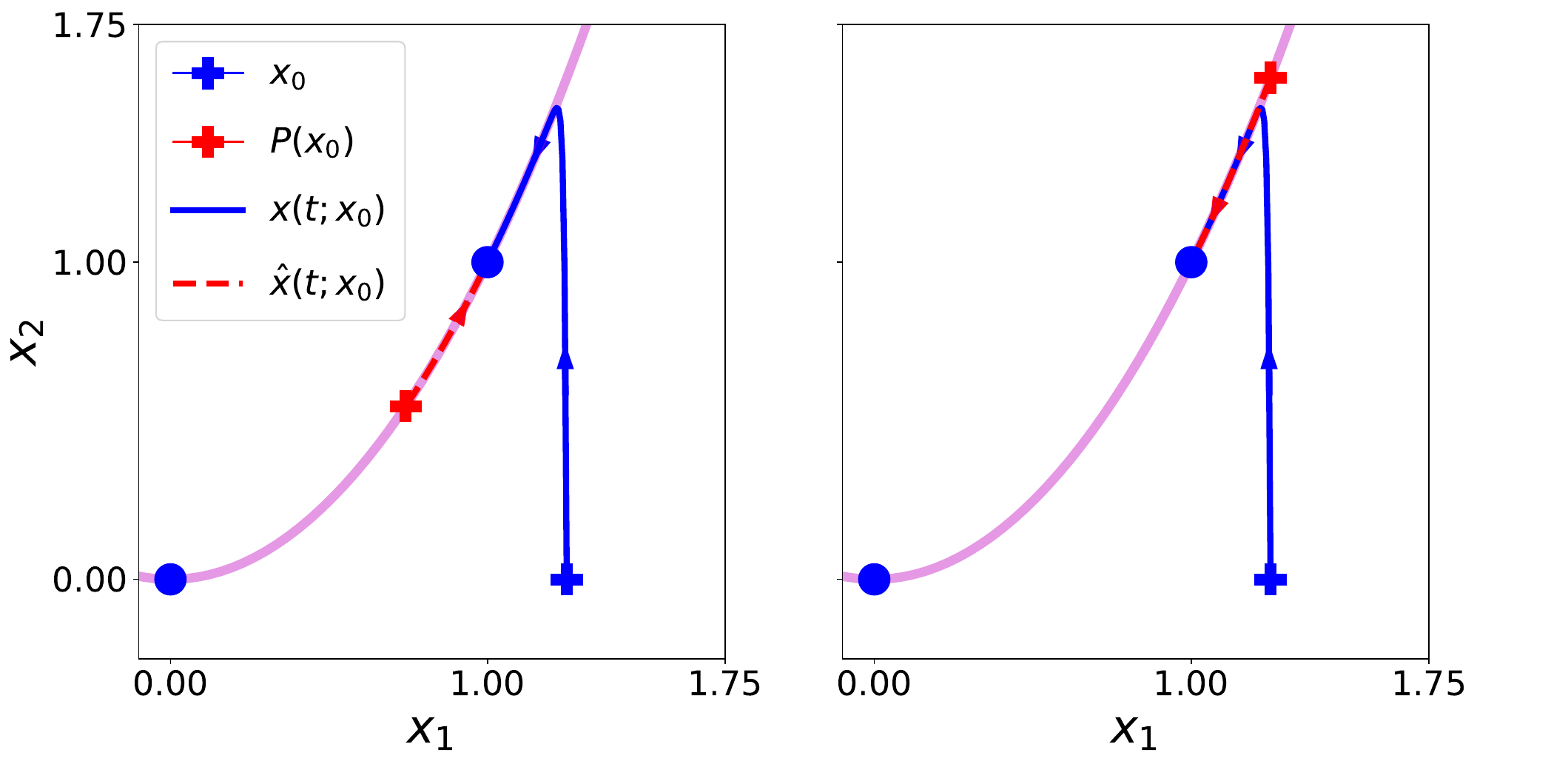}
        \caption{Projected dynamics for Example~1, comparing orthogonal projection onto the slow manifold (left) with projection along the direction of the fast dynamics (right).  The initial condition is shown as a blue plus, while the projected initial condition is shown as a red plus.  Equilibrium points are indicated by blue dots.}
        \label{fig:ic_proj_prob}
    \end{figure}
    \begin{figure}
        \centering
        \includegraphics[scale=0.21]{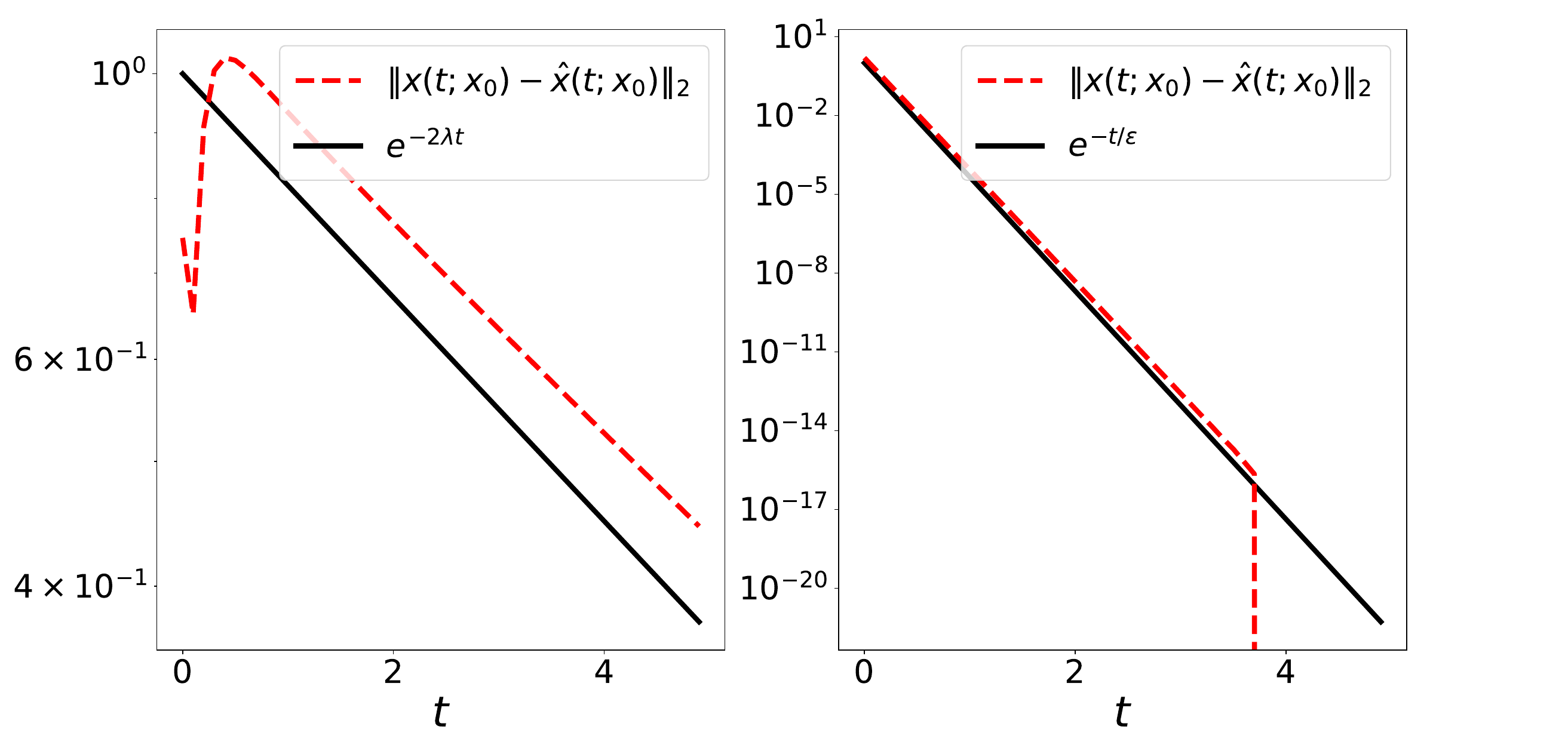}
        \caption{Absolute trajectory error for projection onto the slow manifold.
        \emph{Left:} Absolute error for the orthogonally projected initial condition is drawn in red and the expected asymptotic behavior in black.
        \emph{Right:} Absolute error and expected asymptotic behavior for the obliquely projected initial condition. Note the vastly different vertical scales.}
        \label{fig:ic_proj_rates}
    \end{figure}

    In Figure~\ref{fig:dyn_proj_prob} we consider two methods of projecting (\ref{eq:pitchfork_slow_manifold}) onto the tangent space of an approximate manifold lying near the true slow manifold.
    This mimics the typical situation when a manifold is learned from data.
    The vector field in (\ref{eq:pitchfork_slow_manifold}) evaluated along the approximate manifold (black arrows) has a large vertical component due to the approximation error and fast dynamics.
    Orthogonally projecting this vector field onto the approximate manifold in the left panel of Figure~\ref{fig:dyn_proj_prob} yields dynamics (red arrows) that incorrectly capture the dynamics on the nearby slow manifold (blue arrows).
    Even the stability types of the fixed points on the approximate manifold are the opposites of their counterparts in the true system.
    On the other hand, obliquely projecting the vector field onto the approximate manifold along vertical fibers cancels out the large contribution of the fast dynamics as shown in the right panel of Figure~\ref{fig:dyn_proj_prob}.
    The resulting projected system closely approximates the dynamics on the slow manifold and correctly captures the stability types of the fixed points.
    
    \begin{figure}
        \centering
        \includegraphics[scale=0.25]{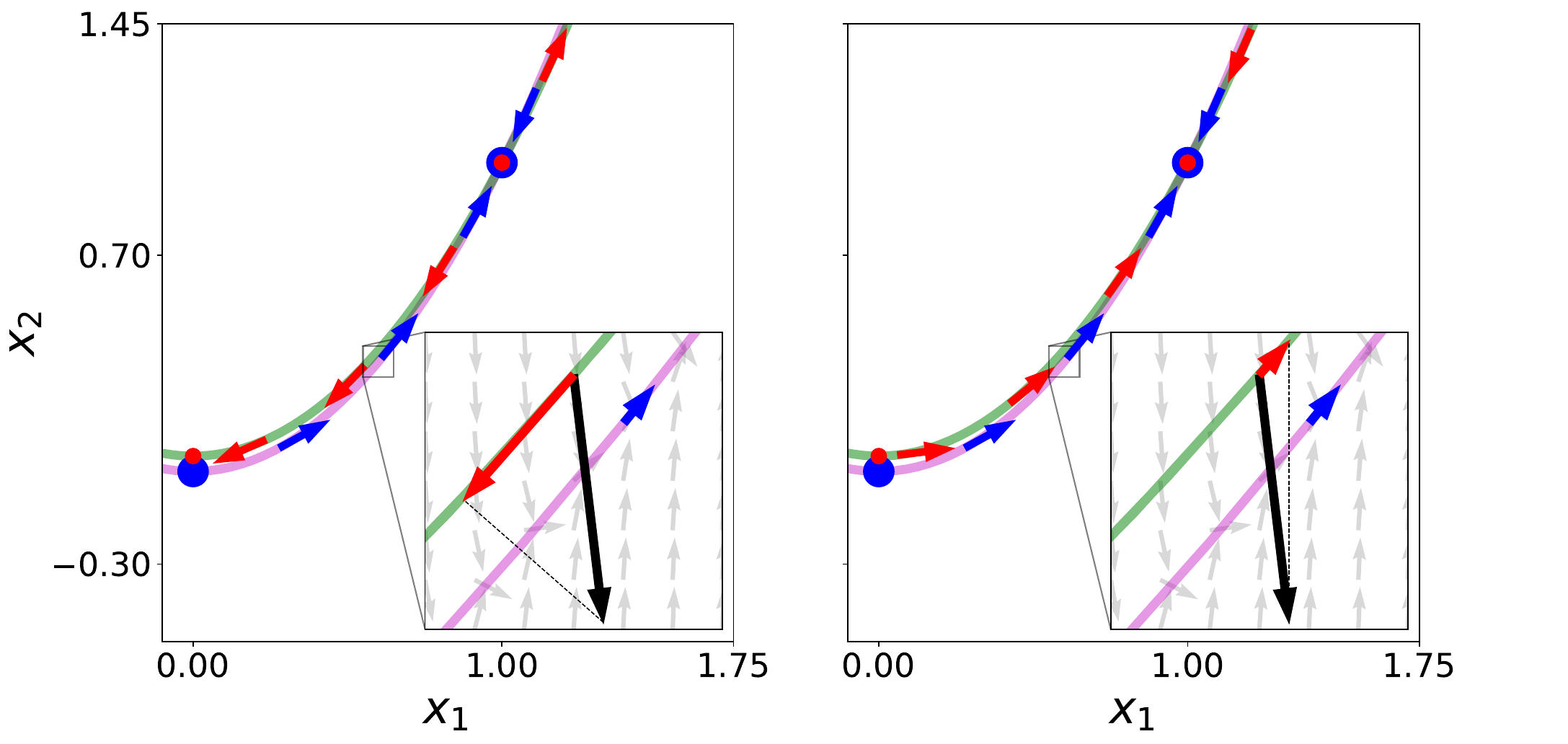}
        \caption{Projection onto an approximation of the slow manifold, given by $(x_1, 0.95x_1^2 + 0.05)$, shown in green, comparing orthogonal projection (left) with oblique projection (right).  The true slow manifold is shown in violet.  The direction of dynamics on the true slow manifold are shown as blue arrows, and the projected dynamics are shown as red arrows.
        Insets show the direction of dynamics of the full model in black, along with the corresponding projections. Note that projecting orthogonally reverses the stability types of the fixed points, even though the manifolds are so close.}
        \label{fig:dyn_proj_prob}
    \end{figure}
    
\end{example}

The importance of learning the correct direction of projection, which may be oblique to the learned manifold, motivates the development of a large parametric class of nonlinear projections based on autoencoders in the next section.
The choice of optimization objectives for training these autoencoders is also crucial and will be pursued in Section~\ref{sec:cost_funs}.

\section{Autoencoder architecture}
\label{sec:architecture}

\begin{figure}
    \centering
    \includegraphics[scale=0.7]{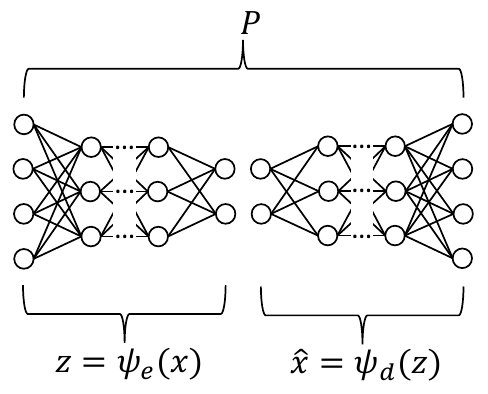}
    \caption{The architecture of an autoencoder, consisting of two component neural networks, the encoder $\psi_e$ and the decoder~$\psi_d$.}
    \label{fig:autoencoder}
\end{figure}

An autoencoder (in particular, an ``undercomplete'' autoencoder) is a neural network architecture depicted in Figure~\ref{fig:autoencoder} commonly used for dimension reduction and feature extraction in machine learning \cite{Goodfellow2016deep}.
It consists of an ``encoder'' $\psi_e$, which maps a data vector $x\in\R^n$ into a lower-dimensional representation or ``latent state'' $z\in \R^r$, $r < n$, and a ``decoder'' $\psi_d$ which reconstructs an approximation of $x$ from the extracted latent variables.
By optimizing the weights defining the encoder and decoder to accurately reconstruct data from a given distribution, the encoder learns a reduced set of features that describe the data.
If the encoder and decoder are smooth maps and the process of decoding and encoding through $\psi_e\circ \psi_d$ is the identity on the latent space, then, per our discussion in discussion in Section~\ref{sec:nlprom} the autoencoder $P = \psi_d \circ \psi_e$ is a smooth projection onto its range $\hat{\mcal{M}} = \Range(P) = \Range(\psi_d)$, which is a smooth manifold.
Moreover, the direction of projection is determined by the preimage fibers of the encoder $P^{-1}(\psi_d(z)) = \psi_e^{-1}(z)$.
In the context of model reduction, we can describe the dynamics of the projection-based reduced-order model (\ref{eqn:P_ROM}) with state $\hat{x} = \psi_d(z)$ in the latent space according to
% \begin{widetext}
\begin{empheq}[box=\fbox]{equation}
\begin{aligned}
    \ddt z &= \tilde{f}(z,u) := \D \psi_e(\psi_d(z)) f( \psi_d(z), u), \quad z(0) = \psi_e(x_0) \\
    y &= \tilde{g}(z) := g(\psi_d(z)).
\end{aligned}
\label{eqn:ROM_in_latent_space}
\end{empheq}
% \end{widetext}

In this setup, we can take advantage of the features learned by the encoder to define the crucial direction of projection for reduced-order modeling.
However, the constraint
\begin{equation}
    \psi_e\circ \psi_d = Id
    \label{eqn:decoder_encoder_constraint}
\end{equation}
has yet to be enforced in the design of autoencoders.
Instead, recent projection-based reduced-order modeling methods using autoencoders have followed the approach of \citet{Lee2020model}, in which the encoder is discarded and the dynamics are projected orthogonally onto the image manifold parametrized by the decoder.

\begin{figure*}
    \centering
    \includegraphics[scale=0.7]{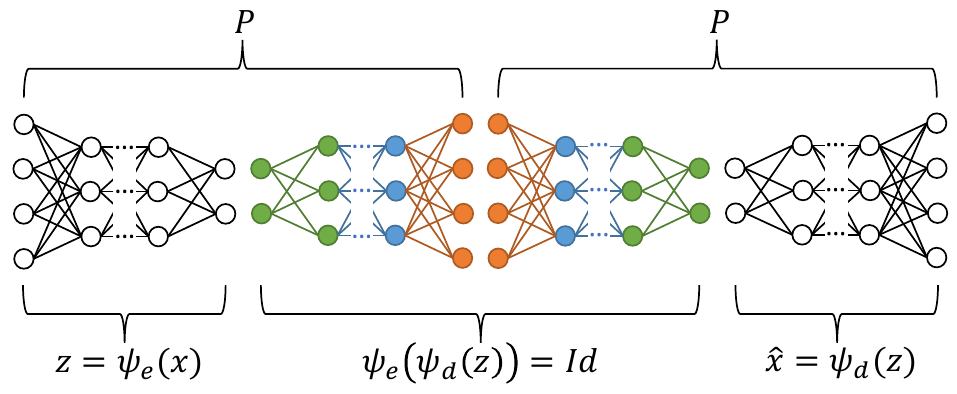}
    \caption{An autoencoder defines an idempotent map, i.e., a projection, as long as the process of decoding and then encoding any latent state $z$ is the identity.
    This constraint can be imposed layer-wise with the last layer of the decoder being ``undone'' by the first layer of the encoder (in orange), the second to last layer of decoder being undone by the second layer of the encoder (in blue), and so on.
    The corresponding layers of the decoder and encoder form a collapsing ``telescope'' that produces the identity.}
    \label{fig:decoder_encoder_composition}
\end{figure*}

Here, we design an auotencoder architecture in which the constraint (\ref{eqn:decoder_encoder_constraint}) is automatically satisfied.
This is accomplished layer-wise, as illustrated in Figure~\ref{fig:decoder_encoder_composition} by defining the encoder and decoder as compositions of layers
\begin{equation}
    \psi_{e} = \psi_{e}^{(1)} \circ \cdots \circ \psi_{e}^{(L)}, \qquad
    \psi_{d} = \psi_{d}^{(L)} \circ \cdots \circ \psi_{d}^{(1)},
\end{equation}
with the property that $\psi_e^{(l)} \circ \psi_d^{(l)} = Id$ for each layer $l$.
This ensures that the composition telescopes to produce the identitiy, that is,
\begin{equation}
\begin{aligned}
    \psi_{e} \circ \psi_{d} &= \psi_{e}^{(1)} \circ \cdots \circ \psi_{e}^{(L-1)} \circ \psi_{e}^{(L)} \circ \psi_{d}^{(L)} \circ \psi_{d}^{(L-1)} \circ \cdots \circ \psi_{d}^{(1)} \\
    & = \psi_{e}^{(1)} \circ \cdots \circ \psi_{e}^{(L-1)} \circ \psi_{d}^{(L-1)} \circ \cdots \circ \psi_{d}^{(1)} \\
    & \vdots \\
    &= \psi_{e}^{(1)} \circ \psi_{d}^{(1)} = Id.
\end{aligned}
\end{equation}
We note that if $\psi_{e}^{(l)}:\R^{n_l}\to \R^{n_{l-1}}$ and $\psi_{d}^{(l)}:\R^{n_{l-1}}\to \R^{n_{l}}$, then the dimensions $n_l$ of the layers must be non-decreasing with $r = n_0 \leq n_1 \leq \cdots \leq n_L = n$.

There are two main ingredients in our approach to constructing layers with the desired properties.  The first is a pair of smooth activation functions $\sigma_+$ and $\sigma_-$ that act element-wise on vectors and satisfy $\sigma_- \circ \sigma_+ = Id$.  The second is a constraint on the weight matrices $\Phi_l, \Psi_l \in \R^{n_{l}\times n_{l-1}}$, such that they satisfy the biorthogonality condition
$\Psi_l^T \Phi_l = I$.  These two ingredients are explained in the following subsections.
Once these are defined, we construct the layers of the encoder and decoder according to
\begin{empheq}[box=\widefbox]{equation}
\begin{aligned}
    \psi_{e}^{(l)}(x^{(l+1)}) &= \sigma_-\big(\Psi_l^T (x^{(l+1)} - b_l)\big), \\
    \psi_{d}^{(l)}(z^{(l-1)}) &= \Phi_l \sigma_+(z^{(l-1)}) + b_l,
\end{aligned}
    \label{eqn:layer_definition}
\end{empheq}
where $b_l$ are bias vectors.  The resulting layer transformation then satisfies
$\psi_e^{(l)} \circ \psi_d^{(l)} = Id$, as desired.

\begin{remark}[Parameter-dependent projections]
    Intrinsic manifolds often depend on system parameters.
    A parameter-dependent projection can be obtained by allowing the biases $b_l$ to be functions of a vector of parameters $q$.
    Specifically, we can define $b_l = W_l q + \tilde{b}_l$ where $W_l$ and $\tilde{b}_l$ are trainable weights and biases.
\end{remark}

% \end{widetext}

\begin{remark}
    \label{rem:outer_layer_width}
    By definition, the decoder reconstructs states in an affine subspace of dimension $n_{L-1}$.
    Therefore, $n_{L-1}$ should be chosen based on Kolmogorov $n$-width considerations so that state data from the system can be accurately reconstructed in an affine subspace of dimension $n_{L-1}$.
\end{remark}

\subsection{Invertible, smooth activation functions}
\label{subsec:activation_funs}
\begin{figure}
    \centering
    % \begin{tikzonimage}[trim=71 5 59 9, clip=true, width=0.3\textwidth]{Figures/activation_fun.eps}%[tsx/show help lines]
    %     \node[rotate=0] at (0.55, 0.8) {$\sigma_+$};
    %     % \node[rotate=0] at (0.2, 0.45) {$\sigma_+$};
    %     \node[rotate=0] at (0.8, 0.55) {$\sigma_-$};
    %     % rotated coordinates
    %     % \draw[->] (.52, .52) -- (.82, .82);
    %     % \node[rotate=45, anchor=north] at (0.67, 0.67) {$\tilde{x}$};
    %     % \draw[->] (.52, .52) -- (.22, .82);
    %     % \node[rotate=45, anchor=east] at (0.37, 0.67) {$\tilde{y}$};
    %     % original coordinates
    %     \node[rotate=0] at (0.5, 0.01) {$x$};
    %     \node[rotate=90] at (0.01, 0.5) {$y$};
    %     % asymptote angle
    %     \draw[->] (0.7,0.7) arc[radius=0.2*sqrt(2), start angle=45, end angle=67.5];
    %     \node[rotate=0] at (0.69,0.76) {$\alpha$};
    % \end{tikzonimage}
    \begin{tikzonimage}[trim=71 5 59 9, clip=true, width=0.3\textwidth]{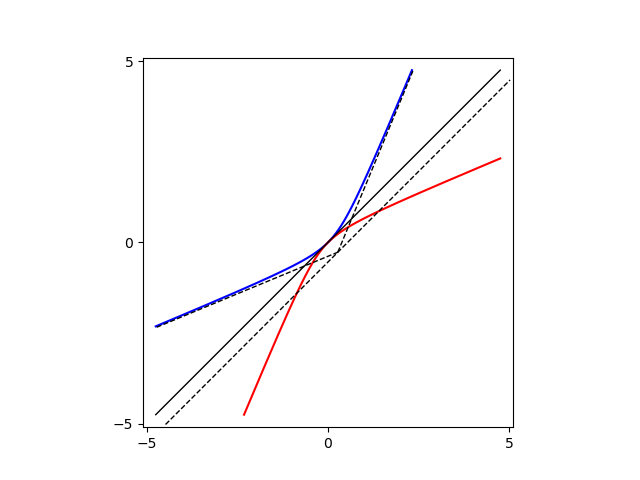}%[tsx/show help lines]
        \node[rotate=0] at (0.57, 0.8) {$\sigma_+$};
        % \node[rotate=0] at (0.2, 0.45) {$\sigma_+$};
        \node[rotate=0] at (0.8, 0.57) {$\sigma_-$};
        % rotated coordinates
        % \draw[->] (.52, .52) -- (.82, .82);
        % \node[rotate=45, anchor=north] at (0.67, 0.67) {$\tilde{x}$};
        % \draw[->] (.52, .52) -- (.22, .82);
        % \node[rotate=45, anchor=east] at (0.37, 0.67) {$\tilde{y}$};
        % original coordinates
        \node[rotate=0] at (0.5, 0.01) {$x$};
        \node[rotate=90] at (0.01, 0.5) {$y$};
        % asymptote angle
        \draw[->] (0.7+0.02, 0.7-0.02) arc[radius=0.2*sqrt(2), start angle=45, end angle=67.5];
        \node[rotate=0] at (0.69, 0.76) {$\alpha$};
    \end{tikzonimage}
    \caption{The smooth, invertible activation functions $\sigma_{\pm}$ are constructed geometrically from a hyperbola with conjugate axis (dashed black line) parallel to $y=x$ (black line) with asymptotes (dashed black lines) at angle $\alpha$ from the conjugate axis. The graph of $\sigma_+$ is the upper branch of this hyperbola (blue curve) and the graph of $\sigma_-$ (red curve) is obtained by reflecting across the line $y=x$.}
    \label{fig:activation_functions}
\end{figure}
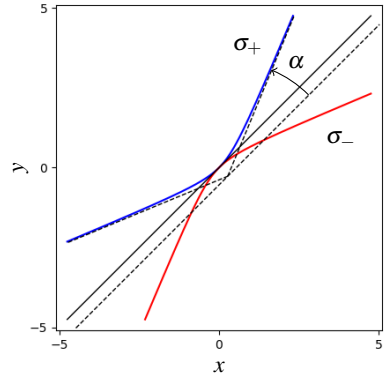

Here we define the smooth, invertible activation functions $\sigma_{\pm} :\R \to \R$ to be used in the encoder and decoder.
Geometrically, the condition that $\sigma_+$ and $\sigma_-$ are inverses is equivalent to the condition that their graphs are reflections about the line $y=x$ in $\R^2$.
In rotated coordinates $(\tilde{x}, \tilde{y}) = \frac{\sqrt{2}}{2}(x+y, y-x)$ where the line $y=x$ corresponds with $\tilde{y} = 0$, we let the graph of $\sigma_+$ be the upper branch ($\tilde{y} > 0$) of the hyperbola defined by
\begin{equation}
    \frac{\big(\tilde{y} + \sin(\alpha)\big)^2}{\sin^2(\alpha)} - \frac{\tilde{x}^2}{\cos^2(\alpha)} = 1,
    \label{eqn:rotated_activation_function}
\end{equation}
where $0 < \alpha < \pi/4$.
To form $\sigma_-$, we flip the sign of $\tilde{y}$.
In (\ref{eqn:rotated_activation_function}), $\tilde{y}$ is shifted by $\sin(\alpha)$ in order to ensure that $\sigma_{\pm}(0) = 0$.
By symmetry, the derivatives satisfy $\sigma_{\pm}'(0) = 1$.
As shown in Figure~\ref{fig:activation_functions}, the upper and lower branches of this hyperbola are reflections about the axis $y=x$ with asymptotes at angle $\alpha$ from this axis.
The condition that $0 < \alpha < \pi/4$ ensures that these branches are graphs of well-defined functions $\sigma_{\pm}$.
In the results shown in Section~\ref{sec:noack_case_study}, we take $\alpha=\pi/8$.
Rotating back to $(x,y)$ coordinates, the activation functions are given by
% \begin{widetext}
\begin{empheq}[box=\fbox]{multline}
    \sigma_{\pm}(x) 
    =  \frac{b x}{a} \mp \frac{\sqrt{2}}{a\sin(\alpha)} \\
    \pm \frac{1}{a} \sqrt{ \left( \frac{2x}{\sin(\alpha)\cos(\alpha)} \mp \frac{\sqrt{2}}{\cos(\alpha)} \right)^2 + 2a }, \\ 
    \mbox{where} \quad
    \left\{ \begin{matrix} a = \csc^2(\alpha) - \sec^2(\alpha) \\ b = \csc^2(\alpha) + \sec^2(\alpha) \end{matrix} \right..
    \label{eqn:activation_functions}
\end{empheq}
% \end{widetext}
Since $0 < a < b$, these functions are well-defined for all $x\in\R$ and are infinitely continuously differentiable.
Examining their graphs in Figure~\ref{fig:activation_functions}, we also observe that they resemble smooth, symmetric versions of ``leaky'' rectified linear units (ReLU) \cite{He2015delving} common in deep learning applications.

\subsection{Weight matrix biorthogonality}
\label{subsec:biorthogonality}

The layers of the encoder and decoder in (\ref{eqn:layer_definition}) are defined using biorthogonal weight matrices, that is, pairs of matrices $\Phi, \Psi \in \R^{n\times r}$, $n \geq r \geq 1$, satisfying $\Psi^T \Phi = I$.
Here, we describe how to enforce this constraint during training.
In Appendix~\ref{app:biorthogonal_mfd} we show that these matrices form a smooth, properly embedded submanifold $\mcal{B}_{n, r}$ of $\R^{n\times r}\times \R^{n\times r}$ with dimension $\dim \mcal{B}_{n,r} = 2 n r - r^2$.

A simple way to optimize the weight matrices on the biorthogonal manifold using existing optimizers for Euclidean spaces is to rely on an over-parametrization.
In particular, we over-parametrize $\mcal{B}_{n,r}$ over an open subset
\begin{equation}
    D_+(\Pi_{n,r}) = \left\{ (\tilde{\Phi}, \tilde{\Psi}) \in \R^{n,r} \times \R^{n,r}  \ : \ \det(\tilde{\Psi}^T\tilde{\Phi}) > 0 \right\}
\end{equation}
of the Euclidean space $\R^{n,r} \times \R^{n,r}$ using a projection map $\Pi_{n,r}: D_+(\Pi_{n,r}) \to \mcal{B}_{n,r}$ defined by
\begin{empheq}[box=\widefbox]{equation}
    \Pi_{n,r}(\tilde{\Phi},\tilde{\Psi}) = \big(\tilde{\Phi}(\tilde{\Psi}^T\tilde{\Phi})^{-1}, \ \tilde{\Psi} \big).
    \label{eqn:projection_onto_Bnr}
\end{empheq}
Indeed, one can easily check that this map is smooth, surjective, and idempotent $\Pi_{n,r} \circ \Pi_{n,r} = \Pi_{n,r}$.
By composing an optimization objective function $J:\mcal{B}_{n,r} \to \R$ with the over-parametrization we produce a new objective
\begin{equation}
    \bar{J} := J \circ \Pi_{n,r} : D_+(\Pi_{n,r}) \to \R
\end{equation}
defined on an open subset of the Euclidean space $\R^{n,r} \times \R^{n,r}$.
Theorem~\ref{thm:submersion_onto_Bnr} in Appendix~\ref{app:biorthogonal_mfd} says that this is locally equivalent (by a smooth change of coordinates) to introducing $r^2$ additional optimization variables on which the cost function does not depend. 
Consequently the over-parametrization does not introduce any new critical points into the optimization problem in the sense that the gradient of the original objective $\grad J(\Phi, \Psi)$ vanishes if and only if the gradient of the composition $\grad \bar{J}(\tilde{\Phi}, \tilde{\Psi})$ vanishes at every element $(\tilde{\Phi}, \tilde{\Psi})$ in the preimage fiber $\Pi_{n,r}^{-1}(\Phi, \Psi)$.

During optimization we must ensure that the representatives $(\tilde{\Phi}, \tilde{\Psi})$ of the weight matrices $(\Phi, \Psi) = \Pi_{n,r}(\tilde{\Phi}, \tilde{\Psi})$ remain in the domain $D_+(\Pi_{n,r})$ and do not approach its boundary.
To do this, regularization functions for each layer of the network are added to the cost function minimized during training.
The regularization we use for each layer is given by
\begin{empheq}[box=\widefbox]{equation}
    R(\tilde{\Phi}, \tilde{\Psi}) 
    = \big\Vert \tilde{\Psi}^T\tilde{\Phi} - I \big\Vert_F^2 \big\Vert (\tilde{\Psi}^T\tilde{\Phi})^{-1} \big\Vert_F^2
    \label{eqn:regularization}
\end{empheq}
Evidently, this function is well-defined and smooth on $D_+(\Pi_{n,r})$.
It takes its minimum value of zero if and only if $(\tilde{\Phi}, \tilde{\Psi}) \in \mcal{B}_{n,r}$ and it blows up to $+\infty$ whenever $\tilde{\Psi}^T\tilde{\Phi}$ approaches a singular matrix.
Therefore, including this regularization term in the cost function forces the optimization iterates $\tilde{\Phi}, \tilde{\Psi}$ to remain near (in fact, to approach) $\mcal{B}_{n,r}$ without approaching the boundary of $D_+(\Pi_{n,r})$.
Note that the weight matrices $(\Phi, \Psi) = \Pi_{n,r}(\tilde{\Phi}, \tilde{\Psi})$ of the autoencoder always remain in the biorthogonal manifold.

Our analysis in Appendix~\ref{app:biorthogonal_mfd} also shows that the optimization domain $D_+(\Pi_{n,r})$ is connected when $n > r$.
This means that restricting the optimizer to this domain does not cut off access to any part of the biorthogonal manifold by an optimization algorithm that follows a continuous path or proceeds in small steps.
On the other hand, when $n=r$, the birothogonal manifold consists of pairs $(\Phi, \Phi^{-1})$, where $\Phi$ are invertible $n\times n$ matrices.
In this case, $D_+(\Pi_{n,n})$ and $\mcal{B}_{n,n}$ consist of two disjoint connected components corresponding to matrices with positive and negative determinants.
However, we show in Appendix~\ref{app:biorthogonal_mfd} that this is of no consequence for the optimization of the autoencoder's weights because any choice for the signs of the determinants in the square layers can be achieved without altering the projection $P = \psi_d \circ \psi_e$.
Hence, one does not have to explore other connected components during optimization.

We summarize the training procedure for our autoencoder in Algorithm~\ref{alg:autoencoder_training}.
The specific cost functions and the types of training data we employ will be discussed in Section~\ref{sec:cost_funs}.
These cost functions $J$ can depend directly on the autoencoder $P = \psi_d \circ \psi_e$, its derivatives, the biorthogonal weights $\Phi_l, \Psi_l$ and biases $b_l$ in each layer, the data in the minibatch, the FOM, or other parameters, but not the weight matrix representatives $\tilde{\Phi}_l, \tilde{\Psi}_l$.
In Section~\ref{subsec:noack_train_procedure} we discuss specific details of the training procedure for our main numerical example including the construction of minibatches and the choice of optimizer and optimization parameters such as the learning rate.

\begin{algorithm}[H]
\caption{Autoencoder training procedure}\label{alg:autoencoder_training}
\begin{algorithmic}[1]
\STATE{\textbf{input:} 
layer widths $r=n_0 \leq n_1 \leq \cdots \leq n_L = n$, activation function asymptote angle $0 < \alpha < \pi/4$, training data, cost function $J$, regularization strength $\beta > 0$, number of training epochs,
initial biorthogonal weight matrices $(\Phi_l, \Psi_l)\in\mcal{B}_{n_{l}, n_{l-1}}$, and initial bias vectors $b_l\in\R^{n_l}$.}
\STATE{initialize $(\tilde{\Phi}_l, \tilde{\Psi}_l) = (\Phi_l, \Psi_l)\in D(\Pi_{n_l, n_{l-1}})$ for $l=1, \ldots, L$}
\FOR{$\text{epoch}=1,2, \ldots, (\text{num. epochs})$}
    \STATE{randomly split training data set into minibatches}
    \FOR{each minibatch}
        \STATE{construct autoencoder $P = \psi_d \circ \psi_e$ with weights $(\Phi_l, \Psi_l) = \Pi_{n_{l}, n_{l-1}}(\tilde{\Phi}_l, \tilde{\Psi}_l)$ and bias vectors $b_l$ in layers $l=1, \ldots, L$ defined by \eqref{eqn:layer_definition} with activation functions in \eqref{eqn:activation_functions}.\label{algstep:construct_AE}}
        \STATE{use $P$ and the minibatch to compute the regularized cost 
        $$\bar{J}_{\beta}\big( [\tilde{\Phi}_l, \tilde{\Psi}_l, b_l]_{l=1}^L \big) = J(P, \text{minibatch}) + \beta \sum_{l=1}^L R(\tilde{\Phi}_l, \tilde{\Psi}_l)$$
        % \begin{multline}
        %     \bar{J}_{\beta}\big( [\tilde{\Phi}_l, \tilde{\Psi}_l, b_l]_{l=1}^L \big) 
        %     = J\big([\Phi_l, \Psi_l, b_l]_{l=1}^L, \text{minibatch}\big) \\
        %     + \beta \sum_{l=1}^L R(\tilde{\Phi}_l, \tilde{\Psi}_l)
        % \end{multline}
        % \begin{multline*}
        %     \bar{J}_{\beta}\big( [\tilde{\Phi}_l, \tilde{\Psi}_l, b_l]_{l=1}^L \big) 
        %     = J\big([\Pi_{n_{l}, n_{l-1}}(\tilde{\Phi}_l, \tilde{\Psi}_l), b_l]_{l=1}^L, \text{minibatch}\big) \\
        %     + \beta \sum_{l=1}^L R(\tilde{\Phi}_l, \tilde{\Psi}_l)
        % \end{multline*}
        \label{algstep:compute_cost}}
        \STATE{compute the gradient of $\bar{J}_{\beta}$ with respect to each $\tilde{\Phi}_l, \tilde{\Psi}_l, b_l$\label{algstep:compute_grad}}
        \STATE{use the gradients and an optimizer such as Adam\cite{Kingma2014Adam} to update each $\tilde{\Phi}_l, \tilde{\Psi}_l, b_l$}
    \ENDFOR
\ENDFOR
\RETURN{The autoencoder $P$ along with the final weights $(\Phi_1, \Psi_1) = \Pi_{n_l, n_{l-1}}(\tilde{\Phi}_l, \tilde{\Psi}_l)$ and biases $b_l$ defining each layer $l=1, \ldots, L$.}
\end{algorithmic}
\end{algorithm}

\begin{remark}
    \label{rem:optimizing_on_Bnr_directly}
    Another approach is to optimize the autoencoder's weights directly on the biorthogonal manifold using gradient-based techniques together with an appropriate retraction and vector transport \cite{Absil2009optimization}.
    In fact, the over-parametrization map $\Pi_{n,r}$ yields a ``projection-like retraction'' \cite{Absil2012projection} on $\mcal{B}_{n,r}$.
    The projection map onto the tangent space of $\mcal{B}_{n,r}$ given by Theorem~\ref{thm:biorthogonal_manifold} in Appendix~\ref{app:biorthogonal_mfd} also yields a vector transport on $\mcal{B}_{n,r}$.
    This approach is discussed in Section~3.4 of Otto's thesis \cite{Otto2022advances}.
    % Because this approach does not add additional parameters, it avoids introducing singularities into the objective function's Hessian and is therefore better suited to high-order optimization algorithms that rely on inverting the Hessian or an approximation of the Hessian.
    However, it is difficult to implement in existing neural network optimizers such as PyTorch \cite{PyTorch} and TensorFlow \cite{TensorFlow}, motivating the use of our simple over-parametrization instead.
\end{remark}

\subsection{Preserving an equilibrium point}
\label{subsec:preserving_an_equilibrium}
In certain cases such as in control applications, it is important for the reduced-order model to preserve a known equilibrium point of the system.
% Via a coordinate shift, we assume this equilibrium is at the origin.
To ensure that our nonlinear projection-based ROM has the same equilibrium point, it suffices to ensure that the equilibrium $x_{\text{eq}}$ is contained in the learned manifold parametrized by the decoder.
To do this, we obtain $\psi_d(0) = x_{\text{eq}}$ by constraining the bias vector in the final layer to be
\begin{equation}
    b_{L} = x_{\text{eq}} - \Phi_L \sigma_{+} \circ \psi_d^{(L-1)} \circ \cdots \circ \psi_d^{(1)}(0).
    \label{eqn:bias_constraint}
\end{equation}
The resulting equilibrium point of the ROM (\ref{eqn:ROM_in_latent_space}) is located at the origin in the latent space of the autoencoder.
Note that it is always possible to shift an equilibrium point to the origin $x_{\text{eq}} = 0$ by a change of coordinates in \eqref{eqn:FOM}.

\subsection{Enforcing linear constraints on state vectors}
\label{subsec:linear_constraints}
% \seoremark{This section aims to address comment~1 by reviewer~1.}
Suppose we know that the state vectors $x$ of the system \eqref{eqn:FOM} satisfy a collection of linear constraints $\mcal{L} x = 0$.
Examples included certain boundary conditions for solutions of partial differential equations as well as incompressibility constraints in fluid flows.
To ensure that all projected states $P(x)$ also satisfy these constraints, it suffices to ensure that the weight matrix and bias vector defining the last layer of the decoder satisfy $\mcal{L} \Phi_L = 0$ and $\mcal{L} b_L = 0$.
Examining \eqref{eqn:layer_definition}, we see that this yields $\mcal{L} \psi_d^{(L)}(z^{(L-1)}) = 0$, which implies that $\mcal{L} P(x) = 0$ for every $x$.
During training (see Algorithm~\ref{alg:autoencoder_training}), we optimize representatives $(\tilde{\Phi}_L, \tilde{\Psi}_L) \in D (\Pi_{n_L, n_{L-1}})$ of the weight matrices $(\Phi_L, \Psi_L) = \Pi_{n_L, n_{L-1}}(\tilde{\Phi}_L, \tilde{\Psi}_L)$.
Enforcing the linear constraint $\mcal{L} \tilde{\Phi}_L = 0$ on the representative automatically ensures that $\mcal{L} \Phi_L = 0$, as one can easily verify from \eqref{eqn:projection_onto_Bnr}.
In practice, 
% $\mcal{L} \tilde{\Phi}_L = 0$ and $\mcal{L} b_L = 0$ 
\begin{equation}
    \mcal{L} \tilde{\Phi}_L = 0
    \qquad \mbox{and} \qquad
    \mcal{L} b_L = 0
\end{equation}
can be enforced either by parametrizing $b_L$ and the columns of $\tilde{\Phi}_L$ in a basis for $\Null(\mcal{L})$, or by employing projected gradient descent methods to constrain the iterates within $\Null(\mcal{L})$

\subsection{Initialization}
\label{subsec:initialization}

In Figure \ref{fig:activation_functions}, we see that each activation function $\sigma_+$ and $\sigma_-$ can produce an output of larger magnitude than the input, and repeated activations in deep networks can result in much greater amplification.
In addition, linear layers with operator norm greater that unity will further enlarge the output magnitude.
These effects can lead to very large initial loss, which interferes with training.
To address this issue, we initialize the network's weights, $\Phi_l$ and $\Psi_l$, such that $\|\Phi_l\|_2 = 1$ and $\|\Psi_l^T\|_2=1$. 
In particular, we randomly sample a square matrix from the orthogonal group and take the first $n_l$ columns to construct $\Phi$ and $\Psi = \Phi$.

Regardless of whether we preserve the equilibrium point via a constraint, as discussed in Section~\ref{subsec:preserving_an_equilibrium}, it is usually advantageous for the network to have the property that $P(0)=0$ at initialization. This property is satisfied if we set all biases to zero at initialization since  $\sigma_+(0) = \sigma_-(0) = 0$.

% \subsection{Training procedure}

\section{Optimization objectives}
\label{sec:cost_funs}

Choosing an appropriate optimization objective is crucial for learning projections that yield accurate reduced-order models.
Typically, the parameters $\theta$ consisting of the weights and biases in an autoencoder are optimized in order to minimize the average reconstruction error
\begin{equation}
    J_{\text{Rec}}(P) = \E_x \left[ \left\Vert x - P(x) \right\Vert^2 \right]
    \label{eqn:Rec_loss}
\end{equation}
over some distribution of states $x$.
For example, this might be an empirical distribution of states sampled along trajectories of interest from the full-order model.
However, the loss function~\eqref{eqn:Rec_loss} encourages the projection to simply map each point $x$ in the support of the distribution to its nearest point on the learned manifold $\hat{\mcal{M}} = \Range(P)$.
In a tubular neighborhood of $\hat{\mcal{M}}$ this yields an orthgonal projection in the sense that the line segment in $\R^n$ (in the Riemannian case, the minimizing geodesic) connecting each $x$ in the tubular neighborhood to $P(x)$ lies in the fiber of $P(x)$ and is orthogonal to $T_{P(x)}\hat{\mcal{M}}$ (see \citet{Lee2013introduction} or \citet{Guillemin1974differential}).
As we discussed in Section~\ref{sec:nlprom} (see Figure~\ref{fig:ic_proj_prob}), this is not always ideal for modeling the dynamics since the truncation does not account for coordinates that have a large influence on the future behavior of the system.
In this section, we develop alternative objectives (loss functions) for training the autoencoder that account for this kind of sensitivity.

\subsection{Reconstruction and Velocity Projection (RVP) loss}
\label{subsec:RVP}

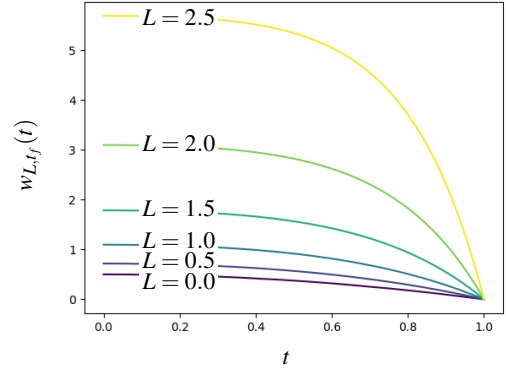
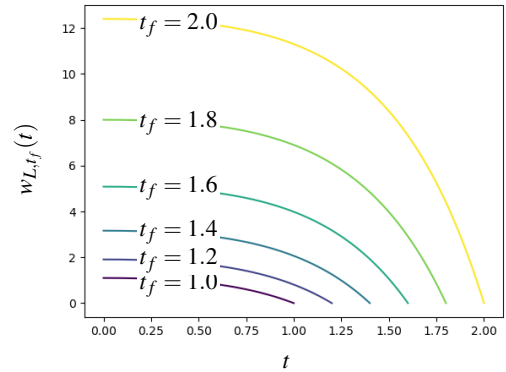
\begin{figure}
    \centering
    % \begin{minipage}{0.48\textwidth}
    \subfloat[$t_f=1$, varying $L$]{
        \centering
        \begin{tikzonimage}[trim=0 0 0 0, clip=true, width=0.40\textwidth]{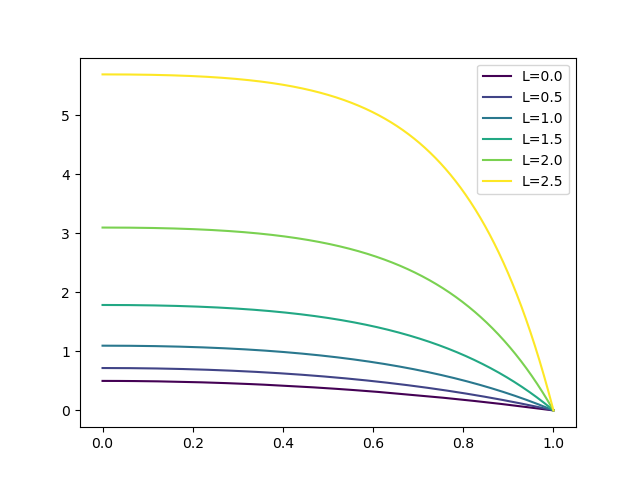}%[tsx/show help lines]
        \node[rotate=0] at (0.5, 0.00) {$t$};
        \node[rotate=90] at (0.02, 0.5) {$w_{L,t_f}(t)$};
        \node [fill=white, inner sep=1pt] at (0.3, 0.19) {$L=0.0$};
        \node [fill=white, inner sep=1pt] at (0.3, 0.24) {$L=0.5$};
        \node [fill=white, inner sep=1pt] at (0.3, 0.29) {$L=1.0$};
        \node [fill=white, inner sep=1pt] at (0.3, 0.37) {$L=1.5$};
        \node [fill=white, inner sep=1pt] at (0.3, 0.53) {$L=2.0$};
        \node [fill=white, inner sep=1pt] at (0.3, 0.84) {$L=2.5$};
        \fill [white] (0.73 ,0.59) rectangle (0.895,0.87);
        \end{tikzonimage}
        }\\
    % \end{minipage} \hfill
    % \hspace{0.025\textwidth}
    % \begin{minipage}{0.48\textwidth}
    \subfloat[$L=1$, varying $t_f$]{
        \centering
        \begin{tikzonimage}[trim=0 0 0 0, clip=true, width=0.40\textwidth]{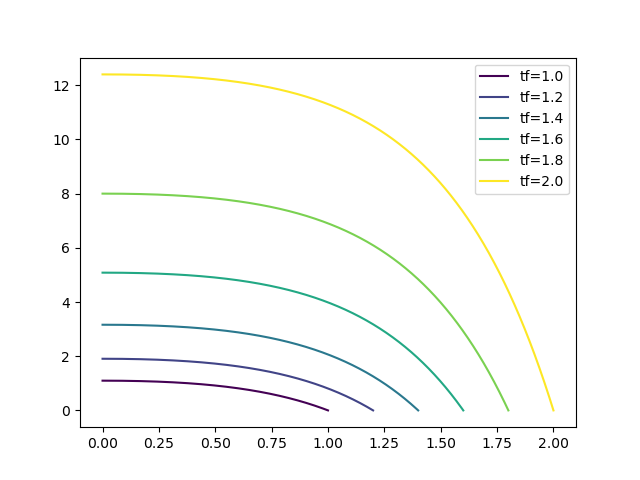}%[tsx/show help lines]
        \node[rotate=0] at (0.5, 0.00) {$t$};
        \node[rotate=90] at (0.02, 0.5) {$w_{L,t_f}(t)$};
        \node [fill=white, inner sep=1pt] at (0.3, 0.19) {$t_f=1.0$};
        \node [fill=white, inner sep=1pt] at (0.3, 0.25) {$t_f=1.2$};
        \node [fill=white, inner sep=1pt] at (0.3, 0.32) {$t_f=1.4$};
        \node [fill=white, inner sep=1pt] at (0.3, 0.43) {$t_f=1.6$};
        \node [fill=white, inner sep=1pt] at (0.3, 0.59) {$t_f=1.8$};
        \node [fill=white, inner sep=1pt] at (0.3, 0.83) {$t_f=2.0$};
        \fill [white] (0.73 ,0.59) rectangle (0.895,0.87);
        \end{tikzonimage}
    }
    % \end{minipage}
    \caption{Velocity projection weight function given by~\eqref{eqn:velocity_error_weight_function} as the parameters $L$ and $t_f$ are varied. The weight function increases exponentially when these parameters are increased.}
    \label{fig:v_proj_wfun}
\end{figure}

One way to account for the dynamics is to penalize the difference between the time derivative of the reduced-order model (\ref{eqn:P_ROM}) and the time derivative along projected trajectories of the full-order model (\ref{eqn:FOM}).
If $x(t)$ is a trajectory of the FOM generating output $y(t)$, then the time derivative of the projected trajectory $x_P(t) = P(x(t))$ is
\begin{equation}
    \ddt x_P(t)
    = \D P(x(t)) \ddt x(t).
    % = \D P(x(t)) f(x(t), u(t)).
\end{equation}
At the same point $x_P(t)$, the time derivative of the ROM (\ref{eqn:P_ROM}) is given by
\begin{equation}
    \hat{f}_P(x_P(t), u(t)) 
    = \D P(x_P(t)) f(x_P(t), u(t) ).
\end{equation}
These two quantities are equal for all $t$ if and only if the trajectory $\hat{x}(t)$ of the ROM agrees with the projected trajectory $x_P(t)$.
The following proposition shows how the integrated square error between these trajectories is bounded by a weighted integral of the square projection error for the time derivatives.
\begin{proposition}[Weighted velocity projection error]
    \label{prop:weighted_VPE}
    Let $x(t)$, $t\in [0, t_f]$ be a trajectory of (\ref{eqn:FOM}) and let $P:\R^n\to\R^n$ be a smooth projection map.
    Suppose that (\ref{eqn:P_ROM}) has a unique solution $\hat{x}(t)$ over the same time interval.
    If $x_P(t) = P(x(t))$ and $\hat{x}(t)$ are contained within a subset $\mcal{U}\subset \R^n$ over which $x \mapsto \hat{f}_P(x,u(t)) = \D P(x) f(x, u(t))$ has Lipshitz constant $L$ for every $t\in [0,t_f]$ then
    \begin{multline}
        \int_{0}^{t_f} \big\Vert x_P(t) - \hat{x}(t) \big\Vert^2 \td t \\
        \leq \int_{0}^{t_f} w_{L,t_f}(t) \Big\Vert \ddt x_P(t) - \hat{f}_P(x_P(t),u(t)) \Big\Vert^2 \td t,
        \label{eqn:velocity_error_weighting}
    \end{multline}
    where
    \begin{equation}
        w_{L,t_f}(t) = \frac{1}{4L^2} \left[ \left( e^{2L t_f} - 2 L t_f \right) - \left( e^{2L t} - 2 L t \right) \right].
        \label{eqn:velocity_error_weight_function}
    \end{equation}
\end{proposition}
\begin{proof}
        % The proof involves using a modified Gr\"{o}nwall-Bellman inequality and Cauchy-Schwarz to bound $\Vert x_P(t) - \hat{x}(t) \Vert^2$ in terms of a $t$-dependent weight times the square integral of $\Vert \ddt x_P(\tau) - \hat{f}_P(x_P(\tau),u(\tau)) \Vert$ over the interval $\tau \in [0, t]$.
        % Integrating this bound and exchanging the order of integration yields the desired result.
        The result essentially follows from a Gr\"{o}nwall-Bellman-type inequality.
        We provide the details in Appendix~\ref{app:proofs}
\end{proof}
The significance of this result is that it tells us how to properly weight the velocity projection error in formulating optimization objectives.
While we are primarily interested in the error between the trajectory of the ROM and the projected trajectory of the FOM, 
velocity projection error is a more convenient quantity to optimize because it does not involve integrating the ROM forward in time.
Since it is difficult to determine the Lipschitz constant $L$ in practice, we treat it as a parameter when using Proposition~\ref{prop:weighted_VPE} as a guide to formulate objective functions for optimization.
In this case, our choice of $L$ reflects the rate at which we expect nearby trajectories of the ROM to diverge.
The weight function is plotted in Figure~\ref{fig:v_proj_wfun} over a range of values for its parameters $L$ and $t_f$.
We observe that in the limit as $L\to 0$, the weight function becomes
\begin{equation}
    \lim_{L\to 0} w_{L,t_f}(t) = \frac{1}{2} \left( t_f^2 - t^2 \right).
\end{equation}
On the other hand, the weight function increases exponentially with $L t_f$, so we must be somewhat careful that $L t_f$ is not too large.

If the projected FOM trajectory $x_P(t)$ agrees with the trajectory of the ROM $\hat{x}(t)$, then the error between the output of the ROM $\hat{y}(t) = g(\hat{x}(t))$ and the output of the FOM $y(t) = g(x(t))$ is due only to the difference between $g(x_P(t))$ and $y(t)$.
We can measure this using a reconstruction loss resembling (\ref{eqn:Rec_loss}).
Therefore, we combine this reconstruction error with the bound on the trajectory error from Proposition~\ref{prop:weighted_VPE} along trajectories $x(t)$ drawn from a given distribution over initial conditions and input signals.
Combining the reconstruction error and a constant $\gamma \geq 0$ times the weighted velocity projection error into a single loss function, we seek to minimize
\begin{empheq}[box=\fbox]{multline}
    J_{\text{RVP}, \gamma}(P) = \E_{x_0, u}\Bigg[ \frac{1}{t_f} \int_{0}^{t_f} \bigg( \big\Vert y(t) - g(x_P(t)) \big\Vert^2 \\
    + \gamma w_{L,t_f}(t) \Big\Vert \ddt x_P(t) - \hat{f}_P(x_P(t), u(t)) \Big\Vert^2 \bigg) \td t \Bigg].
    \label{eqn:RVP_loss}
\end{empheq}
Here, we select $\gamma \geq 0$ to determine the strength of the velocity projection error term.
This parameter can be selected based on the amplification of state errors through the output map $g$ in~(\ref{eqn:FOM}).
For example, if the state is observed though a linear map $g: x \mapsto C x$ with operator norm $\Vert C \Vert_{\text{op}}$, then setting $\gamma = \Vert C \Vert_{\text{op}}^2$ can be used to bound the square error of the output using the RVP loss:
\begin{align*}    
    \int_{0}^{t_f} \frac{1}{2} &\Vert y(t) - \hat{y}(t) \Vert^2 \td t \\
    &\leq \int_{0}^{t_f} \bigg( \Vert y(t) - g(x_P(t)) \Vert^2 + \Vert C \Vert_{\text{op}}^2 \Vert x_P(t) - \hat{x}(t) \Vert^2 \bigg) \td t \\
    &\leq \int_{0}^{t_f} \bigg( \big\Vert y(t) - g(x_P(t)) \big\Vert^2 \\
    &\qquad + \Vert C \Vert_{\text{op}}^2 w_{L,t_f}(t) \Big\Vert \ddt x_P(t) - \hat{f}_P(x_P(t), u(t)) \Big\Vert^2 \bigg) \td t.
\end{align*}
Here, the first inequality follows from the convexity of squared Euclidean norm and the second inequality follows from Proposition~\ref{prop:weighted_VPE}.
In cases where the relative (rather than absolute) square error is of interest, the terms in (\ref{eqn:RVP_loss}) inside the expectation or the integrand can be normalized by magnitudes of ground truth values for $y(t)$ and $\dot{x}(t)$.

Unlike the reconstruction loss, the reconstruction and velocity projection (RVP) loss (\ref{eqn:RVP_loss}) requires us to evaluate the governing equations of the FOM.
In order to compute the gradient of this loss function, we must be able to act on vectors with the transposes (adjoints) of Jacobians derived from the FOM, i.e., to compute $\left( \frac{\partial}{\partial x}f(x,u) \right)^T v$ and $\left( \frac{\partial}{\partial x}g(x) \right)^T w$ for vectors $v\in\R^n$ and $w\in\R^m$.

\begin{figure*}
    \centering
    \includegraphics[scale=0.4]{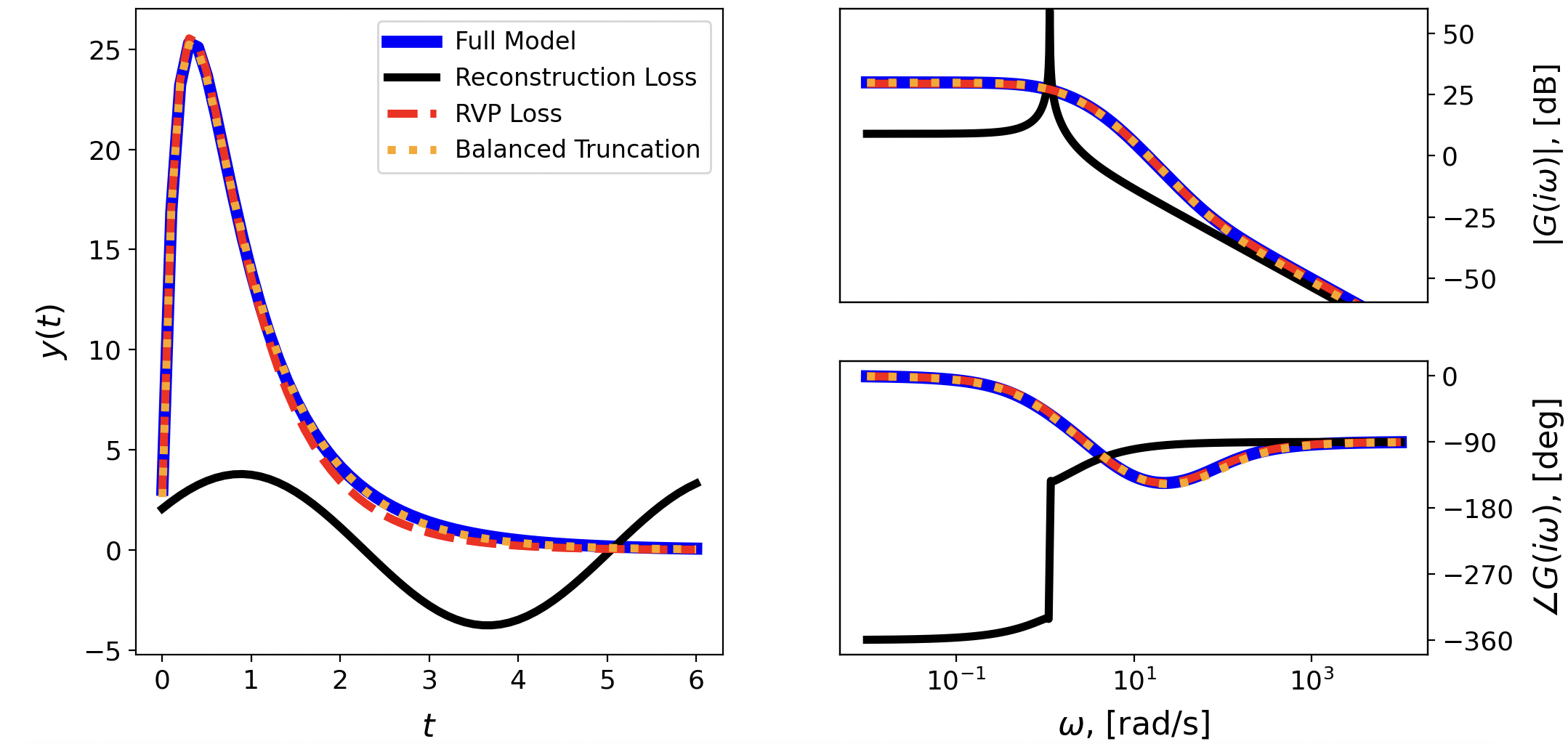}
    \caption{Impulse response (left) and frequency response (right) for Example~\ref{ex:linear}, comparing the full model and three second-order reduced-order modeling approaches: reconstruction loss (POD), reconstruction and velocity projection loss (RVP), and balanced truncation.}
    \label{fig:linear_step_response}
\end{figure*}

The upshot of this added complexity is that the RVP loss can account for system nonnormality, as the following example illustrates.
\begin{example}[RVP loss for a nonnormal linear system]
    \label{ex:linear}
    We consider the problem of finding a two-dimensional linear projection for the nonnormal linear system
    \begin{equation}
    \begin{aligned}
        \dot{x}_1 &= - x_1 + 100 x_3 + u \\
        \dot{x}_2 &= - 2 x_2 + 100 x_3 + u \\
        \dot{x}_3 &= - 5 x_3 + u \\
        y &= x_1 + x_2 + x_3,
    \end{aligned}
    \label{eqn:nonnormal_linear_toy_model}
    \end{equation}
    discussed as an example in \citet{Holmes2012turbulence}.
    In response to an impulse, the state $x_3$ decays rapidly to zero and exerts a large influence on $x_1$ and $x_2$, causing them to experience a large transient growth before eventually decaying.
    in \citet{Holmes2012turbulence} it is shown that POD, while being optimal with respect to reconstruction loss (\ref{eqn:Rec_loss}), yields an orthogonal projection subspace closely aligned with the $x_1, x_2$ coordinate plane and therefore ignores the important influence of $x_3$.
    The resulting model does not experience the large transient growth present in the impulse response of (\ref{eqn:nonnormal_linear_toy_model}).
    To see why optimizing the projection with respect to RVP loss can improve this situation, consider the orthogonal projection $P_{1,2}$ onto the $x_1, x_2$ coordinate plane in $R^3$ and a state $x = (x_1, x_2, x_3)$ along the impulse-response trajectory of (\ref{eqn:nonnormal_linear_toy_model}).
    While the reconstruction error $x - P_{1,2}x = (0,0,x_3)$ is small, the velocity projection error
    \begin{equation}
        P_{1,2} f(x,0) - P_{1,2} f(P_{1,2} x,0) =
        % \begin{bmatrix}
        %     - x_1 + 100 x_3 \\
        %     - 2 x_2 + 100 x_3 \\
        %     0
        % \end{bmatrix} - 
        % \begin{bmatrix}
        %     - x_1 \\
        %     - 2 x_2 \\
        %     0
        % \end{bmatrix} = 
        \begin{bmatrix}
            100 x_3 \\
            100 x_3 \\
            0
        \end{bmatrix}
    \end{equation}
    is over $100$ times larger in magnitude.
    By adding the velocity projection term to the loss function with a positive constant $\gamma$, we force the learned projection to account for the influence of $x_3$ on the dynamics. To substantiate our claims, we recreated the results presented in Holmes \textit{et al.} alongside an RVP loss trained ROM, where $P=\Phi \Psi^T$ and $\Psi^T \Phi = I$, and the results are shown in Figure~\ref{fig:linear_step_response}. 
    % \seoremark{The following sentence partly addresses comment~2 by reviewer~2.}
    The weight function (\ref{eqn:velocity_error_weight_function}) was used in (\ref{eqn:RVP_loss}) with $\gamma = \|C\|_{\text{op}} = 1$, $L = 1 / t_f$, and $t_f = 6$. 
    As expected, reconstruction loss (POD) performs poorly, while RVP loss performs nearly as well as balanced truncation.
\end{example}

RVP loss also resembles the loss function used to train SINDy-autoencoders\cite{Champion2019data}.
However, there are two main differences.
First, we determine the dynamics in the latent space via nonlinear projection using \eqref{eqn:ROM_in_latent_space}, whereas SINDy-autoencoders fit a model of the latent space dynamics during training.
Second, RVP loss measures the error between the ROM and the FOM time derivative projected onto the learned manifold.
In contrast, SINDy-autoencoders use a loss term measuring the difference between ROM and FOM time derivatives directly, i.e., without projection, together with another term measuring the difference between ROM and FOM time derivatives in the latent space.
While the RVP loss depends only on the projection $P$, the SINDy-autoencoder loss depends on the latent space, which can be scaled arbitrarily depending on the weights learned during training.
Using projected time derivatives to formulate RVP loss prevents the fast dynamics of the FOM from dominating the loss function, which can cause the learned manifold to become aligned with the fast dynamics, rather than capturing slow dynamics.
Incorporating our neural network architecture into SINDy-autoencoders where the latent space dynamics are learned is an interesting avenue of future work.
Variants of RVP loss could also be formulated in this setting and compared to the original SINDy-autoencoder loss.
We do not pursue this further here.

\subsection{Gradient-Aligned Projection (GAP) loss}
\label{subsec:GAP}

In order to quantify how well a given (nonlinear) projection $P:\R^n \to \R^n$ on the state space of a dynamical system preserves information about future outputs, we follow \citet{Otto2022model} and consider the map
\begin{equation}
    F_u: x_0 \mapsto (y(t_0), y(t_1), \ldots, y(t_L))
    \label{eqn:forecasting_map}
\end{equation}
defined by simulating the full-order model (\ref{eqn:FOM}) and sampling the output at times $0 \leq t_0 < t_1 < \cdots < t_L$.
We aim to find a projection so that $F_u(P(x))$ closely approximates $F_u(x)$ over a distribution of states $x$ and input signals $u$ drawn from trajectories of the full-order model.
If we are willing to simulate the FOM during the process of optimizing the projection, then we could form a loss function simply by computing the mean square error of these quantities.
However, this will be costly for high-dimensional systems of interest and we prefer a method that uses simulation data obtained from the FOM prior to optimizing the projection.

We construct a cost function that can be computed using a fixed set of samples from the FOM obtained ahead of time by expanding the difference $F_u(x) - F_u(P(x))$ in a first-order Taylor series about $x$.
Under mild boundedness and continuity assumptions, the following lemma says that we can use these first-order terms to bound the square error when $x - P(x)$ is small.
\begin{lemma}
    \label{lem:GAP_Taylor_estimate}
    Let $\mcal{X}$ be a compact convex subset of $\R^n$ and let $\mcal{U}$ be a compact topological space containing input signals $u$ defined on the interval $[0, t_L]$. 
    We assume that $(x, u) \mapsto F_u(x)$ is twice continuously differentiable with respect to $x$ on $\mcal{X}\times\mcal{U}$.
    Then there is a constant $C \geq 0$ so that
    \begin{equation}
        \big\Vert F_u(x) - F_u(P(x)) \big\Vert^2 \leq \big\Vert \D F_u(x) (x - P(x)) \big\Vert^2 + C \big\Vert x - P(x) \big\Vert^3
        \label{eqn:first_order_difference_estimate}
    \end{equation}
    holds whenever $x \in \mcal{X}$, $P(x) \in \mcal{X}$, and $u \in \mcal{U}$.
\end{lemma}
\begin{proof}
    This is a consequence of Taylor's theorem. We give the detailed proof in Appendix~\ref{app:proofs}.
\end{proof}
Taking the expectation over a distribution of states and input signals over sets satisfying the hypotheses of the lemma, the mean square approximation error is bounded by
\begin{multline}
    \E_{x, u} \Big[ \big\Vert F_u(x) - F_u(P(x)) \big\Vert^2\Big] 
    \leq \underbrace{\E_{x, u} \Big[\big\Vert \D F_u(x) (x - P(x)) \big\Vert^2\Big]}_{J_{\text{GAP}}(P)} \\
    + C \E_{x} \Big[\big\Vert x - P(x) \big\Vert^3\Big].
\end{multline}
We use a sample-based approximation of the leading-order term as a cost function for optimizing $P$ since, at least in principle, $\D F_u(x)$ can be computed prior to optimization given a collection of states and input signals.
For reasons that will become clear, we refer to this cost function as the gradient-aligned projection (GAP) loss.

In many practical applications the dimension $(L+1)m$ of the output sequences is large enough to make computing $\D F_u(x)$ impractical.
Instead, we can rely on randomized projections of the output sequences in a similar manner to the output projection method introduced by \citet{Rowley2005model}.
Specifically, we select an independent, zero mean, isotropic random vector $\xi \in \R^{(L+1)m}$ and compute the univariate gradients
\begin{equation}
    g = \grad (\xi^T F_u)(x)
    \label{eqn:random_univariate_gradient}
\end{equation}
using the adjoint of the full-order model linearized about the time-$t_L$ trajectory starting at $x$ as described in \citet{Otto2022model}.
These randomized univariate gradients allow us to write the GAP loss as
\begin{empheq}[box=\widefbox]{equation}
    J_{\text{GAP}}(P)
    % = \E_{x, u} \Big[\big\Vert \D F_u(x) (x - P(x)) \big\Vert^2\Big]
    = \E_{x, g} \Big[\big\langle g, \ x - P(x) \big\rangle^2\Big].
\end{empheq}
Collecting samples $\{(x_i, g_i)\}_{i=1}^s$ drawn from the joint distribution of $(x, g)$, we can compute a projection by minimizing the empirical GAP loss
\begin{equation}
    \hat{J}_{\text{GAP}}(P) = \frac{1}{s} \sum_{i=1}^s \big\langle g_i, \ x_i - P(x_i) \big\rangle^2.
    \label{eqn:GAP_loss}
\end{equation}

% \begin{proposition}
%     Suppose that the hypotheses of Lemma~\ref{lem:GAP_Taylor_estimate} hold and let $\D F_u(x)$ be injective for every $(x,u) \in \mcal{X}\times \mcal{U}$.
%     Letting $\sigma_n \big( \D F_u(x) \big)$ denote the $n$th singular value of $\D F_u(x)$, the constant
%     \begin{equation}
%         \sigma = \min_{(x,u) \in \mcal{X}\times \mcal{U}} \sigma_n \big( \D F_u(x) \big)
%     \end{equation}
%     is positive.
%     Letting $D = \sup$
%     \begin{equation}
%         \E_{x, u} \Big[ \big\Vert F_u(x) - F_u(P(x)) \big\Vert^2\Big]
%         \leq \E_{x, g} \Big[\big\langle g, \ x - P(x) \big\rangle^2\Big] + C \E_{x} \Big[\big\Vert x - P(x) \big\Vert^3\Big]
%     \end{equation}
% \end{proposition}

% \begin{equation}
%     \E_{x} \Big[\big\Vert x - P(x) \big\Vert^3\Big] \leq \frac{1}{\sigma^3} \E_{x, g} \Big[\big\vert\big\langle g, \ x - P(x) \big\rangle\big\vert^3\Big]
% \end{equation}

Minimizing GAP loss over linear projections for linear time-invariant (LTI) systems becomes equivalent to balanced truncation (BT) for certain limits and distributions of $x$.
For example, let
\begin{equation}
\begin{aligned}
    \dot{x} &= A x + B u \\
    y &= C x
\end{aligned}
\label{eqn:LTI_system}
\end{equation}
be an asymptotically stable LTI system with $\dim(u) = d_u$.
Suppose we sample $x$ uniformly from impulse response trajectories with $0 \leq t \leq t_f$ and initial conditions $x(0) = B e_j$, $j=1, \ldots, d_u$.
If we choose uniformly spaced sample times $t_k = k t_f / L$ to form $F_u$, then it is straightforward to show that
\begin{equation}
    \lim_{t_f \to \infty} \lim_{L\to\infty} \frac{t_f^2 d_u}{L} J_{\text{GAP}}(P)
    = \Tr\left[ W_{o} (I-P) W_{c} (I-P)^T \right]
\end{equation}
where $W_{o}$ and $W_{c}$ are the observability and controllability Gramians of (\ref{eqn:LTI_system}).
The quantity on the right is minimized by the balanced truncation projection \cite{Otto2022model, Singler2010optimality, Singler2015optimality}.
The performance of BT on the nonnormal system in Example~\ref{ex:linear} is shown in Figure~\ref{fig:linear_step_response}, 
providing evidence that minimizing GAP loss is appropriate for modeling such systems.
% This provides evidence that minimizing GAP loss is appropriate for nonnormal systems.

\subsection{Orthogonality-promoting regularization}
Regularization is often employed in over-parametrized neural networks to prevent over-fitting.
A commonly used method is to penalize the squared Frobenius norm of the weights in each layer.
It turns out that applying this penalty to weights $(\Phi,\Psi)$ in the biorthogonal manifold $\mcal{B}_{n,r}$ drives them towards orthogonality, that is $\Phi = \Psi$ having orthonormal columns.
Specifically, we have the following result:
\begin{theorem}
    \label{thm:orthogonality_promoting_regularization}
    The minimum value of the function
    \begin{equation}
        R_F(\Phi, \Psi) = \Vert \Phi \Vert_F^2 + \Vert \Psi \Vert_F^2
    \end{equation}
    over $(\Phi, \Psi) \in \mcal{B}_{n,r}$ is $2r$ and this value is achieved if and only if $\Phi = \Psi$ has orthonormal columns.
    Moreover, $R_F(\Phi, \Psi) \to \infty$ if any of the principal angles\cite{Afriat1957orthogonal, Bjorck1973numerical} between $\Range(\Phi)$ and $\Range(\Psi)$ approach $\pi/2$.
\end{theorem}
\begin{proof}
    We give a proof in Appendix~\ref{app:proofs}.
\end{proof}

\section{Case study of a simplified fluid model}
\label{sec:noack_case_study}

In this section, we compare our reduced-order modeling approach to several other methods, on a highly simplified model of a fluid flow.
The studied dynamical system is a three-state model of vortex shedding behind a circular cylinder, as described by Noack \textit{et al.}\cite{Noack2003hierarchy}.
In particular, the system is defined by the following set of equations:
\begin{equation}
    \begin{aligned}
        \dot{x}_1 = \mu x_1 - \omega x_2 + A x_1 x_3 \\
        \dot{x}_2 = \omega x_1 + \mu x_2 + A x_2 x_3 \\
        \varepsilon \dot{x}_3 = - (x_3 - x_1^2 - x_2^2)
        \label{eq:noack}
    \end{aligned}
\end{equation}
with $\mu = 0.1$, $\omega = 1$, $A = -0.1$, and $\varepsilon = 0.1$.
This system possesses an unstable fixed point at the origin and a global asymptotically stable limit cycle of radius $1$ about $x_1=x_2=0$ in the plane $x_3=1$. 
Additionally, the system's slow manifold is situated a distance $\mathcal{O}(\varepsilon)$ away from the critical manifold $x_3 = x_1^2 + x_2^2$.
An asymptotic approximation of the slow manifold to second order in $\varepsilon$ can be found in Otto's thesis \cite{Otto2022advances} along with the recurrence relation needed to obtain the higher-order terms.
In this work we use the fourth-order approximation in~$\varepsilon$ computed using this relation. 
The slow manifold calculation follows the same procedure discussed in Example~\ref{ex:toy_slow_fast}. 
In the following section, we outline the network architectures responsible for learning the nonlinear projection described in Section~\ref{sec:nlprom}.

\subsection{Autoencoder Architectures}
\label{subsec:autoencoder_architectures}

We compare two autoencoder architectures in this manuscript.
The first architecture, which we refer to as ProjAE, is the projection-constrained autoencoder described in Section~\ref{sec:architecture}.
The second architecture, referred to as StandAE, is a standard state-of-the-art differentiable autoencoder.
In particular, this architecture's encoder and decoder are modeled as fully connected networks using the GeLU activation function, denoted~$\sigma$ \cite{hendrycks2016gaussian}, which satisfies the requirement of differentiability discussed in Section~\ref{sec:nlprom}.
The encoder and decoder layer structure, denoted as $\psi_e^{(l)}$ and $\psi_d^{(l)}$, follow the standard feed-forward neural network form: $\sigma(A^{(l)} x + b^{(l)})$\cite{Goodfellow2016deep}.
As a final note, we attached a linear output layer to both the encoder and decoder of StandAE, i.e., $\psi_e = W_e \psi_e^{(L)} \circ \cdots \circ \psi_e^{(1)}$ and $\psi_d = W_d \psi_d^{(L)} \circ \cdots \circ \psi_d^{(1)}$ where $W_e$ and $W_d$ are trainable weight matrices.

To initialize the weights and biases of ProjAE, we follow the procedure outlined in Section~\ref{subsec:initialization}.
StandAE was initialized using the procedure discussed in Section 2.2 of Hendrycks \textit{et al.}\cite{hendrycks2016adjusting}.
In particular, the rows of each weight matrix were uniformly sampled from unit hypersphere. StandAE's weights were then scaled by a GeLU dependent factor designed too maintain both activation and back-propagated gradients variance as one forward or backward through the network.
For both architectures, the biases are set to zero at initialization.
Both architectures have a $5$ layer encoder and $5$ layer decoder where $\psi_e^{(5)}: \R^3 \rightarrow \R^2$, $\psi_d^{(5)}: \R^2 \rightarrow \R^3$, and $\psi_e^{(i)}, \psi_d^{(i)} : \R^3 \rightarrow \R^3$ for $i = 1, 2, 3, 4$.
We do not use the constraint described in Section~\ref{subsec:preserving_an_equilibrium} to preserve the equilibrium at the origin.

\subsection{Autoencoder-Based Reduced-Order Models}
\label{subsec:autoencoder_based_ROMs}

When defining the reduced-order model, we must select a method by which we project the dynamics onto the learned manifold.
As discussed in Section~\ref{sec:architecture}, one approach is to use the encoder to define the reduced-order model~\eqref{eqn:ROM_in_latent_space}.
We denote this type of reduced-order model by EncROM.
The approach used by \citet{Lee2020model} has instead projected the dynamics orthogonally onto the tangent space of the manifold parameterized by the decoder.
We denote this type of reduced-order model by DecROM.

\subsection{Data Collection}
\label{subsec:noack_data_collection}

Two separate data sets were generated to examine the effect on training.
The first data set, which we call the Fine Data Set, consisted of 1000 trajectories with initial conditions given on a $10\times 10\times 10$ grid evenly spaced in the cube $[-1,1]^3$.  The second data set, which we call the Coarse Data Set, consisted of 216 trajectories with initial conditions $\big\{-1, -\tfrac{1}{3}, -\tfrac{1}{9}, \tfrac{1}{9}, \tfrac{1}{3}, 1\big\}^3$.  For each training set, we created a validation data set with same number of trajectories, with initial conditions sampled uniformly from the cube $[-1,1]^3$. The testing data set consisted of 1000 trajectories with initial conditions sampled uniformly in the cube.

Trajectories were generated by numerically integrating the governing equations with a 4th-order Runge-Kutta method, over the time interval $[0,20]$, using a time step $\Delta t = 0.1$.
In order to generate the gradient samples for GAP loss, we used the method of long trajectories discussed by \citet{Otto2022model} with parameters $s_g=10$ and~$L=20$.
The hyperparameter $L$ was chosen such that if a gradient sample was based at the initial condition, then the adjoint would be sampled before and after transients have decayed.
In this example, transients decay after about $0.2$ time units and trajectories reach the limit cycle by about $2$ time units.
Using the aforementioned parameters, the fine data set has a total of $201,000$ state samples and $191,558$ gradient samples, and the coarse data set has a total of $43,416$ state samples and $41,307$ gradient samples.
The hyperparameter $s_g$ was chosen such that the number of state and gradient samples were roughly the same to give all loss functions a fair chance to perform.

\subsection{Training Procedure}
\label{subsec:noack_train_procedure}

In total, $12$ training sessions were carried out, with each session corresponding to a unique combination of data set, architecture, and loss function.
During each session, $64$ networks were trained simultaneously, each with a different choice of initial parameters (weights and biases).
All 12 training sessions used the same 64 sets of initial parameters.
The weights and biases of each network were saved during training if the lowest loss-function evaluation on the validation data set was achieved. 
Due to the computational cost of simulating the autoencoder-based reduced-order model, we used the loss function to determine which model to save, instead of simulating the reduced-order model explicitly.
The computational challenge of simulating the reduced-order model is addressed in Section~\ref{sec:ROM_assembly}.
After each session's training phase, the most effective EncROM and DecROM models were chosen from the saved networks, and this selection process was based on the true ROM prediction error (rather than the loss function), using the fine or coarse validation data sets.

To ensure a fair comparison across network architectures and loss functions, each training session employed mostly identical hyperparameters.
All training sessions implemented the PyTorch ReduceLROnPlateau class with a patience of 50, an initial learning rate of $10^{-3}$, and a validation loss equal to the loss-function evaluation on the validation data set. 
Using PyTorch's built-in Adam\cite{Kingma2014Adam} optimizer with default settings, each network was trained for a total of $900$ epochs.

For reconstruction loss and GAP loss, a batch size of $400$ was employed.
In the case of RVP loss, we utilized a prediction horizon of $t_f=20$ and a trajectory batch size of $2$ (with a time step $\Delta t=0.1$, as mentioned previously), so that each mini-batch looks at the same number of sample points.
Trapezoidal integration was used to discretize the integral in \eqref{eqn:RVP_loss} defining the RVP loss.
% \seoremark{The above sentence addresses comment~3 by reviewer~1.}
Since the full state is being observed, we set $\gamma = 1$ per the discussion in Section~\ref{subsec:RVP}.
We use $L = 1 / t_f$ to define the weight function in \eqref{eqn:velocity_error_weight_function}.
Finally, all ProjAEs were trained using the regularization
in \eqref{eqn:regularization} with a factor $\beta = 10^{-5}$.
% for ProjAE, all loss functions incorporated the regularization in \eqref{eqn:regularization} with a factor $\beta = 10^{-5}$.
% \seoremark{The above partly addresses comment~2 by reviewer~2.}

\subsection{Results}
\label{subsec:noack_results}

We expect a successful reduced-order model to learn and capture three fundamental features of this example.
First, the autoencoder’s range should closely approximate the system’s slow manifold.
Second, the projected dynamics should approximate the dynamics on the slow manifold.
Finally, the fibers of projection should align with the direction of fast dynamic transients.

In order to quantitatively analyze these features, we define two performance metrics.
The first metric measures the proximity between the autoencoder’s range and the slow manifold.
In particular, \textit{manifold reconstruction error} is defined by
\begin{equation}
    \label{eq:manif_recon_error}
    \frac{1}{|\tilde{\mathcal{M}}|}\sum_{x \in \tilde{\mathcal{M}}}\|x - P(x)\|_2^2,
\end{equation}

\noindent
where $\tilde{\mathcal{M}}$ is a finite subset of the system's slow manifold.
In this study, $\tilde{\mathcal{M}}$ contains points of the form $(x_1, x_2, h_\epsilon(x_1, x_2))$ where $h_\epsilon(x_1, x_2)$ is the slow manifold's graph representation to fourth order.
Furthermore, coordinates $(x_1, x_2)$ were sampled on a $20\times 20$ grid evenly spaced in the square $[-1, 1]^2$.
The second metric, called \textit{ROM prediction error}, quantifies a ROM's ability to predict an initial condition's future, and is defined by
\begin{equation}
    \label{eq:rom_pred_error}
    \frac{1}{N} \sum_{n=1}^N \|\hat{x}_n - x_n\|_2^2,
\end{equation}
where $x_n$ corresponds to a state sample from either the validation or test data set and $\hat{x}_n$ denotes the corresponding state predicted by the ROM.
Note that the prediction error above depends on both the autoencoder $P$ and the chosen method of projection, EncROM or DecROM.
Employing these metrics, alongside other qualitative techniques, let us now examine how the various methods presented here perform, relative to existing methods.

First, let us explore how closely the autoencoder's range approximates system's slow manifold.
Looking at Table~\ref{tab:noack_res}, we find that all models trained on reconstruction loss are able to consistently capture the slow manifold, with a manifold reconstruction error of at most $0.005$.
The majority of models trained on GAP loss also have small manifold reconstruction error.
We observe a large manifold reconstruction error for RVP loss, possibly because the models were selected on forecasting, and not reconstruction.

Next, let us explore how well the reduced-order models make forecasts from new initial conditions, as quantified by the ROM prediction error~\eqref{eq:rom_pred_error}. 
As shown in the ``Pred'' columns of Table~\ref{tab:noack_res}, we find that for all 12 autoencoder-loss combinations, the EncROM models exhibit a lower prediction error than the DecROM models.
Therefore, at least for this example, an encoder-based ROM provides a benefit over the traditionally-used decoder-based ROM.
This effect is more pronounced for our new cost functions, GAP and RVP.
Furthermore, some of the DecROM models blow up or have very large error.  
This is because of an effect we observed in Example~\ref{ex:toy_slow_fast} in Section~\ref{sec:nlprom}: in particular, orthogonally projecting onto the tangent space of an approximate manifold can yield incorrect stability types for fixed points and periodic orbits.

For prediction of dynamics, the traditional Reconstruction loss performs poorly across the board, for reasons we have explained in Example~\ref{ex:toy_slow_fast}.

The lowest error was obtained for RVP loss, with EncROM projection, and ProjAE architecture, with GAP loss having similar results.  The degree to which the constraints imposed by the ProjAE architecture are beneficial depend both on the cost function and the size of the training data set.
Enforcing constraints significantly improved performance when training with RVP loss, and this benefit was more pronounced when the size of training data set was smaller.
% When training with GAP loss on the larger data set, enforcing the
When training with GAP loss on a large data set, the standard autoencoder was able to achieve high forecasting accuracy without additional constraints.
These constraints were beneficial when training with GAP loss on a smaller data set.

\begin{figure*}
    \centering
    \includegraphics[scale=0.7]{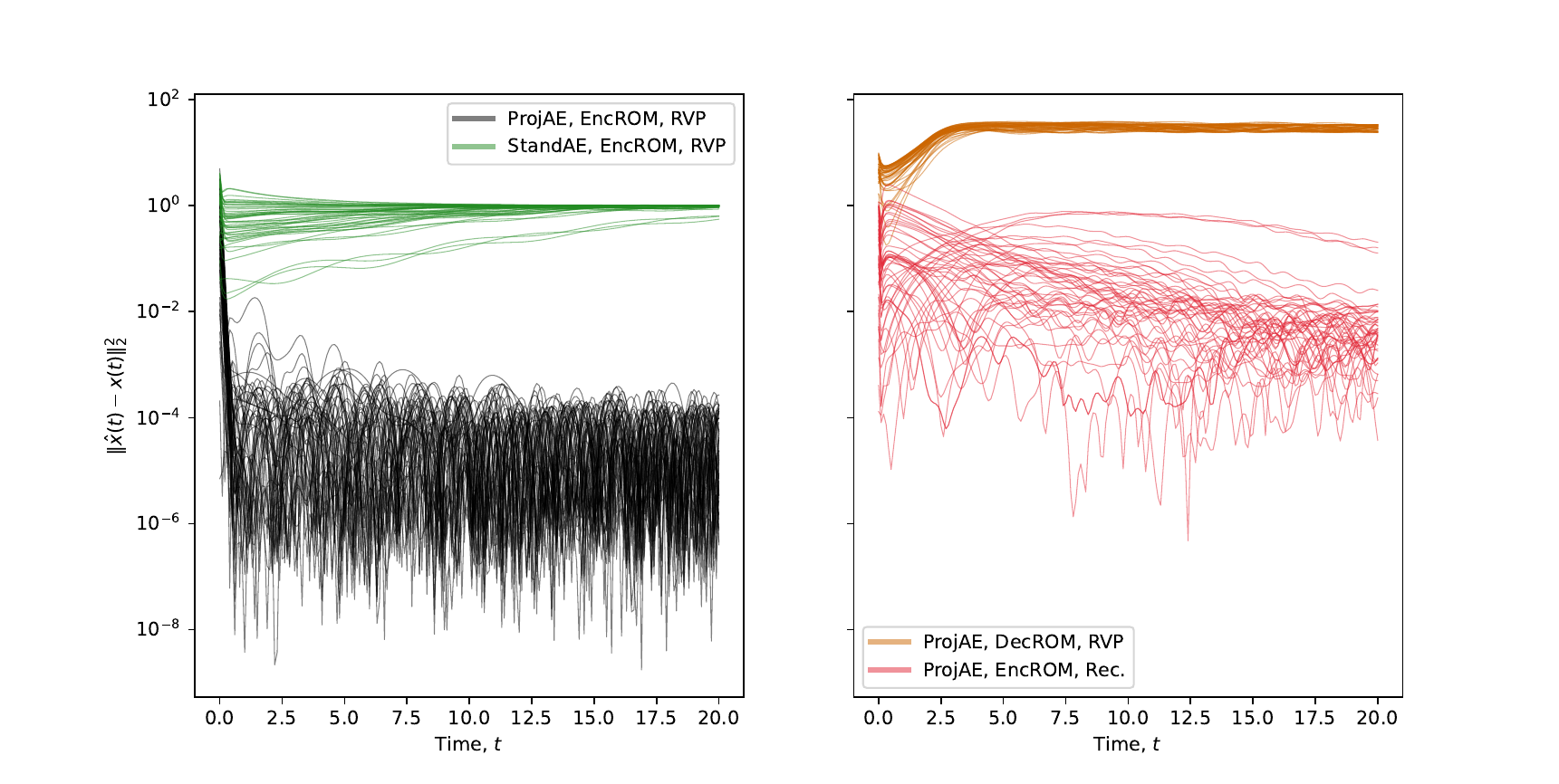}
    \caption{Comparative study of the model that performs best for prediction (ProjAE architecture, RVP loss, and EncROM projection, trained on the fine data set), changing one component at a time.
    \emph{Left:} ROM prediction error for 50 test trajectories, comparing the two network architectures (ProjAE and StandAE).
    Our projection-constrained autoencoder reduces average prediction error by three orders of magnitude.
    \emph{Right:} ROM prediction error for the typical projection approach (DecROM) and loss function (Reconstruction). In both cases, our approach significantly decreases the model prediction error.}
    \label{fig:noack_error}
\end{figure*}

These observations are illustrated further in Figure~\ref{fig:noack_error}, which shows the error on all 50 test trajectories.  We start with our best architecture (ProjAE architecture, with RVP loss and EncROM projection), and change one component at a time.

\begin{figure*}
  \centering
  \begin{tabular}{c}
    \subfloat[(ProjAE, EncROM, Rec.) Model\label{fig:noack_RECtraj}]{%
      \includegraphics[scale=0.515]{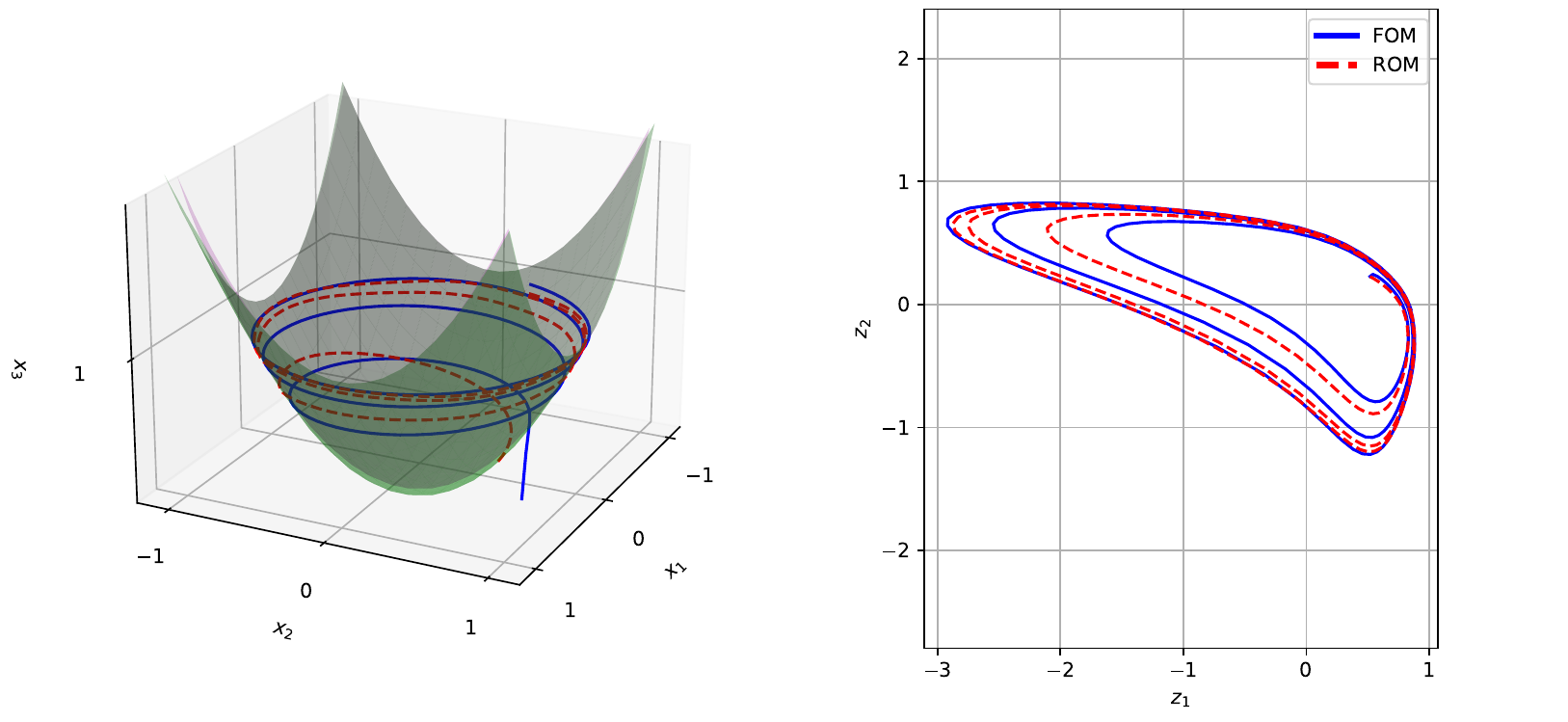}%
    } \\
    \subfloat[(ProjAE, EncROM, GAP) Model\label{fig:noack_GAPtraj}]{%
      \includegraphics[scale=0.515]{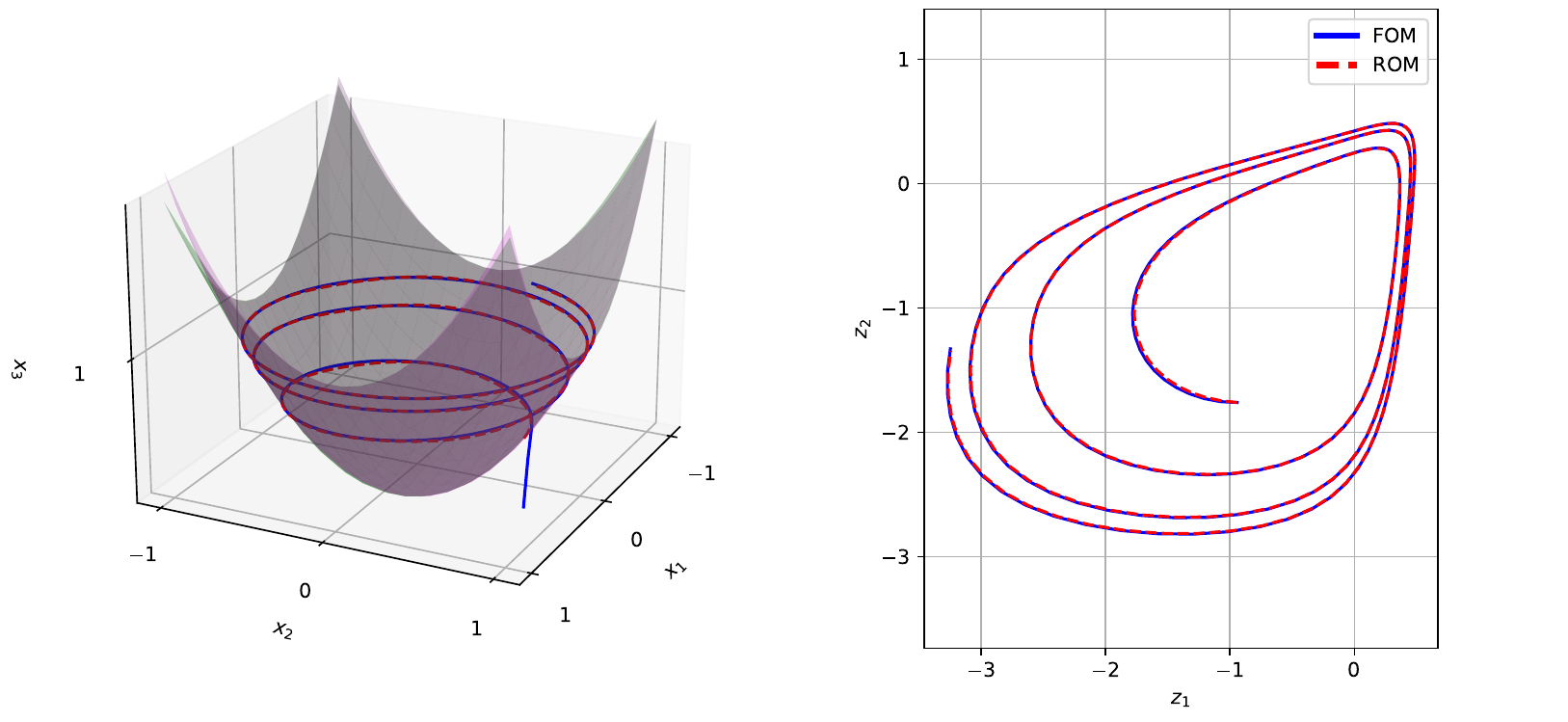}%
    }
    \\
    \subfloat[(ProjAE, EncROM, RVP) Model\label{fig:noack_RVPtraj}]{%
      \includegraphics[scale=0.515]{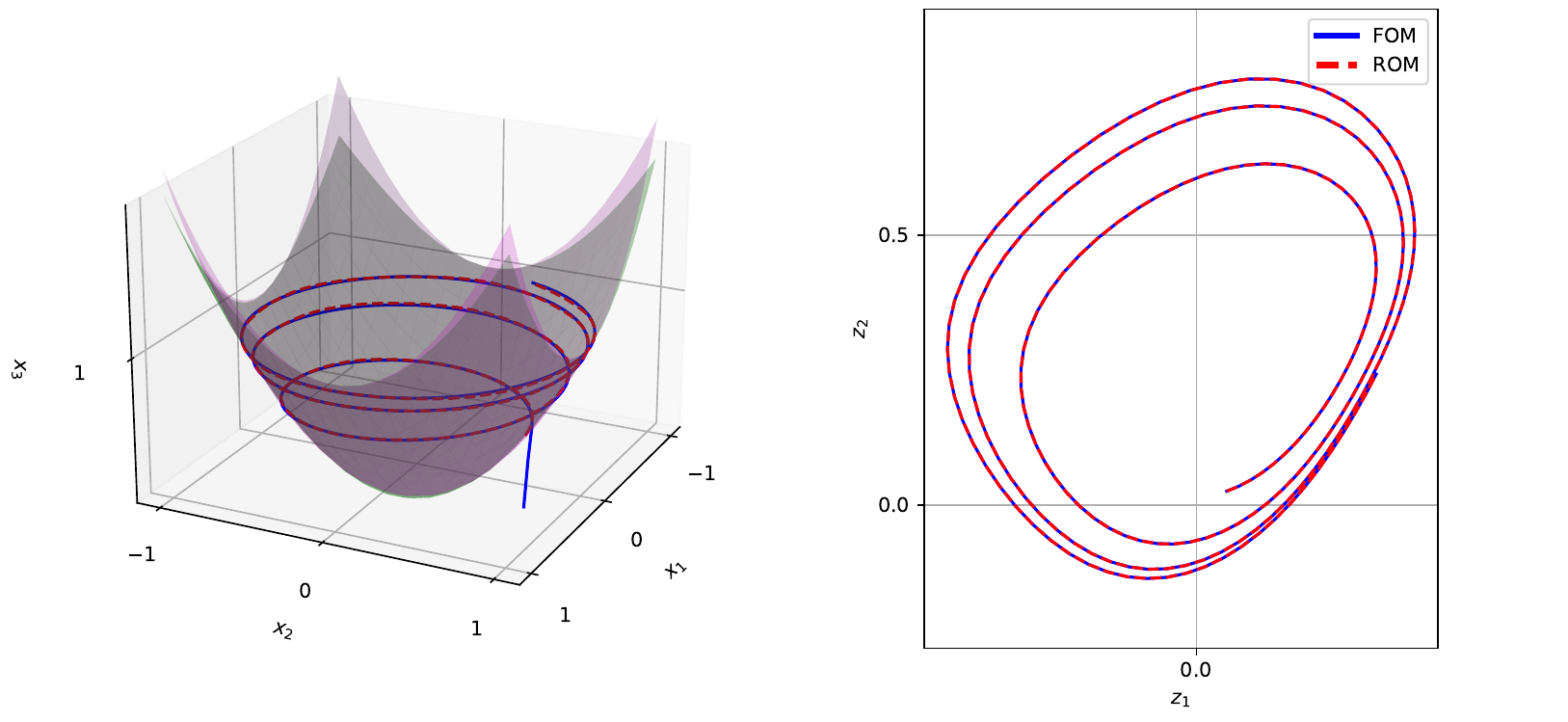}%
    }
  \end{tabular}
  \caption{Visualizing a sample test trajectory for each loss function. Models seen here are trained on the Fine Data Set.
  \emph{Left:} The slow and learned manifold (shaded violet and green, respectively), along with a trajectory from the full and reduced models (in blue and red); \emph{Right:} Corresponding trajectories in the latent space. In the top left figure, the initial condition is projected almost orthogonally onto the learned manifold, while in the other cases, the initial condition is projected vertically along the direction of fast dynamics.}
  \label{fig:noack_traj}
\end{figure*}

Figure~\ref{fig:noack_traj} shows a typical test trajectory in both the 3-dimensional state space, as well as the latent space, for ProjAE architecture and EncROM projection, comparing the three loss functions (Reconstruction, GAP, and RVP).  When reconstruction loss is used, the projection approximates an orthogonal projection, while the other loss functions result in oblique projection, accounting for the fast dynamics.

\begin{table*}
  \centering
  \subfloat[Fine Data Set]{%
\begin{tabular}{|c|cc|cc|cc|}
\hline
                & \multicolumn{2}{c|}{Rec.}              & \multicolumn{2}{c|}{GAP}               & \multicolumn{2}{c|}{RVP}               \\ \hline
                & \multicolumn{1}{c|}{Manif.}  & Pred.   & \multicolumn{1}{c|}{Manif.}  & Pred.   & \multicolumn{1}{c|}{Manif.}  & Pred.   \\ \hline
ProjAE, EncROM  & \multicolumn{1}{c|}{0.00031} & 0.08339 & \multicolumn{1}{c|}{0.00030} & 0.00584 & \multicolumn{1}{c|}{0.00218} & 0.00566 \\ \hline
StandAE, EncROM & \multicolumn{1}{c|}{0.00469} & 0.07078 & \multicolumn{1}{c|}{0.00196} & 0.00646 & \multicolumn{1}{c|}{0.78596} & 0.84948 \\ \hline
ProjAE, DecROM  & \multicolumn{1}{c|}{0.00031} & 0.08565 & \multicolumn{1}{c|}{0.00019} & 0.02084 & \multicolumn{1}{c|}{3.67335} & 27.1888 \\ \hline
StandAE, DecROM & \multicolumn{1}{c|}{0.00068} & 0.08484 & \multicolumn{1}{c|}{0.00165} & $\infty$ & \multicolumn{1}{c|}{0.79029} & $\infty$ \\ \hline
\end{tabular}
    \label{tab:noack_res10}
  } \\
  \subfloat[Coarse Data Set]{%
\begin{tabular}{|c|cc|cc|cc|}
\hline
                & \multicolumn{2}{c|}{Rec.}              & \multicolumn{2}{c|}{GAP}               & \multicolumn{2}{c|}{RVP}                \\ \hline
                & \multicolumn{1}{c|}{Manif.}  & Pred.   & \multicolumn{1}{c|}{Manif.}  & Pred.   & \multicolumn{1}{c|}{Manif.}  & Pred.    \\ \hline
ProjAE, EncROM  & \multicolumn{1}{c|}{0.00091} & 0.05780 & \multicolumn{1}{c|}{0.00023} & 0.00658 & \multicolumn{1}{c|}{0.00281} & 0.00570 \\ \hline
StandAE, EncROM & \multicolumn{1}{c|}{0.00201} & 0.07553 & \multicolumn{1}{c|}{0.00390} & 0.01381 & \multicolumn{1}{c|}{0.78988} & 0.84743  \\ \hline
ProjAE, DecROM  & \multicolumn{1}{c|}{0.00091} & 0.06229 & \multicolumn{1}{c|}{0.00024} & $\infty$ & \multicolumn{1}{c|}{0.00109} & 0.3701  \\ \hline
StandAE, DecROM & \multicolumn{1}{c|}{0.00201} & 0.12626 & \multicolumn{1}{c|}{0.22420} & $\infty$ & \multicolumn{1}{c|}{0.80293} & $\infty$ \\ \hline
\end{tabular}
\label{tab:noack_res4}
}
\caption{Manifold reconstruction error \eqref{eq:manif_recon_error} and ROM prediction error \eqref{eq:rom_pred_error} for the 12 training sessions.  Each row corresponds to an autoencoder-ROM combination described in Sections~\ref{subsec:autoencoder_architectures}~and~\ref{subsec:autoencoder_based_ROMs}.
Each column identifies the type of loss function used during training (Reconstruction, GAP, or RVP) as well as the performance metric (Manifold or Prediction).
The lowest prediction error is achieved by the ProjAE network, with RVP loss, and with EncROM projection.}
\label{tab:noack_res}
\end{table*}

\section{Assembling efficient ROMs}
\label{sec:ROM_assembly}
The example discussed in the previous section began with a system that was already low dimensional, with only 3 states.  For higher dimensional systems, significant computational challenges arise when simulating the reduced-order model.  In this section, we discuss three possible methods for addressing these challenges.

After training the autoencoder, we obtain a nonlinear projection-based reduced-order model \eqref{eqn:ROM_in_latent_space} in the autoencoder's latent space coordinates.
Even though the latent space is low-dimensional,
evaluating the right-hand side, $\tilde{f}$, of (\ref{eqn:ROM_in_latent_space}) involves evaluating the right-hand-side, $f$, of the full-order model (\ref{eqn:FOM}).
Hence, we cannot expect speedups when simulating the ROM by evaluating $\tilde{f}$ in this manner.
This section presents three methods for obtaining computationally-efficient ROMs that can be evaluated more quickly than the FOM.
However, even in cases when it is more costly to evaluate the ROM than the FOM, we note that it may be possible to use larger time steps when simulating the ROM due to the removal of dynamics with fast time scales.

\subsection{Fitting the model in latent space}
\label{subsec:model_fitting}

A simple approach to construct an efficient ROM in the latent space is to fit a surrogate model for $\tilde{f}$ and $\tilde{g}$ in (\ref{eqn:ROM_in_latent_space}) using sample-based interpolation or regression.
Specifically, given a collection of samples $z_i\in \R^r$ in the latent space of the autoencoder and samples of the input $u_i$, we can evaluate $\dot{z}_i = \tilde{f}(z_i, u_i)$ and $y_i = \tilde{g}(z_i)$ using the definitions in (\ref{eqn:ROM_in_latent_space}), which rely on the FOM.
Once the time derivatives and outputs at the samples have been evaluated, we can fit surrogates for $\tilde{f}$ and $\tilde{g}$ that can be evaluated more efficiently.
Since $\tilde{f}$ and $\tilde{g}$ can be evaluated at arbitrary pairs $(z,u)$, 
we can choose the samples $(z_i, u_i)$ to achieve a desired level of accuracy for the surrogates of $\tilde{f}$ and $\tilde{g}$.
In particular, we are not limited to the encoded snapshots used to train the autoencoder.

Appropriate sampling and fitting procedures to construct the surrogates will depend on the dimension of the latent space.
For very low-dimensional latent spaces ($\leq 5$-dimensional) it is possible to construct a grid of sample locations and use spline-based interpolation.
For higher-dimensional latent spaces, one can rely on random sampling and radial basis function interpolation or Gaussian process regression.
The distribution from which the samples are drawn can be based on a density estimate from the encoded snapshot data collected from the FOM.
More samples can also be added in an iterative manner until a desired level of accuracy for the surrogates of $\tilde{f}$ and $\tilde{g}$ is achieved.

\subsection{Assembling tensors using the outer layer}
\label{subsec:tensors}
In certain cases when the full-order model has polynomial nonlinearities, we can improve the efficiency of the reduced-order model by pre-computing the linear projection of the full-order model associated with the outer-most layer of the autoencoder.
Our method is similar to the approach described in Section~4.2 of \citet{Holmes2012turbulence} for assembling Petrov-Galerkin models.
To illustrate, suppose that the right-hand side of the full-order model (\ref{eqn:FOM}) has a term $f_2:\R^n \to \R^n$ that can be expressed as
\begin{equation}
    f_2(x) = h_2(x, x),
\end{equation}
where $h_2:\R^n \times \R^n \to \R^n$ is a symmetric bilinear form.
In the Navier-Stokes equations, such terms arise from discretization of the convective term $u\cdot \nabla u$ and from the solution of the pressure-Poisson equation $\grad P = - \nabla \Delta^{-1} \nabla \cdot (u\cdot \nabla u)$, where $u$ denotes the velocity field.
More generally, a system with polynomial nonlinearities can always be
converted into a system with quadratic nonlinearities evolving on an invariant submanifold in a higher-dimensional state space via a lifting process called ``quadratization'' \cite{Bychkov2023exact,Gu2011qlmor,Kramer2022balanced}.
The quadratic nonlinearity can then be expressed using a symmetric bilinear form.
% \begin{remark}
%     Since the original system corresponds to an invariant submanifold of the lifted system, our model reduction approach may prove useful after a system has been quadratized.
% \end{remark}

Isolating the linear operations in the outer-most layer of the autoencoder, we observe that the encoder and decoder can be written as
\begin{equation}
    \psi_e: x \mapsto \tilde{\psi}_e\big(\Psi_L^T(x - b_L) \big), \qquad \psi_d: z \mapsto b_L + \Phi_L \tilde{\psi}_d(z),
\end{equation}
where 
% $\tilde{\psi}_{e} = \psi_{e}^{(1)} \circ \cdots \circ \psi_{e}^{(L-1)}\circ \sigma_{-}$ and $\tilde{\psi}_{d} = \sigma_{+}\circ\psi_{d}^{(L)} \circ \cdots \circ \psi_{d}^{(1)}$.
\begin{equation*}
    \tilde{\psi}_{e} = \psi_{e}^{(1)} \circ \cdots \circ \psi_{e}^{(L-1)}\circ \sigma_{-}, \qquad
    \tilde{\psi}_{d} = \sigma_{+}\circ\psi_{d}^{(L-1)} \circ \cdots \circ \psi_{d}^{(1)}.
\end{equation*}
This allows us to express the ROM given by (\ref{eqn:ROM_in_latent_space}) as
\begin{equation}
\begin{aligned}
    \tilde{f}(z,u) &= \D \tilde{\psi}_e\big( \tilde{\psi}_d(z) \big) \Psi_L^T f\big( b_L + \Phi_L \tilde{\psi}_d(z), \ u \big)  \\
    \tilde{g}(z) &= g\big(b_L + \Phi_L \tilde{\psi}_d(z)\big).
\end{aligned}
\label{eqn:outer_layer_ROM}
\end{equation}
The contribution of the bilinear term to the ROM expressed element-wise is given by
% \begin{multline}
%     \big[\Psi_L^T f_2\big(b_L + \Phi_L \tilde{\psi}_d(z)\big) \big]_i
%     = \underbrace{\col_i[\Psi_L]^T h_2\big( b_L, b_L \big)}_{a_i} \\
%         + 2\sum_{j_1 = 1}^{n_{L-1}} \underbrace{\col_i[\Psi_L]^T h_2\big(b_L, \col_{j_1}[\Phi_L] \big)}_{b_{i, j_1}} \big[\tilde{\psi}_d(z)\big]_{j_1} \\
%         + \sum_{j_1 = 1}^{n_{L-1}} \sum_{j_2 = 1}^{n_{L-1}} \underbrace{\col_i[\Psi_L]^T h_2\big(\col_{j_1}[\Phi_L], \col_{j_2}[\Phi_L] \big)}_{c_{i, j_1, j_2}} \big[\tilde{\psi}_d(z)\big]_{j_1} \big[\tilde{\psi}_d(z)\big]_{j_2}.
% \end{multline}
\begin{multline}
    \big[\Psi_L^T f_2\big(b_L + \Phi_L \tilde{\psi}_d(z)\big) \big]_i
    = \underbrace{\Psi_L[:,i]^T h_2\big( b_L, b_L \big)}_{a_i} + \\
    2\sum_{j_1 = 1}^{n_{L-1}} \underbrace{\Psi_L[:,i]^T h_2\big(b_L, \Phi_L[:,j_1] \big)}_{b_{i, j_1}} \big[\tilde{\psi}_d(z)\big]_{j_1} + \\
    \sum_{j_1,j_2 = 1}^{n_{L-1}} \underbrace{\Psi_L[:,i]^T h_2\big(\Phi_L[:,j_1], \Phi_L[:,j_2] \big)}_{c_{i, j_1, j_2}} \big[\tilde{\psi}_d(z)\big]_{j_1} \big[\tilde{\psi}_d(z)\big]_{j_2}.
\end{multline}
We observe that the elements of the tensors $[a_i]$, $[b_{i, j_1}]$, and $[c_{i,j_1, j_2}]$ can be computed and stored prior to simulating the ROM.

Even if one does not employ quadratization, the above approach applies analogously to any term in the governing equations that can be expressed as a sum of multilinear forms, that is, any polynomial term of finite degree.
The rank of the tensors to be assembled is $d+1$ where $d$ is the degree of the polynomial nonlinearity.
The dimensions of these tensors are all equal to the layer width $n_{L-1}$.
In general, these tensors are dense.
Therefore, the amount of storage and number of operations required to act with the pre-assembled tensors on $\tilde{\psi}_d(z)\in\R^{n_{L-1}}$ both scale as $\mcal{O}\big((n_{L-1})^{d+1}\big)$.
This should be compared against the $\mcal{O}(n)$ scaling typically required to act with sparse finite difference operators of the FOM acting on $\psi_d(z)\in\R^{n}$.
Therefore, simulating a ROM based on pre-assembled tensors will likely be advantageous only when $(n_{L-1})^{d+1} \ll n$, that is, when the degree of the polynomial nonlinearity and the width $n_{L-1}$ are both sufficiently small.
For example, when the FOM has $n=10^5$ state variables coming from a finite difference discretization of the incompressible Navier-Stokes equations ($d=2$), we only expect to see advantages from pre-assembling tensors when $n_{L-1}$ is in the low tens.

Since the decoder reconstructs states in an affine subspace of dimension $n_{L+1}$, making this parameter too small can impair the decoder's ability to reconstruct state data with slowly decaying Kolmogorov $n$-widths (see Remark~\ref{rem:outer_layer_width}).
The trade-off between computational efficiency and representational power associated with the choice of $n_{L-1}$ limits the scope of applications in which pre-assembling tensors will be advatageous for reduced-order modeling.

\subsection{Sparsifying the encoder}
\label{subsec:sparse_encoder}

When the state variables in the governing equations of the full-order order model are sparsely coupled, computational speedups for the reduced-order model can be achieved by sparsifying the weight matrices in the encoder.
Here we rely on essentially the same principle as the Discrete Empirical Interpolation Method (DEIM) \cite{Chaturantabut2010nonlinear}.
That is, if the time-derivative of the reduced-order model can be determined based on the time derivatives of a small collection of state variables in the full-order model, then we need only reconstruct the neighboring variables to evolve the reduced-order model.

Given a collection of state variable indices $\mscr{I} = \{ i_1, \ldots, i_K \}$, we defined the selection operator $S_{\mscr{I}}:\R^n \to \R^K$ by
\begin{equation}
    S_{\mscr{I}}:\big( [x]_1, \ldots, [x]_n \big) \mapsto
    \big( [x]_{i_1}, \ldots, [x]_{i_K} \big).
\end{equation}
The time derivative of the selected states $S_{\mscr{I}} x$ under (\ref{eqn:FOM}) depend on a collection of state variables with indices $\mscr{N}(\mscr{I})$ that we refer to as the ``neighbors'' of $\mscr{I}$.
In other words, there is a function $f_{\mscr{I}}$ so that
\begin{equation}
    \ddt (S_{\mscr{I}} x) 
    = S_{\mscr{I}} f(x,u) 
    = f_{\mscr{I}}\big( S_{\mscr{N}(\mscr{I})} x, u \big).
\end{equation}
In sparsely coupled systems, the time derivative of each state $[x]_i$ depends only on a small number of neighbors, meaning that if $\mscr{I}$ is small compared to the state dimension $n$, then $\mscr{N}(\mscr{I})$ is also small compared to $n$.

Suppose that the weight matrix $\Psi_L\in \R^{n\times n_{L-1}}$ describing the input layer of the encoder has nonzero entries only in the rows indexed by $\mscr{I}=\{i_1, \ldots, i_K\}$.
Assembling the sub-matrix $\tilde{\Psi}_L\in \R^{K\times n_{L-1}}$ from these nonzero rows, we have $\Psi_L^T = \tilde{\Psi}_L^T S_{\mscr{I}}$.
In the notation of Section~\ref{subsec:tensors}, this means that the reduced-order model (\ref{eqn:outer_layer_ROM}) can be written in terms of $f_{\mscr{I}}$ as
\begin{equation}
\begin{aligned}
    \tilde{f}(z,u) &= \D \tilde{\psi}_e\big( \tilde{\psi}_d(z) \big) \tilde{\Psi}_L^T f_{\mscr{I}}\big( S_{\mscr{N}(\mscr{I})} b_L + S_{\mscr{N}(\mscr{I})} \Phi_L \tilde{\psi}_d(z), \ u \big)  \\
    \tilde{g}(z) &= g\big(b_L + \Phi_L \tilde{\psi}_d(z)\big).
\end{aligned}
\label{eqn:sparsified_ROM}
\end{equation}
If the number of neighbors described by the set $\mscr{N}(\mscr{I})$ is small compared to the original state dimension $n$, then we can obtain computational speedups by evaluating $f_{\mscr{I}}$ instead of $f$.

Note that because the columns of $\Psi_L$ are linearly independent, we must have $K \geq n_{L-1}$, and in general there are at least $K$ elements in $\mscr{N}(\mscr{I})$.
This means that the dimension $n_{L-1}$ must be chosen to be much smaller than the state dimension, and so Remark~\ref{rem:outer_layer_width} applies.
However, the cost to evaluate the time derivative of the ROM in (\ref{eqn:sparsified_ROM}) does not grow as rapidly with $n_{L-1}$ as in the tensor-based method described in Section~\ref{subsec:tensors}.
In the case of PDEs discretized in space using finite-difference schemes with small stencils, the number of neighboring elements in $\mscr{N}(\mscr{I})$ and the cost to evaluate $f_{\mscr{I}}$ will grow linearly with the size of $\mscr{I}$.
Therefore, in the best-case scenario where the number of nonzero rows of $\Psi_L$ grows linearly with $n_{L-1}$, then the cost to evaluate the time derivative of the ROM will also scale linearly with $n_{L-1}$.

The simplest way obtain a sparse $\Psi_L$ is to constrain which rows can have nonzero entries prior to training.
The row indices can be chosen using methods such as random selection, coarsening a spatial grid, or QR-pivoting-based DEIM \cite{Drmac2016new}.
However, choosing the nonzero rows of $\Psi_L$ prior to optimization may prevent the encoder from learning a useful direction of projection.

Better performance can likely be achieved by learning a sparse $\Psi_L$ during the training process for the autoencoder.
One option is to add a sparsity-promoting penalty on $\Psi_L$ to the cost function used to train the autoencoder.
This penalty should not introduce additional biases including weight matrix shrinkage into the optimization problem since this can affect the learned manifold and projection fibers.
For example, an $\ell^1$ penalty (see \citet{Tibshirani1996regression}) with a large weight factor will shrink the encoder weight matrix $\Psi_L$ towards zero, while pushing the corresponding decoder weight matrix $\Phi_L$ towards infinity due to the biorthogonality constraint.
Other sparsity-promoting penalties such as those in \cite{Yuan2006model, Scardapane2017group, Koneru2019sparse, Wang2017novel} have this same issue in our setting.
Instead, for a matrix $\Psi \in \R^{n\times r}$ (dropping the subscript $L$) with linearly independent columns, we construct $U\in \R^{n\times r}$ having orthonormal columns spanning $\Range(\Psi)$, for example 
% via $U = \Psi (\Psi^T \Psi)^{-1/2}$ or 
via QR factorization $U = \qf(\Psi)$.
Our proposed penalty function is then defined by
\begin{empheq}[box=\widefbox]{multline}
     R_{1,2}\big(\Range(\Psi)\big)
    % = \big\Vert \Range(\Psi) \big\Vert_{1,2} 
    = \big\Vert U \big\Vert_{1,2} - r \\
    = \sum_{i=1}^{n} \big\Vert \row_i(U) \big\Vert_2 - r,
    \label{eqn:sparsity_penalty}
\end{empheq}
where $\Vert \cdot \Vert_{1,2}$ denotes the sum of Euclidean norms of the rows of a matrix.
This function does not depend on the choice of $U$ since $\Vert \cdot \Vert_{1,2}$ is invariant under multiplication on the right by $r\times r$ orthonormal matrices.
Most importantly, the penalty defined by (\ref{eqn:sparsity_penalty}) depends only on the range of $\Psi$ since it remains invariant under changes of basis, i.e., when $\Psi$ is replaced by $\Psi A$ for any invertible matrix $A\in\R^{r\times r}$.
Indeed, it defines a continuous function on the Grassmann manifold $\mcal{G}_{n,r}$ consisting of $r$-dimensional subspaces of $\R^n$ (see \citet{Bendokat2020grassmann, Absil2004riemannian, Wong1967differential}).
The following theorem shows that this penalty does in fact promote sparsity of the rows of $\Psi$.
\begin{theorem}
    The minimum value of the penalty function defined by (\ref{eqn:sparsity_penalty}) over the space $\R_*^{n\times r}$ of real $n\times r$ matrices with linearly independent columns is zero.
    This value is attained by $\Psi\in \R_*^{n\times r}$ if and only if $\Psi$ has precisely $r$ rows with nonzero entries.
    \label{thm:minimizers_of_sparse_penalty}
\end{theorem}
\begin{proof}
    We give the proof in Appendix~\ref{app:Grassmannian_sparsity}.
\end{proof}
Moreover, the penalty function increases sharply (in much the same way as $x \mapsto \vert x \vert$) in the neighborhood of its minimizers.
Specifically, we have Corollary~\ref{cor:local_behavior_of_12fun} in Appendix~\ref{app:Grassmannian_sparsity}, which we have not stated here as it requires machinery for the Grassmann manifold that is beyond the scope of this paper.
This result implies that the penalty produces sparse minimizers when it is added with a sufficiently large, but finite factor to smooth optimization objectives.
More precisely, we have the following theorem.
\begin{theorem}
    \label{thm:solutions_of_12_regularized_minimization}
    Let $\mcal{M}$ be a smooth manifold and let $D(J_0)$ be an open subset of $\mcal{M}\times \mcal{G}_{n,r}$ on which a real non-negative-valued function $J_0$ is defined and continuously differentiable.
    Suppose that there is a finite constant $M$ so that the preimage set $S_M = \{ (x,\mcal{V}) \in D(J_0) \ : \ J_0(x,\mcal{V}) \leq M \}$ is compact and contains a point $(x_0, \mcal{V}_0)$ so that $ R_{1,2}(\mcal{V}_0) = 0$.
    Then for any $\gamma \geq 0$, the function on $D(J_0)$ defined by
    \begin{equation}
        J_{\gamma}(x,\mcal{V}) = J_0(x,\mcal{V}) + \gamma  R_{1,2}(\mcal{V})
    \end{equation}
    attains its minimum and all such minimizers lie in $S_M$.
    Furthermore, there is a constant $\Gamma \geq 0$ so that when $\gamma > \Gamma$, every minimizer $(x_*, \mcal{V}_*)$ of $J_{\gamma}$ satisfies $ R_{1,2} ( \mcal{V}_* ) = 0$.
\end{theorem}
\begin{proof}
    The proof of this result uses tools from Grassmannian geometry that are beyond the scope of this paper. 
    We provide the necessary background, lemmata, and proof in Appendix~\ref{app:Grassmannian_sparsity}.
\end{proof}
As a consequence of this theorem, a sparse matrix $\Psi_L$ in the encoder with precisely $n_{L-1}$ nonzero rows can be obtained by minimizing a cost function to which (\ref{eqn:sparsity_penalty}) has been added with a sufficiently large factor.
Increasing the factor beyond this point has no further affect on the minimizers; specifically, there is no additional shrinkage of the weight matrix $\Psi_L$.
In practice, we suggest first optimizing the network without the sparsity-promoting penalty, then activating the penalty during a subsequent optimization stage to sparsify $\Psi_L$.

\section{Conclusion}
\label{sec:conclusion}

In this paper we develop a nonlinear projection-based model reduction framework in which it is possible to learn both a low-dimensional manifold and appropriate projection fibers for capturing transient dynamics away from the manifold.
% As we demonstrate, utilizing appropriate projection fibers is key to capturing transient dynamics in the neighborhood of a learned manifold.
To do this, we introduce a new autoencoder neural network architecture defining a parametric class of nonlinear projections along with new dynamics-aware cost functions for training.
In order to define a nonlinear projection, we ensure that the encoder is a left inverse of the decoder by utilizing a new pair of invertible activation functions and enforcing a biorthogonality constraint on the weight matrices.
The biorthogonality constraint defines a smooth matrix manifold on which the optimization during training takes place.

As we demonstrate, optimizing the autoencoder on standard reconstruction-based loss does not generally yield appropriate projection fibers for capturing transient dynamics.
To address this problem, we introduce two new cost functions based on additional information from the full-order model.
The first cost function, which we call Reconstruction and Velocity Projection (RVP) loss, is based on a Gr\"{o}nwall-Bellman-type error analysis of the reduced-order model.
It entails adding a time-derivative (``velocity'') projection loss to the usual reconstruction-based loss.
The second cost function, which we call Gradient-Aligned Projection (GAP) loss, is based on a first-order Taylor expansion of projection-based forecasting error.
This analysis yields a cost function measuring the alignment of state projection errors with randomized gradient samples along trajectories.
Both of these new loss function require us to be able to query the adjoint of the full-order model, acting on vector, and thus are not suitable if only experimental data is available.

We present a detailed study comparing our framework to state-of-the-art methods on a simple three-state model, introduced by \citet{Noack2003hierarchy}, of vortex shedding in the wake of a bluff body.
Regardless of the cost function and neural network architecture, the autoencoders we trained were able to accurately locate the two-dimensional slow manifold in this problem.
Nonetheless, the cost function used to train the networks had a large effect on the resulting model's ability to forecast trajectories with initial conditions lying away from the slow manifold.
Training on reconstruction loss consistently produced inaccurate models with projection fibers failing to cancel the fast coordinate.
Both of our new cost functions were able to remedy this issue, with RVP loss yielding slightly better performance than GAP loss and suffering from less deterioration in performance on a smaller training data set.
For the forecasting task, our proposed architecture trained using either GAP or RVP loss significantly outperformed standard architectures and loss functions.

% \seoremark{The following paragraph aims to address comments by both reviewers about prospects for scaling to high-dimensional systems.}
While we have discussed several methods for constructing computationally efficient reduced-order models, we have not yet applied our method to high-dimensional systems. 
This will be an important direction for future work.
In particular, we will be interested in studying whether, or to what extent the amount of training data required to obtain an accurate ROM scales with the state dimension of the FOM.
We have reason to expect favorable scaling behavior because the data requirements for computing CoBRAS\cite{Otto2022model} projections, which minimize a loss similar to GAP, do not scale with the dimension of the FOM, but rather with the effective ranks of covariance matrices for states and gradient data.
This suggests that using a loss function like GAP, or the gradient-weighted CoBRAS loss, might allow for dimension-independent scaling of the training data set for certain systems with low-dimensional underlying manifolds and few directions of high sensitivity in the state space.
We will also be interested in the performance of our proposed encoder sparsification technique, which may also reduce data requirements when the added bias towards sparsity is appropriate.
Finally, in follow-up work we aim to provide some practical guidelines for choosing the number of layers and their widths in applications to high-dimensional systems.

Other directions for future work include developing convolutional autoencoders with similar constraints, as well as applying our autoencoder architecture for other tasks such as preprocessing data from dynamical systems, or as part of a SINDy-autoencoder\cite{Champion2019data}.
% Other constraints such as boundary conditions and incompressibility in fluid flows could be enforced by further constraining the weight matrix $\Phi_L$ and bias vector $b_L$ defining the final layer of the decoder.
% In principle any set of linear constraints on state vectors can be enforced in this manner.
Further investigation into data sampling strategies, especially in the presence of unstable structures in state space may also lead to practical guidelines for reduced-order modeling using our framework.
Another exciting direction for future work will be to use our autoencoder to approximate solutions of the equations derived by \citet{Roberts1989appropriate, Roberts2000computer} for the correct spatially-varying affine projections.
For this, a method analogous to physics-informed neural networks (PINNs)\cite{Raissi2019physics} could be employed.

\begin{acknowledgments}
This work was supported by the Air Force Office of Scientific Research, award FA9550-19-1-0005.
\end{acknowledgments}

\section*{Author declarations}

\subsection*{Conflict of interest}
% required, see AIP instructions https://publishing.aip.org/resources/researchers/author-instructions/#prep
The authors have no conflicts to disclose.

\subsection*{Author contributions}
% required, see AIP instructions https://publishing.aip.org/resources/researchers/author-instructions/#prep

% Example (see https://publishing.aip.org/resources/researchers/policies-and-ethics/authors/)
% Amanda Green: review and editing (equal). Kerry Jones: Conceptualization (lead); writing – original draft (lead); formal analysis (lead); writing – review and editing (equal). Roberto Solis: Software (lead); writing – review and editing (equal). Hao Wang: Methodology (lead); writing – review and editing (equal). Jenny Wu: Conceptualization (supporting); Writing – original draft (supporting); Writing – review and editing (equal).

\textbf{Samuel E. Otto:} conceptualization (lead); formal analysis (lead); methodology (lead); writing -- original draft (lead); supervision (supporting).
\textbf{Gregory R. Macchio:} software (lead); visualization (lead); writing -- original draft (supporting).
\textbf{Clarence W. Rowley:} funding acquisition (lead); supervision (lead); resources (lead); writing -- review \& editing (lead); conceptualization (supporting).

% S.E.O developed the autoencoder architecture and proved all stated theoretical results.
% He also wrote the majority of the paper with the exception being Section~\ref{sec:noack_case_study}, which was primarily written by G.R.M.
% All of the code and numerical results were written and generated by G.R.M.
% He also proposed the illustrative model in Example~\ref{ex:toy_slow_fast}.
% C.W.R Acquired funding, oversaw the research, and first proposed the idea of an autoencoder in which the encoder is a left inverse of the decoder in a discussion with S.E.O regarding the work in \citet{Otto2019linearly}.

\section*{Data Availability Statement}
% required, see AIP instructions https://publishing.aip.org/resources/researchers/author-instructions/#prep

Data sharing is not applicable to this article as no new data were created or analyzed in this study.
Our code was written in Python and is available at  \url{https://github.com/grmacchio/romnet_chaos2023} (Gregory R. Macchio's GitHub).

\bibliography{jfull,references}% Produces the bibliography via BibTeX.

\appendix

\section{The biorthogonal manifold}
\label{app:biorthogonal_mfd}

In this appendix we provide supplementary information about the biorthogonal manifold and our over-parametrization which are relevant for optimization.
As we mentioned in Remark~\ref{rem:optimizing_on_Bnr_directly}, these results can be used to implement alternative Riemannian optimization algorithms (see \citet{Absil2009optimization}) on the biorthogonal manifold.
These algorithms require two ingredients called retraction and vector transport, which we provide below.
We also discuss topological properties of the biorthogonal manifold, the over-parametrization domain, and the ramifications of these properties for optimization.

% In particular the over-parametrization function yields a retraction on the birothogonal manifold, and (\ref{eqn:orthogonal_projection_onto_TBnr}), below, can be used for vector transport \cite{Otto2022advances, Absil2009optimization}.

The biorthogonal manifold along with its tangent and normal spaces are characterized by the following theorem.
\begin{theorem}[The biorthogonal manifold]
\label{thm:biorthogonal_manifold}
Let $n \geq r \geq 1$ be integers.
The set of biorthogonal matrices 
\begin{equation}
    \mcal{B}_{n,r} = \left\lbrace (\Phi, \Psi)\in \R^{n\times r}\times \R^{n\times r} \ : \ \Psi^T\Phi = I_r  \right\rbrace,
    \label{eqn:biorthogonal_manifold_defn}
\end{equation}
is a smooth, closed $2nr-r^2$ dimensional submanifold of $\R^{n\times r}\times \R^{n\times r}$, with tangent space at a point $(\Phi, \Psi) \in \mcal{B}_{n,r}$ given by
\begin{equation}
    T_{(\Phi,\Psi)}\mcal{B}_{n,r} = \left\{ (X,Y)\in\R^{n\times r}\times \R^{n\times r} \ : \ Y^T \Phi + \Psi^T X = 0 \right\}.
\end{equation}
When $\R^{n\times r}\times \R^{n\times r}$ is endowed with the Euclidean inner product
\begin{equation}
    \big\langle (X_1,Y_1),\ (X_2, Y_2) \big\rangle = \Tr(X_1^T X_2) + \Tr(Y_1^T Y_2),
\end{equation}
then the normal space of the biorthogonal manifold is given by
\begin{equation}
    \left(T_{(\Phi,\Psi)}\mcal{B}_{n,r}\right)^{\perp} = \left\{ (\Psi A, \ \Phi A^T)\in\R^{n\times r}\times \R^{n\times r} \ : \ A \in \R^{r\times r} \right\}.
\end{equation}
The orthogonal projection of any $(X,Y)\in\R^{n\times r}\times \R^{n\times r}$ onto $T_{(\Phi,\Psi)}\mcal{B}_{n,r}$ is given by
\begin{equation}
    P_{(\Phi,\Psi)}(X,Y) = \left( X - \Psi A, \ Y - \Phi A^T \right),
    \label{eqn:orthogonal_projection_onto_TBnr}
\end{equation}
where $A\in\R^{r\times r}$ is the unique solution of the Sylvester equation
\begin{equation}
    A (\Phi^T \Phi) + (\Psi^T \Psi) A = Y^T \Phi + \Psi^T X.
\end{equation}
% The Sylvester equation is equivalent to the symmetric, positive-definite linear system
% \begin{equation}
%     \left[ \left(\Phi^T \Phi\right)\otimes I_r +  I_r\otimes\left(\Psi^T \Psi\right) \right]\vct(A)
%     = \vct\left( Y^T \Phi + \Psi^T X \right), \quad
%     \vct(A) = \begin{bmatrix} \col_1(A) \\ \vdots \\ \col_r(A) \end{bmatrix}.
% \end{equation}
\end{theorem}
\begin{proof}
We construct $\mcal{B}_{n,r}$ as the preimage of the regular value $0$ under the map $F:(X,Y) \mapsto Y^T X - I_r$ using the preimage theorem \cite{Guillemin1974differential}.
The tangent space to $\mcal{B}_{n,r}$ at $(\Phi,\Psi)$ is given by the null-space of $\D F(\Phi,\Psi)$ according to the local submersion theorem \cite{Guillemin1974differential}.
We give the details of the proof in Appendix~\ref{app:proofs}.
\end{proof}
We observe that the orthogonal projection onto the tangent space given by (\ref{eqn:orthogonal_projection_onto_TBnr}) can be used to provide a vector transport on $\mcal{B}_{n,r}$.
Specifically, if $R_{(\Phi,\Psi)}: T_{(\Phi,\Psi)} \mcal{B}_{n,r} \to \mcal{B}_{n,r}$ is a retraction (see \citet[Definition~4.1.1]{Absil2009optimization}) then
\begin{equation}
    \mcal{T}_{(\Phi,\Psi), (X_1,Y_1)}: (X_2, Y_2) \mapsto P_{R_{(\Phi,\Psi)}(X_1,Y_1)} (X_2, Y_2)
\end{equation}
is easily seen to satisfy the required conditions (\citet[Definition~8.1.1]{Absil2009optimization}) to be a vector transport.

The following theorem characterizes the local structure of the over-parametrization function in its domain.
It says that smooth local coordinates can be chosen about any point in the domain $D_+(\Pi_{n,r})$ so that the first $2 n r - r^2$ are coordinates for a patch of the biorthogonal manifold.
The over-parametrization does not depend on the remaining $r^2$ coordinates.
\begin{theorem}[Over-parametrization]
\label{thm:submersion_onto_Bnr}
    The map $\Pi_{n,r}:D_+(\Pi_{n,r}) \to \mcal{B}_{n,r}$ defined by (\ref{eqn:projection_onto_Bnr}) is a surjective submersion, i.e., its tangent map $\D\Pi_{n,r}:T D_+(\Pi_{n,r}) \to T\mcal{B}_{n,r}$ is surjective.
    There are several consequences:
    \begin{enumerate}
        \item The preimage set $\Pi_{n,r}^{-1}(\Phi, \Psi)$ of each $(\Phi, \Psi)\in \mcal{B}_{n,r}$ is a smooth, closed $r^2$-dimensional submanifold of $D_{+}(\Pi_{n,r})$ intersecting $\mcal{B}_{n,r}$ transversally.
        \item For any $(\tilde{\Phi},\tilde{\Psi}) \in D_+(\Pi_{n,r})$ there is an open neighborhood $\mcal{U}$ of $(\tilde{\Phi},\tilde{\Psi})$ in $D_+(\Pi_{n,r})$ so that $\mcal{V} = \Pi_{n,r}(\mcal{U})$ is an open neighborhood of $(\Phi, \Psi) = \Pi_{n,r}(\tilde{\Phi},\tilde{\Psi})$ in $\mcal{B}_{n,r}$ and on these neighborhoods there are smooth coordinates $\phi:\mcal{U}\to\R^{2nr}$ and $\psi:\mcal{V}\to\R^{2nr-r^2}$ such that
        \begin{equation}
            (\psi \circ \Pi_{n,r} \circ \phi^{-1})(x_1, \ldots, x_{2nr}) = (x_1, \ldots, x_{2nr-r^2})
        \end{equation}
        on the open subset $\phi(\mcal{U}) \subset \R^{2nr}$.
        \item In these coordinate neighborhoods, a function $J:\mcal{B}_{n,r}\to\R$ and its composition with $\Pi_{n,r}$ are related by
        \begin{equation}
            J\circ\Pi_{n,r}\circ\phi^{-1}(x_1, \ldots, x_{2nr}) = J \circ \psi^{-1}(x_1, \ldots, x_{2nr-r^2}).
        \end{equation}
        \item The function $J$ is smooth if and only if $J\circ\Pi_{n,r}$ is smooth, and if so, their gradients (with Riemannian metrics inherited from the ambient Euclidean spaces) are related by
        \begin{equation}
            \grad J( \Phi, \Psi ) = 
            G(\tilde{\Phi}, \tilde{\Psi})^{-1}
            \D \Pi_{n,r}(\tilde{\Phi}, \tilde{\Psi}) \grad(J \circ \Pi_{n,r})(\tilde{\Phi}, \tilde{\Psi}),
            \label{eqn:overparametrized_gradient_relation}
        \end{equation}
        where $G(\tilde{\Phi}, \tilde{\Psi}) = \D \Pi_{n,r}(\tilde{\Phi}, \tilde{\Psi}) \D \Pi_{n,r}(\tilde{\Phi}, \tilde{\Psi})^*$ is invertible on $T_{(\Phi, \Psi)}\mcal{B}_{n,r}$ and $\grad(J \circ \Pi_{n,r})(\tilde{\Phi}, \tilde{\Psi}) \in \big( T_{(\tilde{\Phi},\tilde{\Psi})}\Pi_{n,r}^{-1}(\Phi, \Psi) \big)^{\perp}$ is orthogonal to the fiber.
    \end{enumerate}
\end{theorem}
\begin{proof}
    Direct computation using the formula for differentiating the matrix inverse shows that the tangent map is surjective.
    Transversality of the fiber and the biorthogonal manifold follows from the fact that $\Pi_{n,r}$ is idempotent and an argument resembling Theorem~1.15 in \citet{Michor2008topics}.
    The remaining properties follow from standard results characterizing smooth surjective submersions that can be found in \citet{Lee2013introduction} or in \citet{Guillemin1974differential}.
    We give the details in Appendix~\ref{app:proofs}.
\end{proof}
We observe that the over-parametrization function can be used to provide a retraction on the biorthogonal manifold.
% In fact, it is a ``projection-like retraction'' \cite{Absil2012projection}.
Specifically, it is easily verified that
\begin{equation}
    R_{(\Phi, \Psi)}: (X,Y) \mapsto \Pi_{n,r}(\Phi + X, \Psi + Y)
\end{equation}
satisfies the required conditions (\citet[Definition~4.1.1]{Absil2009optimization}) to be a retraction.
In fact, this is a projection-like retraction \cite{Absil2012projection}.

The following result shows that the domain of the over-parametrization $D_+(\Pi_{n_{l}, n_{l-1}})$ for the $l$th layer of our autoencoder is connected when the dimensions of the layer decrease, i.e., when $n_{l} > n_{l-1}$.
This means that restricting the optimizer to this domain does not cut off access to any part of the biorthogonal manifold by an optimization algorithm that follows a continuous path or proceeds in small steps.
\begin{proposition}
    \label{prop:overparametrization_domain}
    If $n > r\geq 1$ then $D_+(\Pi_{n,r})$ is connected.
    If $n = r\geq 1$ then $D_+(\Pi_{n,n})$ is a union of two disjoint connected components
    \begin{subequations}
    \begin{align}
        D_{+,+}(\Pi_{n,n}) &= \left\{ (\tilde{\Phi}, \tilde{\Psi}) \in \R^{n,n} \times \R^{n,n}  \ : \ \det(\tilde{\Phi}), \det(\tilde{\Psi}) > 0 \right\}, \\
        D_{-,-}(\Pi_{n,n}) &= \left\{ (\tilde{\Phi}, \tilde{\Psi}) \in \R^{n,n} \times \R^{n,n}  \ : \ \det(\tilde{\Phi}), \ \det(\tilde{\Psi}) < 0 \right\}.
    \end{align}
    \label{eqn:components_of_domain}
    \end{subequations}
\end{proposition}
\begin{proof}
    The $n=r$ case follows from the fact that the general linear group has two disjoint connected components corresponding to matrices with positive and negative determinants.
    To prove connectedness when $n > r$, we rely on the path constructed in the proof of Theorem.~3.1 in \citet{Otto2022optimizing} between biorthogonal matrix pairs.
    The details are provided in Appendix~\ref{app:proofs}.
\end{proof}

When the $l$th layer is square, i.e., when  $n_{l} = n_{l-1}$, then the domain $D(\Pi_{n_{l}, n_{l-1}})$, and hence the biorthogonal manifold $\mcal{B}_{n_{l}, n_{l-1}}$, has two disjoint components.
% In this case, the sequence of iterates produced by an optimizer following a continuous path will all lie within one of these components and not the other.
However, this is of little consequence for the network because
any choice of signs for the determinants of the square weight matrices in the network can be achieved without changing the overall projection $P = \psi_d \circ \psi_e$ as we now explain.
The key observation is that swapping a pair of nodes in layer $l-1$ will flip the sign of both determinants without changing $P$.
In particular, suppose $S$ is an $n_{l-1}\times n_{l-1}$ permutation matrix.
We rearrange the nodes in layer $l-1$ by replacing the weight matrices in layer $l$ with
\begin{equation}
    (\Phi_l,\ \Psi_l) \mapsto (\Phi_l S,\ \Psi_l S),
\end{equation}
and if $l > 1$, we also replace the weight matrices and bias vector in layer $l-1$ with
\begin{equation}
    (\Phi_{l-1},\ \Psi_{l-1}) \mapsto (S^T \Phi_{l-1},\ S^T \Psi_{l-1}), \quad b_{l-1} \mapsto S^T b_{l-1}.
\end{equation}
Since the activation functions act element-wise, they commute with permutation, i.e., $S\circ \sigma_{\pm} = \sigma_{\pm}\circ S$.
Recalling that the layers of the autoencoder are defined by (\ref{eqn:layer_definition}), we observe that performing this operation on any layer $l > 1$ leaves $\psi_d$ and $\psi_e$ unchanged.
If the permutation is performed on the first layer $l=1$, then the operation yields
\begin{equation}
    \psi_d \mapsto \psi_d\circ S, \qquad \psi_e \mapsto S^T \circ \psi_e,
\end{equation}
which leaves $\psi_d \circ \psi_e = \psi_d \circ S^T S \circ \psi_e$ unchanged.
We could swap nodes using this procedure beginning with layer $L$ and proceed in reverse order down to layer $1$ to achieve any desired sequence of signs of determinants in the square layers of the network while leaving $P=\psi_d \circ \psi_e$ unchanged.

\section{Proofs}
\label{app:proofs}

We use the following handy modification of the Gr\"{o}nwall-Bellman inequality in the proof of Proposition~\ref{prop:weighted_VPE}.
While its proof employs an argument similar to the standard Gr\"{o}nwall-Bellman inequality (see any book on nonlinear dynamical systems such as Khalil \cite{Khalil2002nonlinear} or Kelly and Peterson \cite{Kelly2004theory}), it is apparently absent from the standard literature.
Though the inequality was given in Otto's thesis \cite{Otto2022advances}, we reproduce it here for completeness.
\begin{lemma}[Inhomogeneous Gr\"{o}nwall-Bellman inequality]
\label{lem:inhomogeneous_Gronwall}
Suppose that $w:[0,T] \to \R$ and $b:[0,T] \to \R$ are integrable functions satisfying
\begin{equation}
    w(t) \leq a + \int_0^t\left[ L w(\tau) + b(\tau) \right]  \td \tau \qquad \forall t\in [0,T]
\end{equation}
for some constants $a, L \in \R$.
Then, $w$ is bounded according to
\begin{equation}
    w(t) \leq a e^{L t} + \int_0^t e^{L(t-\tau)} b(\tau)  \td \tau \qquad \forall t\in [0,T].
\end{equation}
\end{lemma}
\begin{proof}[Proof of Lemma~\ref{lem:inhomogeneous_Gronwall}]
We modify the proof of the Gr\"{o}nwall-Bellman inequality given in Khalil \cite{Khalil2002nonlinear}.
Let us define the function
\begin{equation}
    v(t) = e^{-Lt}\int_0^t \left[ L w(\tau) + b(\tau) \right]  \td \tau
\end{equation}
and observe that
\begin{multline}
    v'(t) = e^{-L t}\left\{ -L \int_0^t\left[ L w(\tau) + b(\tau) \right]  \td \tau + L w(t) + b(t) \right\} \\
    \leq e^{-L t}\left\{ L a + b(t) \right\}.
\end{multline}
Integrating, and noting that $v(0) = 0$ we find
\begin{equation}
    v(t) \leq a - a e^{-L t} + \int_0^t e^{-L \tau} b(\tau)  \td  \tau
\end{equation}
and so we obtain
\begin{equation}
    w(t) \leq a + e^{L t} v(t) \leq a e^{L t} + \int_0^t e^{L(t-\tau)} b(\tau)  \td \tau.
\end{equation}
\end{proof}

\begin{proof}[Proof of Proposition~\ref{prop:weighted_VPE}: Weighted velocity projection error]
    To simplify the notation, we denote
    \begin{equation}
        h(t) := \Big\Vert \ddt x_P(t) - \hat{f}_P(x_P(t),u(t))\Big\Vert
    \end{equation}
    and we let $e(t) := x_P(t) - \hat{x}(t)$.
    We have $e(0) = 0$ by definition of the initial condition in (\ref{eqn:P_ROM}) and
    \begin{equation}
        e(t) = \int_{0}^{t} \left[ \ddt x_P(\tau) - \hat{f}_P(\hat{x}(\tau),u(\tau)) \right] \td \tau.
    \end{equation}
    Since $x \mapsto \hat{f}_P(x,u(t))$ is Lipschitz, we obtain
    \begin{equation}
        \Vert e(t) \Vert \leq \int_{0}^{t} \big[ h(\tau) 
        + L \Vert e(\tau) \Vert \big] \td \tau.
    \end{equation}
    By a simple modification of the Gr\"{o}nwall-Bellman inequality stated in Lemma~\ref{lem:inhomogeneous_Gronwall}, it follows that
    \begin{equation}
        \Vert e(t) \Vert \leq \int_{0}^{t} e^{L(t - \tau)} h(\tau) \td \tau.
    \end{equation}
    Applying the Cauchy-Schwarz inequality to the above gives
    \begin{equation}
        \Vert e(t) \Vert^2 \leq 
        % \int_{0}^t e^{2 L (t-\tau)} \td \tau \int_{0}^t h(\tau)^2 \td \tau = 
        \frac{1}{2L}\left( e^{2 L t} - 1 \right) \int_{0}^t h(\tau)^2 \td \tau.
    \end{equation}
    We now integrate over the time interval $[0,t_f]$ and exchange the order of integration to obtain
    \begin{equation}
    \begin{aligned}
        \int_{0}^{t_f} \Vert e(t) \Vert^2 \td t \leq & \int_{0}^{t_f} \int_{0}^t \frac{1}{2L}\left( e^{2 L t} - 1 \right) h(\tau)^2 \td \tau \td t \\
        &= \int_{0}^{t_f} \int_{\tau}^{t_f} \frac{1}{2L}\left( e^{2 L t} - 1 \right) h(\tau)^2 \td t \td \tau \\
        &= \int_{0}^{t_f} w_{L,t_f}(\tau) h(\tau)^2 \td \tau.
    \end{aligned}
    \end{equation}
    Since this is (\ref{eqn:velocity_error_weighting}), the proof is complete.
\end{proof}

\begin{proof}[Proof of Lemma~\ref{lem:GAP_Taylor_estimate}]
    By compactness of $\mcal{X}\times\mcal{U}$ and continuity of the derivatives of $(x, u) \mapsto F_u(x)$ with respect to $x$ up to second order, the maximum absolute value of these derivatives is bounded.
    By Taylor's theorem (in particular, Corollary~C.16 in \citet{Lee2013introduction}) there is a constant $\tilde{C} \geq 0$ so that
    \begin{equation}
        \big\Vert F_u(x) - F_u(P(x)) - \D F_u(x) (x - P(x)) \big\Vert \leq \tilde{C} \Vert x - P(x) \Vert^2
    \end{equation}
    for every $(x, u) \in \mcal{X}\times\mcal{U}$.
    By the triangle inequality,
    \begin{multline}
        \big\Vert F_u(x) - F_u(P(x)) \big\Vert^2 \\
        \leq \left( \big\Vert \D F_u(x) (x - P(x)) \big\Vert + \tilde{C} \Vert x - P(x) \Vert^2 \right)^2.
    \end{multline}
    Since $\Vert x - P(x) \Vert$ is bounded and $\D F_u(x)$ is bounded in operator norm, there is a constant $C \geq 0$ so that (\ref{eqn:first_order_difference_estimate}) holds.
\end{proof}

\begin{proof}[Proof of Theorem~\ref{thm:orthogonality_promoting_regularization}]
    We begin by proving the second claim, that $\Vert \Phi \Vert_F^2 + \Vert \Psi \Vert_F^2 \to \infty$ if any of the principal angles between $\Range(\Phi)$ and $\Range(\Psi)$ approach $\pi/2$.
    The columns of $U = \Phi G_{\Phi}^{-1/2}$ and $V = \Psi G_{\Psi}^{-1/2}$ are orthonormal bases for $\Range(\Phi)$ and $\Range(\Psi)$. 
    As shown by \citet{Bjorck1973numerical}, the cosines of the principal angles between these subspaces are equal to the singular values of
    \begin{equation}
        A := U^T V = G_{\Phi}^{-1/2} G_{\Psi}^{-1/2}.
    \end{equation}
    If any of these principal angles approaches $\pi/2$, then the corresponding singular value $\sigma_i$ of $U^T V$ approaches zero.
    Therefore,
    \begin{equation}
    \begin{aligned}
        \Vert \Phi \Vert_F^2 + \Vert \Psi \Vert_F^2
        &= \big\Vert G_{\Phi}^{1/2} \big\Vert_F^2 + \big\Vert G_{\Psi}^{1/2} \big\Vert_F^2 \\
        &\geq 2 \big\Vert G_{\Phi}^{1/2} \big\Vert_F \big\Vert G_{\Psi}^{1/2} \big\Vert_F \\
        &\geq 2 \big\Vert G_{\Psi}^{1/2} G_{\Phi}^{1/2} \big\Vert_F 
        = 2 \big\Vert (U^T V)^{-1} \big\Vert \to \infty. 
    \end{aligned}
    \end{equation}
    This proves the second claim.

    It is evident by substitution that $\Vert \Phi \Vert_F^2 + \Vert \Psi \Vert_F^2 = 2r$ when $\Phi = \Psi$ have orthonormal columns.
    We now show that this is the minimum possible value of $\Vert \Phi \Vert_F^2 + \Vert \Psi \Vert_F^2$ on the biorthogonal manifold.
    We observe that $(A A^T)^{-1} \succeq I$ with respect to the positive-semidefinite Loewner ordering because all of the singular values of $A$ are in the interval $(0,1]$ by Cauchy-Schwarz.
    Hence, we have
    \begin{equation}
        % G_{\Psi}^{-1/2} = G_{\Phi}^{1/2} A
        % G_{\Psi}^{-1} = G_{\Phi}^{1/2} A A^T G_{\Phi}^{1/2}
        G_{\Psi} 
        = G_{\Phi}^{-1/2} (A A^T)^{-1} G_{\Phi}^{-1/2} 
        \succeq G_{\Phi}^{-1}.
        \label{eqn:GPsi_in_terms_of_GPhi_and_A}
    \end{equation}
    If $\lambda_1, \ldots, \lambda_r$ denote the eigenvalues of $G_{\Phi}$, this means that
    \begin{equation}
        \Vert \Phi \Vert_F^2 + \Vert \Psi \Vert_F^2
        \geq \Tr(G_{\Phi}) + \Tr(G_{\Phi}^{-1})
        = \sum_{i=1}^r\left( \lambda_i + \frac{1}{\lambda_i} \right) \geq 2 r.
        \label{eqn:Frobenius_regularization_lower_bound}
    \end{equation}
    Now we show that equality in the above equation implies that $\Phi = \Psi$ has orthonormal columns.
    Equality in the above implies that $\lambda_1 = \cdots = \lambda_r = 1$, meaning that $G_{\Phi} = I$.
    Therefore, $\Phi$ has orthonormal columns.
    An identical argument swapping the roles of $\Phi$ and $\Psi$ shows that $\Psi$ also has orthonormal columns.
    Since $\Psi^T \Phi = I$, we have
    \begin{equation}
        (\Psi v)^T (\Phi v) = v^T v = \Vert \Psi v \Vert \Vert \Phi v \Vert
    \end{equation}
    for any $v\in\R^r$.
    By Cauchy-Schwarz, $\Phi v$ and $\Psi v$ are linearly dependent.
    Since they have the same magnitude and $(\Psi v)^T (\Phi v) \geq 0$ we must have $\Phi v = \Psi v$.
    As this holds for all $v\in\R^r$, it follows that $\Phi = \Psi$, which completes the proof.
\end{proof}

\begin{proof}[Proof of Theorem~\ref{thm:biorthogonal_manifold}: The biorthogonal manifold]
Consider the smooth map $F:\R^{n\times r}\times \R^{n\times r}\to\R^{r\times r}$ defined by
\begin{equation}
    F(X, Y) = Y^T X - I_r.
\end{equation}
We observe that $\mcal{B}_{n,r} = F^{-1}(0)$ is the preimage of the zero matrix under the map $F$.
If $(\Phi,\Psi)\in \mcal{B}_{n,r}$, then the derivative of $F$ at this point is given by the map
\begin{equation}
    \D F (\Phi,\Psi):(X,Y) \mapsto Y^T \Phi + \Psi^T X, \qquad (X,Y)\in\R^{n\times r}\times \R^{n\times r}.
\end{equation}
It is easy to see that the derivative is surjective at every $(\Phi,\Psi)\in \mcal{B}_{n,r}$ because 
\begin{equation}
    \D F(\Phi,\Psi)(\Phi A, 0) = A, \qquad \forall A\in\R^{r\times r}.
\end{equation}
Hence, the identity matrix $I_r$ is a regular value for the smooth map $F$.
By the regular level set theorem (Corollary~5.14 in \citet{Lee2013introduction}), it follows that $\mcal{B}_{n,r}$ is a smooth, properly embedded (i.e., closed by Proposition~5.5 in \citet{Lee2013introduction}) sub-manifold of codimension $r^2$ in $\R^{n\times r}\times \R^{n\times r}$.
Since $\R^{n\times r}\times \R^{n\times r}$ is $2 n r$-dimensional, the dimension of $\mcal{B}_{n,r}$ is $2 n r - r^2$.

Since $F$ is a global defining function for $\mcal{B}_{n,r}$, Proposition~5.38 in \citet{Lee2013introduction} shows that the tangent space at $(\Phi,\Psi)\in \mcal{B}_{n,r}$ is characterized by the nullspace of the derivative, that is, 
\begin{multline}
    T_{(\Phi, \Psi)}\mcal{B}_{n,r} 
    = \Null \big(\D F(\Phi,\Psi)\big) \\ 
    = \left\lbrace (X,Y)\in\R^{n\times r}\times \R^{n\times r}\ :\ Y^T \Phi + \Psi^T X = 0 \right\rbrace.
\end{multline}

Since $\D F(\Phi,\Psi)$ is a finite-dimensional linear map between the Euclidean spaces $\R^{n\times r}\times \R^{n\times r}$  and $\R^{r\times r}$, we have
\begin{equation}
    \big(T_{(\Phi, \Psi)}\mcal{B}_{n,r}\big)^{\perp} = \Null \big( \D F(\Phi,\Psi)\big)^{\perp} = \Range \big(\D F(\Phi,\Psi)^*\big),
\end{equation}
where $\D F(\Phi,\Psi)^*:\R^{r\times r}\to\R^{n\times r}\times \R^{n\times r}$ is the adjoint of $\D F(\Phi,\Psi)$.
We claim that the adjoint operator is given by
\begin{equation}
    \D F(\Phi,\Psi)^*:A \mapsto ( \Psi A, \ \Phi A^T ).
\end{equation}
To verify this, we choose any $A\in\R^{r\times r}$ and $(X, Y)\in\R^{n\times r}\times \R^{n\times r}$ and observe that
\begin{align}
    \big\langle A,\ \D F(\Phi,\Psi)(X, Y) \big\rangle &= \Tr\big[A^T\left( Y^T \Phi + \Psi^T X \right) \big] \\
    % &= \Tr\big(A^T\Psi^T X \big) + \Tr\big(\Phi^T Y A\big) \\
    &= \Tr\big(A^T\Psi^T X \big) + \Tr\big(A\Phi^T Y\big) \\
    &= \big\langle (\Psi A, \Phi A^T),\ (X, Y) \big\rangle,
\end{align}
thanks to invariance of the trace under transposition and cyclic permutation.
Therefore, we can conclude that
\begin{equation}
    \left( T_{(\Phi, \Psi)}\mcal{B}_{n,r} \right)^{\perp} 
    % = \Range \left(\D F(\Phi,\Psi)^*\right)
    = \left\lbrace \big( \Psi A, \ \Phi A^T \big)\in\R^{n\times r}\times \R^{n\times r} \ :\ A\in\R^{r\times r} \right\rbrace.
\end{equation}

The orthogonal projection $(\hat{X}, \hat{Y}) = P_{(\Phi,\Psi)}(X, Y)$ is the unique element in $T_{(\Phi,\Psi)}\mcal{B}_{n,r}$ such that $(X - \hat{X},\ Y - \hat{Y})\in\left( T_{(\Phi,\Psi)}\mcal{B}_{n,r}\right)^{\perp}$.
By the characterization of the orthogonal complement of the tangent space shown above, we know that
\begin{equation}
    (X - \hat{X},\ Y - \hat{Y}) = (\Psi A,\  \Phi A^T)
\end{equation}
for some matrix $A\in\R^{r\times r}$.
Using the characterization of the tangent space shown above, the condition that $(\hat{X},\ \hat{Y}) = (X - \Psi A,\ Y - \Phi A^T ) \in T_{(\Phi,\Psi)}\mcal{B}_{n,r}$ means that $A$ must satisfy
\begin{equation}
    0 = (Y - \Phi A^T)^T \Phi + \Psi^T \left(X - \Psi A\right).
\end{equation}
Rearranging, yields the Sylvester equation
\begin{equation}
    A(\Phi^T\Phi) + (\Psi^T\Psi)A = Y^T\Phi + \Psi^T X.
\end{equation}
Since $(\Phi^T\Phi)$ and $(\Psi^T\Psi)$ are positive-definite, the Sylvester equation has a unique solution (see Theorem~4.4.6 in Horn and Johnson \cite{Horn1991topics}).

\end{proof}

\begin{proof}[Proof of Theorem~\ref{thm:submersion_onto_Bnr}: Over-parametrization]
    Clearly $\Pi_{n,r}$ is surjective because $\Pi_{n,r}(\Phi, \Psi) = (\Phi, \Psi)$ for every $(\Phi, \Psi)\in\mcal{B}_{n,r}$.
    To verify that $\D \Pi_{n,r}$ is surjective we compute its action on a vector $(X,Y) \in \R^{n\times r}\times \R^{n\times r}$ at $(\tilde{\Phi}, \tilde{\Psi})\in D_+(\Pi_{n,r})$.
    Using the derivative of the matrix inverse given by Theorem~3 in Section~8.4 of \cite{Magnus2007matrix}, we obtain
    \begin{multline}
        \D \Pi_{n,r}(\tilde{\Phi}, \tilde{\Psi})(X,Y) \\
        = \big(X (\tilde{\Psi}^T\tilde{\Phi})^{-1} - (\tilde{\Psi}^T\tilde{\Phi})^{-1}(Y^T \tilde{\Phi} + \tilde{\Psi}^T X)(\tilde{\Psi}^T\tilde{\Phi})^{-1} , \ Y \big).
    \end{multline}
    Choosing any $(\hat{X}, \hat{Y}) \in T_{(\Phi, \Psi)}\mcal{B}_{n,r}$ with $(\Phi, \Psi) = \Pi_{n,r}(\tilde{\Phi}, \tilde{\Psi})$, we observe that taking $X=\hat{X}(\tilde{\Psi}^T\tilde{\Phi})$ and $Y = \hat{Y}$ yields
    \begin{multline}
    % \begin{equation}
    % \begin{aligned}
        \D \Pi_{n,r}(\tilde{\Phi}, \tilde{\Psi})(X, Y) \\
        = \big(\hat{X} - (\tilde{\Psi}^T\tilde{\Phi})^{-1}[\hat{Y}^T \tilde{\Phi}(\tilde{\Psi}^T\tilde{\Phi})^{-1} + \tilde{\Psi}^T \hat{X}] , \ \hat{Y} \big) \\
        = \big(\hat{X} - (\tilde{\Psi}^T\tilde{\Phi})^{-1}[\hat{Y}^T \Phi + \Psi^T \hat{X}] , \ \hat{Y} \big)
        = (\hat{X}, \ \hat{Y}).
    % \end{aligned}
    % \end{equation}
    \end{multline}
    Here, the first equality is by substitution, the second equality uses (\ref{eqn:projection_onto_Bnr}), and the third equality follows from the characterization of $T_{(\Phi, \Psi)}\mcal{B}_{n,r}$ in Theorem~\ref{thm:biorthogonal_manifold}.
    Therefore we have proved that $\Pi_{n,r}$ is a surjective submersion.

    By the submersion level set theorem (Corollary~5.13 in \citet{Lee2013introduction}),
    the preimage set $\Pi_{n,r}^{-1}(\Phi, \Psi)$ is a smooth, properly embedded (i.e., closed by Proposition~5.5 in \citet{Lee2013introduction}) submanifold of $D_{+}(\Pi_{n,r})$ with dimension $r^2$.
    Moreover, since $\Pi_{n,r}$ is a global defining function for each fiber $\Pi^{-1}(\Phi, \Psi)$, Proposition~5.38 in \citet{Lee2013introduction} shows that its tangent space is given by
    \begin{equation}
        T_{(\tilde{\Phi}, \tilde{\Psi})} \Pi^{-1}(\Phi, \Psi) 
        = \Null \D \Pi_{n,r}(\tilde{\Phi}, \tilde{\Psi}).
    \end{equation}
    The fiber $\Pi_{n,r}^{-1}(\Phi, \Psi)$ intersects $\mcal{B}_{n,r}$ only at the point $(\Phi, \Psi)$ thanks to the fact that $\Pi_{n,r}$ restricts to the identity on $\mcal{B}_{n,r} \subset D_+(\Pi_{n,r})$.
    Since $\Pi_{n,r}$ is idempotent (when viewed as a map $D_+(\Pi_{n,r}) \to D_+(\Pi_{n,r})$), its tangent map is also idempotent thanks to the chain rule
    \begin{equation}
        \Pi_{n,r} = \Pi_{n,r} \circ \Pi_{n,r} \quad \Rightarrow \quad \D \Pi_{n,r} = \D\Pi_{n,r} \circ \D \Pi_{n,r}.
    \end{equation}
    Hence, $\D \Pi_{n,r}(\Phi, \Psi)$ is a linear projection onto 
    \begin{equation}
        \Range \D \Pi_{n,r}(\Phi, \Psi) = T_{(\Phi, \Psi)}\mcal{B}_{n,r}
    \end{equation}
    with nullspace
    \begin{equation}
        \Null \D \Pi_{n,r}(\Phi, \Psi)
        = T_{(\Phi, \Psi)} \Pi^{-1}(\Phi, \Psi).
    \end{equation}
    Therefore, by the properties of linear projection described in Section~5.9 of Meyer \cite{Meyer2000matrix} these subspaces form a direct sum decomposition
    \begin{equation}
        T_{(\Phi, \Psi)} D_+(\Pi_{n,r}) 
        % = \R^{n\times r}\times\R^{n\times r}
        = T_{(\Phi, \Psi)}\mcal{B}_{n,r} \oplus T_{(\Phi, \Psi)} \Pi^{-1}(\Phi, \Psi),
    \end{equation}
    proving that the intersection is transversal (see \citet{Guillemin1974differential}).

    Thanks to the properties of smooth submersions given by Proposition~4.28 in \citet{Lee2013introduction}, $\Pi_{n,r}$ is a quotient map.
    In particular, $\Pi_{n,r}$ is an open map.
    The coordinate representation of $\Pi_{n,r}$ follows from the rank theorem (Theorem~4.12 in \citet{Lee2013introduction}) and the fact that $\Pi_{n,r}$ is an open map.
    The relationship between coordinate representations for a function $J:\mcal{B}_{n,r}\to\R$ and $J\circ\Pi_{n,r}$ is deduced immediately from the coordinate representation of $\Pi_{n,r}$ supplied by the rank theorem.
    The fact that $J$ is smooth if and only if $J\circ\Pi_{n,r}$ is given by the characteristic property of surjective smooth submersions stated by Theorem~4.29 in \citet{Lee2013introduction}.

    To compute the gradient, we choose $(X,Y)\in T D_+(\Pi_{n,r}) = \R^{n\times r}\times\R^{n\times r}$ and using the definition of the gradient as the Riesz representative of the differential of a function we obtain
    \begin{multline}
    % \begin{equation}
    % \begin{aligned}
        \left\langle \grad ( J \circ \Pi_{n,r})(\tilde{\Phi}, \tilde{\Psi}), \ (X,Y) \right\rangle
        = \D( J \circ \Pi_{n,r})(\tilde{\Phi}, \tilde{\Psi})(X, Y) \\
        = \D J(\Phi, \Psi) \D \Pi_{n,r}(\tilde{\Phi}, \tilde{\Psi})(X, Y) \\
        = \left\langle \grad J(\Phi, \Psi), \ \D \Pi_{n,r}(\tilde{\Phi}, \tilde{\Psi})(X, Y) \right\rangle \\
        = \left\langle \D \Pi_{n,r}(\tilde{\Phi}, \tilde{\Psi})^* \grad J(\Phi, \Psi), \ (X, Y) \right\rangle.
    % \end{aligned}
    % \end{equation}
    \end{multline}
    Since $(X,Y)$ was arbitrary, we have
    \begin{equation}
        \grad ( J \circ \Pi_{n,r})(\tilde{\Phi}, \tilde{\Psi}) = \D \Pi_{n,r}(\tilde{\Phi}, \tilde{\Psi})^* \grad J(\Phi, \Psi).
        \label{eqn:overparametrized_gradient_relation0}
    \end{equation}
    This shows that the gradient of $\grad ( J \circ \Pi_{n,r})(\tilde{\Phi}, \tilde{\Psi})$ is orthogonal to the fiber because it is an element of
    \begin{equation}
        \Range \D \Pi_{n,r}(\tilde{\Phi}, \tilde{\Psi})^* 
        = \big(\Null \D \Pi_{n,r}(\tilde{\Phi}, \tilde{\Psi}) \big)^{\perp}
        = \big( T_{(\Phi, \Psi)} \mcal{B}_{n,r} \big)^{\perp}.
    \end{equation}
    Moreover, since $\D \Pi_{n,r}(\tilde{\Phi}, \tilde{\Psi}):T D_+(\Pi_{n,r}) \to T_{(\Phi, \Psi)}\mcal{B}_{n,r}$ is a surjective linear map between Euclidean spaces, its adjoint $\D \Pi_{n,r}(\tilde{\Phi}, \tilde{\Psi})^*:T_{(\Phi, \Psi)}\mcal{B}_{n,r} \to T D_+(\Pi_{n,r})$ is injective and the Gram operator $G(\tilde{\Phi}, \tilde{\Psi})$ on $T_{(\Phi, \Psi)}\mcal{B}_{n,r}$ is invertible.
    Acting on both sides of (\ref{eqn:overparametrized_gradient_relation0}) with $\D \Pi_{n,r}(\tilde{\Phi}, \tilde{\Psi})$ and inverting the Gram operator appearing on the right hand side yields (\ref{eqn:overparametrized_gradient_relation}).
\end{proof}

\begin{proof}[Proof of Proposition~\ref{prop:overparametrization_domain}]
    We first consider the case when $n=r$.
    Here, the matrices $(\tilde{\Phi},\tilde{\Psi}) \in D_+(\Pi_{n,n})$ are invertible, and so they belong to the general linear group $GL_n$.
    It is well-known that the general linear group has two components: invertible matrices with positive determinant, denoted $GL_n^+$, and invertible matrices with negative determinant, denoted $GL_n^-$.
    Since $\det(\tilde{\Psi}^T\tilde{\Phi}) > 0$ if and only if $\det(\tilde{\Psi})$ and $\det(\tilde{\Phi})$ are nonzero and have the same sign, it follows that
    \begin{equation}
        D_+(\Pi_{n,n}) = (GL_n^+\times GL_n^+) \cup (GL_n^-\times GL_n^-).
    \end{equation}
    These are the components in (\ref{eqn:components_of_domain}), which are obviously disjoint and connected since they are products of connected spaces.

    To show that the domain is connected when $r < n$, we construct a continuous path between points $(\tilde{\Phi}_0,\tilde{\Psi}_0), (\tilde{\Phi}_1,\tilde{\Psi}_1) \in D_+(\Pi_{n,r})$.
    We use the path between two arbitrary biorthogonal matrix pairs $(\Phi_0, \Psi_0), (\Phi_1, \Psi_1) \in \mcal{B}_{n,r}$ constructed in the proof of Theorem.~3.1 in \citet{Otto2022optimizing}.
    By construction, this path lies in $D_+(\Pi_{n,r})$.
    Therefore, it suffices to show that any $(\tilde{\Phi},\tilde{\Psi})\in D_+(\Pi_{n,r})$ can be connected by a continuous path in $D_+(\Pi_{n,r})$ to a biorthogonal pair $(\Phi, \Psi) \in \mcal{B}_{n,r}$.
    
    Since $GL_r^+$ is connected and $\tilde{\Psi}^T \tilde{\Phi}\in GL_r^+$, there is a continuous path $t \mapsto A_t$ in $GL_r^+$ so that $A_0 = I$ and $A_1 = (\tilde{\Psi}^T \tilde{\Phi})^{-1}$.
    Consequently, 
    \begin{equation}
        t \mapsto \big( \tilde{\Phi} A_t, \ \tilde{\Psi} \big)
    \end{equation}
    defines a continuous path in $D_+(\Pi_{n,r})$ from $(\tilde{\Phi},\tilde{\Psi})$ to the biorthogonal pair $(\Phi, \Psi) = \Pi_{n,r}(\tilde{\Phi},\tilde{\Psi})$.
    This completes the proof.
\end{proof}

\section{Properties of the sparsity-promoting penalty}
\label{app:Grassmannian_sparsity}

\begin{proof}[Proof of Theorem~\ref{thm:minimizers_of_sparse_penalty}]
We denote the Stiefel manifold consisting of real $n\times r$ matrices with orthonormal columns by $\mcal{O}_{n,r}$.
It suffices to prove that the minimum value of 
\begin{equation}
    \Vert U\Vert_{1,2} := \sum_{i=1}^n \Vert \row_i(U) \Vert_2
\end{equation}
over $U\in \mcal{O}_{n,r}$ is $r$, and this value is attained if and only if $U$ has exactly $r$ nonzero rows.

The manifold $\mcal{O}_{n,r}$ coincides with the preimage of the regular value $I_r$ under the map $F:\R^{n\times r} \to \R^{r\times r}$ defined by $F: U\mapsto U^T U$, see \cite{Guillemin1974differential}.
We denote the sub-gradient (also known as the sub-differential) of a real-valued convex function $f$ on a Hilbert space $\mcal{H}$ by 
\begin{equation}
    \partial f(x) = \left\{ v\in\mcal{H} \ : \ f(y) \geq f(x) + \left\langle v, \ y - x \right\rangle \right\}.
\end{equation}

Since the Stiefel manifold $\mcal{O}_{n,r}$ is compact and the map $U \mapsto \Vert U \Vert_{1,2}$ is continuous, it attains its minimum value on $\mcal{O}_{n,r}$.
By Theorem~10.8, Theorem~10.45, and Proposition~10.36 in Clarke \cite{Clarke2013functional}, any minimizer $U_*$ satisfies
\begin{equation}
    0 \in \partial\Vert \cdot\Vert_{1,2}(U_*) + \Range\big( \D F(U_*)^* \big),
    \label{eqn:12_norm_optimality_condition_0}
\end{equation}
where $\D F(U_*):\R^{n\times r} \to \R^{r\times r}$ is the linear map $X \mapsto U_*^T X + X^T U_*$ and its adjoint $\D F(U_*)^*:\R^{r\times r} \to \R^{n\times r}$ is given by $A \mapsto U_*(A + A^T)$.
Therefore, we have
\begin{equation}
    \Range\big( \D F(U_*)^* \big) = \left\{ U_*\Lambda \ : \ \Lambda\in\R^{r\times r}, \ \Lambda = \Lambda^T \right\}.
\end{equation}
Letting $\row_i:\R^{n\times r} \to \R^{r}$ denote the linear map extracting the $i$th row a matrix, we obtain
\begin{multline}
    \row_i\big( \partial\Vert \cdot\Vert_{1,2}(U_*) \big) \\
    % = \partial \Vert \cdot \Vert_2\big( \row_i(U_*) \big)
    = \left\{ \begin{matrix}
    \row_i(U_*)/\Vert \row_i(U_*) \Vert_2, & \row_i(U_*) \neq 0 \\
    \{ v\in\R^r : \Vert v \Vert_2 \leq 1 \}, & \row_i(U_*)=0
    \end{matrix}\right.,
\end{multline}
for $i=1, \ldots,n$ thanks to Theorems~10.13~and~10.19 in \citet{Clarke2013functional} and the well-known expression for the subdifferential of the Euclidean norm.
% Letting $\row_i:\R^{n\times r} \to \R^{r}$ denote the linear map extracting the $i$th row a matrix, Theorems~10.13~and~10.19 in \citet{Clarke2013functional} yield
% \begin{equation}
%     \partial\Vert \cdot\Vert_{1,2}(U_*) = \sum_{i=1}^n \row_i^* \partial \Vert \cdot \Vert_2\big( \row_i(U_*) \big),
% \end{equation}
% where $\row_i^* v$ places the vector $v$ into the $i$th row of an $n\times r$ matrix whose remaining rows are equal to $0$.
% Using the well-known expression for the subdifferential of the Euclidean norm we obtain
% \begin{equation}
%     \row_i\left( \partial\Vert \cdot\Vert_{1,2}(U_*) \right) = 
%     \left\{ \begin{matrix}
%     \row_i(U_*)/\Vert \row_i(U_*) \Vert_2, & \row_i(U_*) \neq 0 \\
%     \{ v\in\R^r : \Vert v \Vert_2 \leq 1 \}, & \row_i(U_*)=0
%     \end{matrix}\right..
% \end{equation}
Letting
\begin{equation}
    D := \begin{bmatrix}
        \Vert \row_1(U_*) \Vert_2 & & \\
        & \ddots & \\
        & & \Vert \row_n(U_*) \Vert_2
    \end{bmatrix},
\end{equation}
the optimality condition (\ref{eqn:12_norm_optimality_condition_0}) yields a symmetric matrix $\Lambda \in \R^{r\times r}$ such that
\begin{equation}
    D U_* = U_* \Lambda.
    \label{eqn:12_norm_optimality_condition_1}
\end{equation}
The matrix $\Lambda = U_* ^T D U_*$ is positive-definite because rescaling the nonzero rows of $U_*$ does not change its rank.

We show that $U_*$ has exactly $r$ nonzero rows by contradiction.
We employ the following technical lemma, whose proof relies on the optimality condition (\ref{eqn:12_norm_optimality_condition_1}).
\begin{lemma}
If $U_*$ has greater than $r$ nonzero rows then $U_*$ has a linearly dependent set of rows with identical nonzero Euclidean norms.
\end{lemma}
\begin{proof}
Let $\tilde{U} \in \R^{(r+1)\times r}$ be a rank-$r$ matrix composed of nonzero rows of $U_*$ and let $\tilde{D}$ be the diagonal matrix with the corresponding entries of $D$.
Then $\tilde{U}^T$ has a one-dimensional nullspace spanned by a vector $\tilde{v} \in \R^{r+1}$.
Moreover, $\tilde{v}$ has at least two nonzero entries, for if it had only a single nonzero entry then the corresponding row of $\tilde{U}$ would be zero, a contradiction.
Since $\Lambda$ is invertible, we have $\tilde{D}\tilde{v} \in \Null(\tilde{U}^T)$ and so $\tilde{D}\tilde{v} = \lambda \tilde{v}$ for some $\lambda \in \R$.
Since $\tilde{D}$ is nonsingular, $\tilde{D}\tilde{v}\neq 0$ and so $\lambda \neq 0$.
Therefore, the entries of the diagonal matrix $\tilde{D}$ corresponding to nonzero entries of $\tilde{v}$ (of which there are at least two) are all equal to $\lambda \neq 0$.
Since the diagonal entries of $\tilde{D}$ are the magnitudes of rows of $U_*$, we have identified a collection of at least two linearly dependent rows of $U_*$ with identical Euclidean norms.
\end{proof}

Let the rows of $U_1$ be a set of $q$ linearly dependent rows of $U_*$ with identical nonzero Euclidean norms, as provided by the lemma.
Letting $U_2$ contain the remaining rows of $U_*$ we observe that
\begin{equation}
    \Vert U_* \Vert_{1,2} = \Vert U_1 \Vert_{1,2} + \Vert U_2 \Vert_{1,2}.
\end{equation}
Let $v\in \Null(U_1^T)$ have unit magnitude, and let $Q\in \R^{q\times q}$ be a unitary matrix whose last row is $v$.
We observe that 
\begin{equation}
    U := \begin{bmatrix}
    Q U_1 \\
    U_2
    \end{bmatrix} \in \mcal{O}_{n,r},
\end{equation}
where the last row of $Q U_1$ is zero by construction of $Q$.
Moreover, denoting $r_i = \Vert \row_i(Q U_1) \Vert_2$, we have
\begin{equation}
    \sum_{i=1}^q r_i^2 = \Vert Q U_1 \Vert_F^2 = \Vert U_1 \Vert_F^2 = \frac{1}{q} \Vert U_1 \Vert_{1,2}^2,
\end{equation}
where the first equality is by definition of the Frobenius norm, the second equality is by unitary invariance of the Frobenius norm, and the third equality holds because the rows of $U_1$ have identical Euclidean norms.
Since $r_q = 0$ by construction of $Q$, the Cauchy-Schwarz inequality gives
\begin{multline}
    \Vert Q U_1 \Vert_{1,2} 
    = \sum_{i=1}^q r_i 
    \leq \sqrt{q-1}\sqrt{\sum_{i=1}^q r_i^2} \\
    = \sqrt{\frac{q-1}{q}} \Vert U_1 \Vert_{1,2} 
    < \Vert U_1 \Vert_{1,2}.
\end{multline}
Therefore, we have constructed a matrix $U \in \mcal{O}_{n,r}$ such that
\begin{equation}
    \Vert U \Vert_{1,2} 
    = \Vert Q U_1 \Vert_{1,2} + \Vert U_2 \Vert_{1,2}
    < \Vert U_1 \Vert_{1,2} + \Vert U_2 \Vert_{1,2}
    = \Vert U_* \Vert_{1,2},
\end{equation}
contradicting the optimality of $U_*$.
Therefore, we have shown that any optimizer $U_*$ has exactly $r$ nonzero rows.

Since the $r\times r$ submatrix of $U_*$ formed by its nonzero rows is unitary, it follows that these rows all have unit magnitude, and so $\Vert U_* \Vert_{1,2} = r$, completing the proof of Theorem~\ref{thm:minimizers_of_sparse_penalty}.
\end{proof}

Let $\mcal{G}_{n,r}$ denote the Grassmann manifold of $r$-dimensional subspaces of $\R^n$ \cite{Bendokat2020grassmann, Absil2004riemannian, Wong1967differential}.
This manifold can be viewed as the quotient of $\mcal{O}_{n,r}$ under the free and proper action of the orthogonal group $\mcal{O}_{r,r}$ with canonical projection $\pi:U \mapsto \Range(U)$.
Thanks to the quotient manifold theorem, the canonical projection is a surjective submersion (see Theorem~21.10 in \citet{Lee2013introduction}).
This allows us to define $\Vert \cdot \Vert_{1,2}$ on $\mcal{G}_{n,r}$ using orthogonal matrix representatives $U\in\mcal{O}_{n,r}$ according to
\begin{equation}
    \Vert \Range(U) \Vert_{1,2} := \Vert U \Vert_{1,2}.
    \label{eqn:definition_of_12fun_on_Grassmann}
\end{equation}
Similarly, we define $\Vert \cdot \Vert_0$ on $\mcal{G}_{n,r}$ to be the number of nonzero rows in any matrix representative, that is
\begin{equation}
    \Vert \Range(U) \Vert_{0} := \Vert U \Vert_{0,2}.
    \label{eqn:definition_of_0fun_on_Grassmann}
\end{equation}
These are well-defined functions on $\mcal{G}_{n,r}$ because they do not depend on the choice of orthogonal representative thanks to rotational invariance from the right, namely that $\Vert U \Vert_{1,2} = \Vert U Q \Vert_{1,2}$ for any $Q\in\mcal{O}_{r,r}$.
% Because of this rotational invariance, Theorem~\ref{thm:minimizers_of_sparse_penalty} immediately yields
% \begin{corollary}
%     \label{cor:global_minimizer_of_12fun}
%     The minimum value of $\Vert \mcal{V} \Vert_{1,2}$ for $\mcal{V}\in \mcal{G}_{n,r}$ is $r$, and a subspace $\mcal{V}\in\mcal{G}_{n,r}$ attains this value $\Vert \mcal{V} \Vert_{1,2} = r$ if and only if there are $r$ row indices among $\{1, 2, \ldots, n \}$ so that every matrix representative $\Phi \in \R_*^{n,r}$ with $\Range(\Phi) = \mcal{V}$ has nonzero entries in each of these rows and the remaining rows are identically zero.
% \end{corollary}

The penalty function defined by (\ref{eqn:sparsity_penalty}) can be expressed as $ R_{1,2}: \mcal{V} \mapsto \Vert \mcal{V} \Vert_{1,2} - r $.
The following lemma characterizes the behavior of this function in the neighborhood of a global minimizer on $\mcal{G}_{n,r}$.
Before stating the result, we provide some machinery needed to study the local geometry of the Grassmannian --- namely the Riemannian metric and notion or horizontal lift for tangent vectors.
For more details, we refer to \cite{Bendokat2020grassmann, Absil2004riemannian}.
Endowing $\mcal{O}_{n,r}$ with the Riemannian metric
\begin{equation}
    \left\langle X, \ Y \right\rangle_{U} = \Tr\left( X^T Y \right),
\end{equation}
the horizontal space $\mcal{H}_{U}\subset T_{U}\mcal{O}_{n,r}$ is defined to be orthogonal to the fiber $\pi^{-1}(U)$ and is given by
\begin{equation}
    \mcal{H}_U = \left\{ X \in \R^{n\times r} \ : \ U^T X = 0 \right\}.
    \label{eqn:horizontal_space}
\end{equation}
Given $U\in \mcal{O}_{n,r}$ and a tangent vector $\xi\in T_{\Range(U)}\mcal{G}_{n,r}$, there is a unique element $\bar{\xi}_U \in \mcal{H}_U$ so that $\D \pi(U) \bar{\xi}_{U} = \xi$.
This element is called the horizontal lift of $\xi$ at $U$.
These elements obey the transformation law
\begin{equation}
    \bar{\xi}_{UQ} = \bar{\xi}_{U} Q
    \label{eqn:transformation_law_for_horizontal_lifts}
\end{equation}
for every $Q\in \mcal{O}_{r,r}$, which enables us to endow $\mcal{G}_{n,r}$ with the inner product
\begin{equation}
    \left\langle \xi, \ \eta \right\rangle_{\Range(U)}
    := \left\langle \bar{\xi}_U,\ \bar{\eta}_{U} \right\rangle_{U}
    = \Tr\left( (\bar{\xi}_{U})^T \bar{\eta}_U \right)
    \label{eqn:metric_on_Grassmannian}
\end{equation}
which is independent of the choice of matrix representative.

\begin{lemma}
    \label{lem:local_behavior_of_12fun}
    Let $U_*\in\mcal{O}_{n,r}$ have exactly $r$ nonzero rows and let $\exp_{\Range(U_*)}: T_{\Range(U_*)}\mcal{G}_{n,r} \to \mcal{G}_{n,r}$ denote the exponential map on the Grassmannian.
    Then for every $\xi \in T_{\Range(U_*)}\mcal{G}_{n,r}$ we have
    \begin{multline}
        \left\vert  R_{1,2}\big( \exp_{\Range(U_*)}(\xi) \big) - \big\Vert \bar{\xi}_{U_*} \big\Vert_{1,2} \right\vert \\
        %\leq C \Vert \xi \Vert_{\Range(U_*)}^2.
        \leq \sqrt{n}\left( \frac{1}{2} \Vert \xi \Vert_{\Range(U_*)}^2 + \frac{1}{6} \Vert \xi \Vert_{\Range(U_*)}^3 \right).
    \end{multline}
\end{lemma}
\begin{proof}
Letting $\bar{\xi}_{U_*} = U \Sigma V^T$ be a singular value decomposition with $V V^T = I$, Theorem~2.3 in Edelman et al. \cite{Edelman1998geometry} shows that
\begin{equation}
    \Phi = U_* V \cos(\Sigma) V^T + U \sin(\Sigma) V^T \in \mcal{O}_{n,r}
    \label{eqn:Grassmann_exponential_formula}
\end{equation}
satisfies $\Range(\Phi) = \exp_{\Range(U_*)}(\xi)$, and is thus an orthonormal representative of the exponential.
By (\ref{eqn:definition_of_12fun_on_Grassmann}), we automatically have $\big\Vert \exp_{\Range(U_*)}(\xi) \big\Vert_{1,2} = \Vert \Phi \Vert_{1,2}$.

By definition of the horizontal space in (\ref{eqn:horizontal_space}), we have $U_*^T \bar{\xi}_{U_*} = 0$.
Because the nonzero rows of $U_*$ form an invertible sub-matrix of $U_*$, we must have $\row_i(\bar{\xi}_{U_*}) = 0$ whenever $\row_i(U_*) \neq 0$.
Since $V$ is invertible, $\row_i(\bar{\xi}_{U_*}) = 0$ implies that $\row_i(U \Sigma) = 0$, which, in turn, implies that $\row_i(U \sin(\Sigma)) = 0$.
Consequently, the two terms in the sum (\ref{eqn:Grassmann_exponential_formula}) have disjoint sets of nonzero rows and so
\begin{multline}
    \big\Vert \exp_{\Range(U_*)}(\xi) \big\Vert_{1,2}
    = \big\Vert U_* V \cos(\Sigma) V^T \big\Vert_{1,2} \\ 
    + \big\Vert U \sin(\Sigma) V^T \big\Vert_{1,2}.
    \label{eqn:exponential_norm_decomp_at_global_min}
\end{multline}
Using the above together with the definition of the penalty function in (\ref{eqn:sparsity_penalty}) and the triangle inequality gives
\begin{multline}
    \left\vert  R_{1,2}\big( \exp_{\Range(U_*)}(\xi) \big) - \big\Vert \bar{\xi}_{U_*} \big\Vert_{1,2} \right\vert \\
    \leq \left\vert \big\Vert U_* V \cos(\Sigma) V^T \big\Vert_{1,2} - r \right\vert \\
    + \left\vert \big\Vert U \sin(\Sigma) V^T \big\Vert_{1,2} - \big\Vert \bar{\xi}_{U_*} \big\Vert_{1,2} \right\vert.
    \label{eqn:sparsity_penalty_bound_at_min}
\end{multline}
We bound the two terms appearing on the right-hand side of this inequality.

Since $\Vert \cdot \Vert_{1,2}$ is a matrix norm, it satisfies the triangle inequality and we have
\begin{equation}
    \Big\vert \big\Vert U_* V \cos(\Sigma) V^T \big\Vert_{1,2} - \underbrace{\big\Vert U_* \big\Vert_{1,2}}_{r} \Big\vert
    \leq \left\Vert U_* V \big(I -\cos(\Sigma)\big) V^T \right\Vert_{1,2}.
\end{equation}
For any $n\times r$ matrix $M$, we have $\Vert M \Vert_{1,2} \leq \sqrt{n} \Vert M \Vert_{F}$ by the Cauchy-Schwartz inequality.
% \begin{equation}
%     \Vert M \Vert_{1,2} = \sum_{i=1}^n \Vert \row_i(M) \Vert_2 \leq \sqrt{n} \sqrt{\sum_{i=1}^n \Vert \row_i(M) \Vert_2^2 } = \sqrt{n} \Vert M \Vert_{F}
% \end{equation}
Using this and the fact that $U_*$ is an isometry and $V$ is unitary the we obtain
\begin{equation}
    \Big\vert \big\Vert U_* V \cos(\Sigma) V^T \big\Vert_{1,2} - r \Big\vert \leq \sqrt{n} \big\Vert I - \cos(\Sigma) \big\Vert_F.
\end{equation}
From the well-known inequality $\vert 1 - \cos(x) \vert \leq x^2 / 2$ (resulting from the double angle formula and $\vert \sin(x) \vert \leq \vert x \vert$) we readily obtain
\begin{equation}
    \Big\vert \big\Vert U_* V \cos(\Sigma) V^T \big\Vert_{1,2} - r \Big\vert \leq \frac{\sqrt{n}}{2} \Vert \Sigma \Vert_F^2 = \frac{\sqrt{n}}{2} \Vert \xi \Vert_{\Range(U_*)}^2.
\end{equation}
By a similar argument, we obtain
\begin{equation}
    \Big\vert \big\Vert U \sin(\Sigma) V^T \big\Vert_{1,2} - \big\Vert \bar{\xi}_{U_*} \big\Vert_{1,2} \Big\vert
    \leq \sqrt{n} \Vert \Sigma - \sin(\Sigma) \Vert_F.
\end{equation}
The inequality
\begin{equation}
     x - \sin(x) = \int_{0}^x \big( 1 - \cos(t) \big) \td t \leq  \frac{1}{2} \int_{0}^x t^2 \td t = \frac{1}{6} x ^3, 
\end{equation}
for $x > 0$ then gives
\begin{equation}
    \Big\vert \big\Vert U \sin(\Sigma) V^T \big\Vert_{1,2} - \big\Vert \bar{\xi}_{U_*} \big\Vert_{1,2} \Big\vert 
    \leq \frac{\sqrt{n}}{6} \Vert \Sigma \Vert_F^3 
    = \frac{\sqrt{n}}{6} \Vert \xi \Vert_{\Range(U_*)}^3.
\end{equation}
Combining these inequalities with (\ref{eqn:sparsity_penalty_bound_at_min}) completes the proof.
\end{proof}

\begin{corollary}
    \label{cor:local_behavior_of_12fun}
    Consider $\mcal{G}_{n,r}$ with fixed dimensions $n,r$ and let $\varepsilon > 0$.
    Then there is a $\delta_{\varepsilon} > 0$ so that for every $U_*\in\mcal{O}_{n,r}$ having exactly $r$ nonzero rows and every $\xi \in T_{\Range(U_*)}\mcal{G}_{n,r}$ with $\Vert \xi \Vert_{\Range(U_*)} \leq \delta_{\varepsilon}$ we have
    \begin{equation}
        (1 - \varepsilon) \big\Vert \bar{\xi}_{U_*} \big\Vert_{1,2} 
        \leq  R_{1,2}\big( \exp_{\Range(U_*)}(\xi) \big)
        \leq (1 + \varepsilon) \big\Vert \bar{\xi}_{U_*} \big\Vert_{1,2}.
    \end{equation}
\end{corollary}

% \begin{theorem}
%     \label{thm:solutions_of_12_regularized_minimization}
%     Let $\cal{M}$ be a smooth manifold and let $D(J_0)$ be an open subset of $\mcal{M}\times \mcal{G}_{n,r}$ on which a real non-negative-valued function $J_0$ is defined and continuously differentiable.
%     Suppose that there is a finite constant $M$ so that the preimage set $S_M = \{ (x,\mcal{V}) \in D(J_0) \ : \ J_0(x,\mcal{V}) \leq M \}$ is compact and contains a point $(x_0, \mcal{V}_0)$ so that $\Vert \mcal{V}_0 \Vert_{0} = r$.
%     Then for any $\gamma \geq 0$, the function on $D(J_0)$ defined by
%     \begin{equation}
%         J_{\gamma}(x,\mcal{V}) = J_0(x,\mcal{V}) + \gamma \Vert \mcal{V}\Vert_{1,2}
%     \end{equation}
%     attains its minimum and all such minimizers lie in $S_M$.
%     Furthermore, there is a constant $\Gamma \geq 0$ so that when $\gamma > \Gamma$, every minimizer $(x_*, \mcal{V}_*)$ of $J_{\gamma}$ satisfies $\Vert \mcal{V}_* \Vert_0 = r$.
% \end{theorem}
\begin{proof}[Proof of Theorem~\ref{thm:solutions_of_12_regularized_minimization}]
    The function $\Vert \cdot \Vert_{1,2}$ defined by (\ref{eqn:definition_of_12fun_on_Grassmann}) is continuous on $\mcal{G}_{n,r}$ thanks to Theorem~4.29 in \citet{Lee2013introduction}.
    This makes $J_{\gamma}$ continuous.
    We consider an infemizing sequence $p_k = (x_k, \mcal{V}_k)\in D(J_0)$, $k=0,1, \ldots$ so that $J_{\gamma}(p_0) \geq J_{\gamma}(p_1) \geq \cdots$ and 
    \begin{equation}
        \lim_{k\to\infty} J_{\gamma}(p_k) = \inf_{p\in D(J_0)} J_{\gamma}(p).
    \end{equation}
    Since the sequence begins with $p_0 = (x_0, \mcal{V}_0)\in S_M$ where $ R_{1,2}( \mcal{V}_0 ) = 0$, we have
    \begin{equation}
        J_0(p_k) 
        = J_{\gamma}(p_k) - \gamma  R_{1,2}( \mcal{V}_k ) 
        \leq J_{\gamma}(p_0) = J_0(p_0) 
        \leq M,
    \end{equation}
    which implies that $\{p_k\}_{k=0}^{\infty} \subset S_M$.
    Since $S_M$ is compact, we may pass to a convergent subsequence $p_{k_l} \to p_* \in S_K$.
    The limit of this subsequence is a minimizer of $J_{\gamma}$ because
    \begin{equation}
        J_{\gamma}(p_*) 
        = \lim_{l\to\infty} J_{\gamma}(p_{k_l}) 
        = \lim_{k\to\infty} J_{\gamma}(p_{k}) 
        = \inf_{p\in D(J_0)} J_{\gamma}(p),
    \end{equation}
    which follows from continuity of $J_{\gamma}$.
    This establishes the existence of a minimizer of $J_{\gamma}$ in $S_M$.
    Moreover, every minimizer $\tilde{p}_* = (\tilde{x}_*, \tilde{\mcal{V}}_*) \in D(J_0)$ of $J_{\gamma}$ belongs to $S_M$ because
    \begin{equation}
        J_0(\tilde{p}_*) 
        = J_{\gamma}(\tilde{p}_*) - \gamma  R_{1,2}( \tilde{\mcal{V}}_* )
        \leq J_{\gamma}(p_0)
        = M.
    \end{equation}
    
    To produce the required $\Gamma$, we first find $\Gamma_0 > 0$ so that any minimizer of $J_{\gamma}$ with $\gamma \geq \Gamma_0$ lies within a neighborhood of a minimizer of $ R_{1,2}$, whose local behavior is described by Corollary.~\ref{cor:local_behavior_of_12fun}.
    We then find $\Gamma_1 > 0$ so that if $\gamma > \Gamma_1$ then the rate of increase in $\gamma  R_{1,2}$ as one moves away from a minimizer of $ R_{1,2}$ dominates the gradient of $J_0$ in the neighborhood.
    Because of this, the only possible minimizers of $J_{\gamma}$ will also be minimizers of $ R_{1,2}$.
    
    We use the following two lemmas.
    The first lemma, stated below, shows that by making $ R_{1,2}(\mcal{V}_*)$ close to $0$, we also force $\mcal{V}_*$ to be close to some $\tilde{\mcal{V}}_*\in\mcal{G}_{n,r}$ with $\Vert \tilde{\mcal{V}}_* \Vert_0 = r$.
    \begin{lemma}
        \label{lem:small_12_implies_clost_to_sparse}
        For every $\varepsilon > 0$ there is a constant $\delta_{\varepsilon} > 0$ so that every $\mcal{V} \in \mcal{G}_{n,r}$ with $ R_{1,2}( \mcal{V}) < \delta_{\varepsilon}$ lies within an $\varepsilon$-normal neighborhood
        \begin{equation}
            B_{\varepsilon}(\tilde{\mcal{V}}_*) = \left\{ \exp_{\tilde{\mcal{V}}_*}(\xi) \ : \ \Vert \xi \Vert_{\tilde{\mcal{V}}_*} < \varepsilon \right\}
        \end{equation}
        of some $\tilde{\mcal{V}}_* \in \mcal{G}_{n,r}$ with $\Vert \tilde{\mcal{V}}_*\Vert_0 = r$.
    \end{lemma}
    \begin{proof}
        Suppose the contrary.
        Then there is an $\varepsilon > 0$ and a sequence of subspaces $\{\mcal{V}_k\}_{k=1}^{\infty} \subset \mcal{G}_{n,r}$ so that 
        \begin{equation}
             R_{1,2}( \mcal{V}_k ) < n^{-1}
        \end{equation}
        and $\mcal{V}_k \notin B_{\varepsilon}(\mcal{V}_*)$ for every $\tilde{\mcal{V}}_*\in\mcal{G}_{n,r}$ with $\Vert \tilde{\mcal{V}}_* \Vert_0 = r$.
        without loss of generality, we assume that $\varepsilon$ is smaller than the injectivity radius of $\mcal{G}_{n,r}$ (which is $\pi/2$, see Bendokat et al. \cite{Bendokat2020grassmann}) so that each $B_{\varepsilon}(\tilde{\mcal{V}}_*)$ is open.
        Since $\mcal{G}_{n,r}$ is a compact manifold and
        \begin{equation}
            B = \bigcup_{\substack{\tilde{\mcal{V}}_*\in\mcal{G}_{n,r} : \\ \Vert \tilde{\mcal{V}}_* \Vert_0 = r}} B_{\varepsilon}(\tilde{\mcal{V}}_*)
        \end{equation}
        is open, the set $K = \mcal{G}_{n,r} \setminus B$ is compact and contains each $\mcal{V}_k$ in the sequence.
        Passing to a convergent subsequence $\mcal{V}_{k_l} \to \bar{\mcal{V}} \in K$ and using the continuity of $ R_{1,2}$ on $\mcal{G}_{n,r}$ we would have
        \begin{equation}
             R_{1,2} ( \bar{\mcal{V}} ) = \lim_{l\to\infty}  R_{1,2}( \mcal{V}_{k_l} ) = 0,
        \end{equation}
        but $\Vert \bar{\mcal{V}} \Vert_0 > r$, contradicting Theorem~\ref{thm:minimizers_of_sparse_penalty}.
    \end{proof}
    We must also ensure that every $(x_*, \mcal{V})$ within a sufficiently small neighborhood of $(x_*, \tilde{\mcal{V}}_*)$ lies within the domain of the objective $D(J_0)$.
    This will be accomplished using the following lemma and the fact that $(x_*, \mcal{V}_*)\in S_M$, a compact set.
    \begin{lemma}
        \label{lem:S_M_padding}
        For sufficiently small $\varepsilon_0 > 0$, we have
        \begin{equation}
            \{x\}\times \overline{B_{\varepsilon_0}(\mcal{V})} \subset D(J_0) \qquad \forall (x,\mcal{V}) \in S_M.
        \end{equation}
    \end{lemma}
    \begin{proof}
        Let $d:\mcal{G}_{n,r}\times\mcal{G}_{n,r}\to [0,\infty)$ denote the geodesic distance on the Grassmannian.
        If the set
        \begin{equation}
            K = \left\{ \big( (x, \mcal{V}_1), \mcal{V}_2 \big) \in S_M \times \mcal{G}_{n,r} \ : \ (x,\mcal{V}_2) \notin D(J_0) \right\},
        \end{equation}
        is empty, then $\{ x\}\times \mcal{G}_{n,r} \subset D(J_0)$ for every $(x,\mcal{V}) \in S_M$.
        Otherwise, if $K$ is nonempty, it suffices to show that
        \begin{equation}
            0 < \inf_{((x, \mcal{V}_1),\ \mcal{V}_2)\in K} d(\mcal{V}_1, \mcal{V}_2).
        \end{equation}
        We proceed by showing that $K$ is compact.
        To see this, consider the continuous function $h:S_M\times\mcal{G}_{n,r} \to (\mcal{M}\times\mcal{G}_{n,r})\times (\mcal{M}\times\mcal{G}_{n,r})$ defined by
        \begin{equation}
            h:\big( (x, \mcal{V}_1), \mcal{V}_2 \big) \mapsto \big( (x, \mcal{V}_1), (x, \mcal{V}_2) \big).
        \end{equation}
        We observe that $K = g^{-1}(S_M\times D(J_0)^c)$, where $S_M\times D(J_0)^c$ is closed.
        This implies that $K$ is also closed, and since $K$ is a closed subset of the compact set $S_M \times \mcal{G}_{n,r}$, it follows that $K$ is compact.
        Thus, there is a point $\big( (\bar{x}, \bar{\mcal{V}}_1), \bar{\mcal{V}}_2\big) \in K$ so that
        \begin{equation}
            d(\bar{\mcal{V}}_1, \bar{\mcal{V}}_2) = \inf_{((x, \mcal{V}_1),\ \mcal{V}_2)\in K} d(\mcal{V}_1, \mcal{V}_2).
        \end{equation}
        If $d(\bar{\mcal{V}}_1, \bar{\mcal{V}}_2) = 0$ then $\bar{\mcal{V}}_1 = \bar{\mcal{V}}_2$ and we have $(\bar{x}, \bar{\mcal{V}}_2) \in S_M \subset D(J_0)$, a contradiction.
        Therefore, $d(\bar{\mcal{V}}_1, \bar{\mcal{V}}_2) > 0$, and we may take $\varepsilon_0 > 0$ so that $\varepsilon_0 < d(\bar{\mcal{V}}_1, \bar{\mcal{V}}_2)$.
    \end{proof}
    
    To carry out the plan described above, we let $\delta_{1/2}$ be determined by Corollary.~\ref{cor:local_behavior_of_12fun}, $\varepsilon_0$ be determined by Lemma.~\ref{lem:S_M_padding}, and we choose $\varepsilon = \min\{ \delta_{1/2}, \varepsilon_0 / 2, \pi/2 \}$, where we recall that $\pi/2$ is the injectivity radius of the exponential map on $\mcal{G}_{n,r}$ \cite{Bendokat2020grassmann}. 
    We then take $\delta_{\varepsilon}$ provided by Lemma.~\ref{lem:small_12_implies_clost_to_sparse}.
    Since $S_M$ is compact,
    \begin{equation}
        m := \min_{p\in S_M} J_0(p),
    \end{equation}
    is finite, and so we take
    \begin{equation}
        \Gamma_0 = \frac{M-m}{\delta_{\varepsilon}}.
    \end{equation}
    Let $\gamma > \Gamma_0$ and let $p_* = (x_*, \mcal{V}_*) \in D(J_0)$ be a minimizer of $J_{\gamma}$.
    By the first part of the theorem, $p_* \in S_M$, and so we have
    \begin{multline}
         R_{1,2}( \mcal{V}_* ) 
        = \frac{1}{\gamma}\left[ J_{\gamma}(p_*) - J_0(p_*) \right] \\
        \leq \frac{1}{\gamma}\left[ J_{\gamma}(p_0) - m \right]
        \leq \frac{1}{\gamma}\left[ M - m \right]
        < \delta_{\varepsilon}.
    \end{multline}
    By Lemma.~\ref{lem:small_12_implies_clost_to_sparse} there is an element $\tilde{\mcal{V}}_* \in \mcal{G}_{n,r}$ with $\Vert \tilde{\mcal{V}}_* \Vert_{0} = r$ so that $\mcal{V}_* \in B_{\varepsilon}(\tilde{\mcal{V}}_*)$.
    First, we claim that $\{ x_* \} \times B_{\varepsilon}(\tilde{\mcal{V}}_*) \subset D(J_0)$.
    To see this, we let $d:\mcal{G}_{n,r}\times\mcal{G}_{n,r}\to [0,\infty)$ denote the geodesic distance on the Grassmannian and we observe that for any $\mcal{V} \in B_{\varepsilon}(\tilde{\mcal{V}}_*)$, we have
    \begin{equation}
        d(\mcal{V}, \mcal{V}_*) \leq d(\mcal{V}, \tilde{\mcal{V}}_*) + d(\tilde{\mcal{V}}_*, \mcal{V}_*) < \frac{\varepsilon_0}{2} + \frac{\varepsilon_0}{2} = \varepsilon_0.
    \end{equation}
    Therefore, we have $B_{\varepsilon}(\tilde{\mcal{V}}_*) \subset B_{\varepsilon_0}(\mcal{V}_*)$.
    Since $(x_*, \mcal{V}_*) \in S_M$, Lemma.~\ref{lem:S_M_padding} gives
    \begin{equation}
        \{ x_* \} \times B_{\varepsilon}(\tilde{\mcal{V}}_*) \subset \{ x_* \} \times B_{\varepsilon_0}(\mcal{V}_*) \subset D(J_0).
    \end{equation}
    
    Furthermore, since $\varepsilon \leq \delta_{1/2}$ and $\mcal{V}_* \in B_{\varepsilon}(\tilde{\mcal{V}_*})$, we can use Corollary.~\ref{cor:local_behavior_of_12fun} to estimate $ R_{1,2}(\mcal{V}_* )$.
    Since $\varepsilon$ is less than the injectivity radius of the exponential map on $\mcal{G}_{n,r}$, there is a unique $\xi \in T_{\tilde{\mcal{V}}_*}\mcal{G}_{n,r}$ with $\Vert \xi \Vert_{\tilde{\mcal{V}}_*} < \varepsilon$ so that $\mcal{V}_* = \exp_{\tilde{\mcal{V}}_*}(\xi)$.
    Letting $\tilde{U}_* \in \mcal{O}_{n,r}$ denote an orthogonal representative of $\tilde{\mcal{V}}_* = \Range(\tilde{U}_*)$, Corollary.~\ref{cor:local_behavior_of_12fun} gives
    \begin{equation}
         R_{1,2}( \mcal{V}_* ) 
        \geq \frac{1}{2} \Vert \bar{\xi}_{\tilde{U}_*} \Vert_{1,2}.
    \end{equation}
    By concavity of the square root, we have $\Vert \bar{\xi}_{\tilde{U}_*} \Vert_{1,2} \geq \Vert \bar{\xi}_{\tilde{U}_*} \Vert_{F} = \Vert \xi \Vert_{\tilde{\mcal{V}}_*}$, and so
    \begin{equation}
         R_{1,2}(\mcal{V}_*) \geq \frac{1}{2} \Vert \xi \Vert_{\tilde{\mcal{V}}_*}.
        \label{eqn:estimate_of_12fun_from_below}
    \end{equation}
    The set
    \begin{equation}
        K = \bigcup_{(x, \mcal{V}) \in S_M} \{ x \} \times \overline{B_{\varepsilon_0}(\mcal{V})}
    \end{equation}
    is compact by a simple subsequence argument.
    Since $K \subset D(J_0)$ by Lemma.~\ref{lem:S_M_padding}, and $J_0$ is continuously differentiable on $D(J_0)$, the constant
    \begin{equation}
        L = \sup_{(x,\mcal{V}) \in K} \Vert \grad_{\mcal{V}} J_0(x, \mcal{V}) \Vert_{\mcal{V}}
    \end{equation}
    is finite.
    Here, $\grad_{\mcal{V}} J_0(x, \mcal{V}) \in T_{\mcal{V}}\mcal{G}_{n,r}$ denotes the gradient of $J_0$ with respect to its second argument at the point $(x,\mcal{V})\in\mcal{M}\times\mcal{G}_{n,r}$.
    Since $\{x_*\} \times B_{\varepsilon}(\tilde{V}_*) \subset \{ x_* \} \times B_{\varepsilon_0}(\mcal{V}_*) \subset K$, applying the fundamental theorem of calculus along the geodesic curve $c(t) = \exp_{\tilde{\mcal{V}}_*}(t\xi)$, $t\in [0,1]$ yields
    \begin{multline}
        \vert J_0(x_*, \mcal{V}_*) - J_0(x_*, \tilde{\mcal{V}}_*) \vert \\
        \leq \int_{0}^{1} \left\vert \big\langle \grad_{\mcal{V}} J_0 (x_*, c(t)),\ c'(t) \big\rangle_{c(t)} \right\vert \td t
        \leq L \Vert \xi \Vert_{\tilde{\mcal{V}}_*}.
        \label{eqn:estimate_of_J0_from_FTC}
    \end{multline}
    Let
    \begin{equation}
        \Gamma_1 = 2 L,
    \end{equation}
    and let $(x_*, \mcal{V}_*)$ be a minimizer of $J_\gamma$ with $\gamma > \max\{ \Gamma_0, \Gamma_1 \}$.
    To prove the theorem, it suffices to show that $\xi = 0$, for this will imply that $\mcal{V}_* = \tilde{\mcal{V}}_*$ with $\Vert \tilde{\mcal{V}}_*\Vert_0 = r$.
    If $\xi \neq 0$ then combining our estimates (\ref{eqn:estimate_of_12fun_from_below}) and (\ref{eqn:estimate_of_J0_from_FTC}) yields
    \begin{align}
        J_0(x_*, \mcal{V}_*) + \frac{\gamma}{2} \Vert \xi \Vert_{\tilde{\mcal{V}}_*}
        &\leq J_0(x_*, \mcal{V}_*) + \gamma  R_{1,2}( \mcal{V}_* ) \\
        &\leq J_0(x_*, \tilde{\mcal{V}}_*) + \gamma \underbrace{ R_{1,2}( \tilde{\mcal{V}}_* )}_{0} \\
        &\leq J_0(x_*, \mcal{V}_*) + L \Vert \xi \Vert_{\tilde{\mcal{V}}_*},
    \end{align}
    which implies that $\gamma / 2 \leq L$.
    Since this contradicts $\gamma > 2 L = \Gamma_1$, we must have $\xi = 0$, which completes the proof of Theorem~\ref{thm:solutions_of_12_regularized_minimization}.
\end{proof}

\end{document}